
\input amstex.tex
\documentstyle{amsppt}
\magnification=1200
\baselineskip 16pt plus 2pt
\parskip 2pt
\NoBlackBoxes

\def\Vol{\text {Vol}}
\topmatter

\title Torsions for manifolds with boundary and glueing formulas \endtitle

\thanks  All three authors were supported by NSF. The first two
authors wish to thank the Erwin-Schr\"odinger-Institute in Vienna
for hospitality and support during the
summer of 1995 when part of this work was done.
\endthanks

\author D. Burghelea (Ohio State University)\\
L. Friedlander (University of Arizona)\\
T. Kappeler (Ohio State University)\endauthor
\leftheadtext {Analytic and Reidemiester torsion ... part II}
\rightheadtext {D. Burghelea, L. Friedlander, T. Kappeler}

\abstract We extend the definition of analytic and Reidemeister torsion
from closed compact Riemannian manifolds to compact Riemannian manifolds
with boundary $(M, \partial M)$, given
a flat bundle $\Cal F$ of $\Cal A$-Hilbert
modules of finite type and a
decomposition of the boundary $\partial M =\partial_- M \cup \partial_+ M$
into disjoint components.  If the system $(
M,\partial_-M, \partial_+M, \Cal F)$ is of determinant class we compute
the quotient of the analytic and the Reidemeister torsion and prove gluing
formulas for both of them. In particular we answer positively
Conjecture 7.6 in [LL]
\endabstract
\toc
\widestnumber\subhead{3.1}
\head 0. Introduction. \endhead
\head 1. Linear homological algebra in the von Neumann sense. \endhead
\head 2. Torsions for compact manifolds with boundary. \endhead
\subhead 2.1. Reidemeister and analytic torsion.\endsubhead
\subhead 2.2. Determinant class.\endsubhead
\subhead 2.3. Witten deformation of the de Rham complex
for manifolds with boundary. \endsubhead
\subhead 2.4. Asymptotic expansions and a comparison theorem.\endsubhead
\head 3. Applications to torsion.\endhead
\subhead 3.1. Comparison of analytic and Reidemeister torsion. \endsubhead
\subhead 3.2. Glueing formulas (PartI). \endsubhead
\subhead 3.3. Glueing formulas (PartII). \endsubhead
\subhead 3.4. Comparison of the analytic torsion for different boundary
conditions. \endsubhead
\head 4. Appendix: On manifolds of determinant class. \endhead

\endtoc

\endtopmatter

\document
\vfill \eject

\proclaim{0. Introduction}\endproclaim

The purpose of this paper is to extend the analysis of $L_2$-analytic and
$L_2$-Reidemeister torsion to compact Riemannian manifolds with boundary
and, among other results, to provide a formula for the quotient of the
$L_2$-analytic torsion and the $L_2$-Reidemeister torsion.  We also establish
glueing formulas for both torsions. In particular we prove Conjecture
7.6 in [LL].

In order to formulate our
results more precisely, we introduce the following notations.  Let $M$ be a
compact smooth manifold, not necessarily connected, of dimension $d,$ and
with boundary $\partial M$.  Let $g$ be a Riemannian metric on $M$.
Throughout the paper, we will always assume that $g$ has a product structure
near the boundary (product-like), i.e. that there exist $\epsilon >0$ and a
diffeomorphism
$\Theta_{\partial M} : \partial M \times[0, \epsilon) \rightarrow M$ so that
$$\leqalignno{
\Theta_{\partial M} \restriction{_{\partial M \times \{ 0 \}}}
= id &&(P1)\cr}$$
$$\leqalignno{
\Theta^*_{\partial M} (g) = g_{\partial M} + dt^2 &&(P2)\cr}$$
where $g_{\partial M}$ denotes the restriction of $g$ to $\partial M$ and
$dt^2$ denotes the Euclidean metric on the half open interval $[0,
\epsilon)$.
Suppose $\partial M$
is a union of two disjoint components, not necessarily connected and
possibly empty. Let us specify an ordering for them, $\partial_-M$ and
$\partial_+M$.  Then $(M, \partial_- M, \partial_+M)$ will be referred to
as a bordism. A closed manifold $M$ can be regarded as a bordism with
$\partial_- M=\partial_+M=\emptyset.$ If it is not ambiguous, we will write
$M$ instead of $(M,
\partial_-M, \partial_+M)$ and denote by $-M$ the bordism obtained from
$M$ by interchanging $\partial_-M$ and $\partial_+M$.  Next, we introduce the
notion of generalized triangulation for a bordism $(M, \partial_-M,
\partial_+M)$.

\definition{Definition}  A pair $\tau = (h, g')$ consisting of a
$C^\infty$-function $h: M\rightarrow {\Bbb R}$ with range $[a, b] \subseteq
{\Bbb R}$, where a and b are elements in $\Bbb Z,$ and a Riemannian metric
$g'$ is said
to be a generalized triangulation for $M$ if the following properties hold:
$$\leqalignno{
&h \big( \Theta_{\partial M} (x, t) \big) = b \, -t \qquad \big( x \in
\partial_+M,
0 \le t < \epsilon\big) &\text{\bf(T1)}\cr
&h \big( \Theta_{\partial M} (x, t) \big) = a \, +t \qquad \big( x \in
\partial_-M,
0 \le t < \epsilon\big)\cr}$$
where $\Theta_{\partial M} : \partial M \times[0, \epsilon) \rightarrow M$
denotes the exponential map at $\partial M$, associated to $g'$ (cf. above).
$$\leqalignno{
h \big( Cr (h) \big) \subseteq {\Bbb Z},&&\text{\bf(T2)}\cr}$$
where $Cr(h)$ denotes the set of critical points of $h$.  (Notice that we do
not require that $h$ is self-indexing.  This will be convenient for the
glueing construction.)

\noindent{\bf (T3)} All critical points of $h$ are nondegenerate and, given
any
critical point $y$ of index $k,\, y \in Cr_k (h),$ there exists a
diffeomorphism
$\phi_y : D_{\epsilon_y} \rightarrow M,$ with $D_{\epsilon_y}$ being the disc
in
${\Bbb R}^d$ of radius $\epsilon_y > 0,$ centered at $0$, so that
$\phi_y(0) = y,\, \phi_y^\# (g')$ is the Euclidean metric on $D_{\epsilon_y}$
and

$$
h \big( \phi_y (x_1, \dots, x_d) \big) = h (y) - \frac 1 2
\sum_{j = 1}^k x_j^2 + \frac 1 2 \sum_{j = k+1}^d x_j^2.$$
If for each critical point of h, such a coordinate system exists,
 the metric $g'$ is called compatible
with $h$.

\noindent{\bf (T4)} The gradient of $h$ with respect to $g'$,
$grad_{g'} h$,
satisfies the Morse-Smale condition (cf. [BFKM]).
\enddefinition
A smooth simplicial triangulation $t$ of $M$ with the property that
$\partial_{\pm}M$ become subcomplexes induces generalized triangulations
on the bordisms $(M,\partial_-M,\partial_+M)$,
$(M,\partial_+M,\partial_-M)$, $(M,\emptyset,\partial M)$,
$(M,\partial M,\emptyset)$ and on the closed manifolds $\partial_{\pm}M.$
\newline Assume that the following data is given:

$(M,\partial_- M, \partial_+ M) \text { a bordism; }\tau = (h, g')
\text { a generalized triangulation;}$

$\Cal F = ( \Cal E, \nabla)$ a parallel flat bundle  where $\Cal E \overset p
\to
\rightarrow M$ is a bundle of $\Cal A$-Hilbert modules of finite type
and $\nabla$ is a flat
connection making the inner products $(\cdot, \cdot)$ of $\Cal E_y = p^{-1}
(y)$
parallel with respect to $\nabla$.
\footnote{Given a choice of base points $(x_j)$ in each of the connected
components $M_j$ of $M(1 \le j \le k_M)$ $\Cal F$ determines and is
determined
uniquely by the $(\Cal A ,\Gamma_j^{op})-$Hilbert modules $\Cal W_j :=
\Cal E_{x_j}$ where the action of $\Gamma_j := \pi_1 (M_j, x_j)$ on $\Cal
W_j$
is given by the parallel transport and where $\Gamma_j^{op}$ denotes the
group
$\Gamma_j$ with opposite multiplication (cf. [BFKM]).}

If $\tau=(h,g')$ is a generalized triangulation of
$(M, \partial_-M, \partial_+M)$ then $\tau_D=(-h,g')$ is a generalized
triangulation of $(M, \partial_+M, \partial_-M).$ It will be refered
to as the dual triangulation.

In section 2 we define, similar as in [BFKM], analytic Laplacians $\Delta_q$,
acting on $q$-forms with coefficients in $\Cal E$ with relative (or
Dirichlet)
boundary conditions on $\partial_-M$ and absolute (or Neumann) boundary
conditions on $\partial_+M$, associated to  $(M, \partial_-M,
\partial_+M)$, g and $\Cal E$.  Further, we define combinatorial Laplacians
$\Delta_q^{comb}$, associated to $(M,\partial_- M),\  \tau$
and $\Cal F$.  We say
that the triple $\{(M, \partial_- M, \partial_+M), g, \Cal F\}$ is of $a$-
determinant class if all the analytic Laplacians $\Delta_q$ are of
determinant
class, i.e. $\log det'_N \Delta_q \in {\Bbb R}.$
For $\{(M,\partial_- M,\partial_+M),g,\Cal F\}$ of a-determinant class
 the analytic
torsion $T_{an} = T_{an} (M, \partial_- M, g, \Cal F) > 0$
is defined by
$$
\log T_{an} := \frac 1 2 \sum_{q=0}^d (-1)^{q+1} q \log det'_N \Delta_q.$$
Similarly, we say that the triple
$\{ (M,\partial_- M, \partial_+M ), \tau, \Cal F\}$
is of $c$-determinant
class if all the combinatorial Laplacians $\Delta_q^{comb}$ are of
determinant class,
\newline i.e. $\log det'_N \Delta_q^{comb} \in {\Bbb R}.$
For $\{(M,\partial_-M,\partial_+M),\tau,\Cal F\}$ of c-determinant
class
 the combinatorial torsion $T_{comb} = T_{comb} (M, \partial_-M,
\tau, \Cal F) >0$ is defined by
$$
\log T_{comb} := \frac 1 2  \sum_{q=0}^d (-1)^{q+1} q \log det'_N
\Delta_q^{comb}.$$
The Reidemeister torsion $T_{Re} = T_{Re} (M, \partial_-M, g, \tau,
\Cal F)$ is then defined by

$$
\log T_{Re} = \log T_{comb} + \log T_{met}$$
where $T_{met}$ is the metric part of the torsion, defined as in [BFKM] (cf.
section 2). $T_{met}$ is equal to 1 if $( M, \partial_-M)$ is $
\Cal F$-
acyclic.
It turns out that $T_{Re}$ does not depend on the triangulation
$\tau.$ The proof of this fact is not difficult and is the same as in
the case of finite dimensional unitary representations. This result is also
contained in Theorem 3.1. Nevertheless, we will continue to
write $T_{Re}(M,\partial_-M,
g,\tau,\Cal F)$ in order to emphasize that the Reidemeister torsion
is calculated with the triangulation $\tau.$

One can show that the notions of $a$-determinant class and $c$-determinant
class are equivalent (cf. [BFKM]) which allows us to introduce the notion of
a pair $\{ ( M, \partial_-M, \partial_+M), \Cal F \}$ being of
determinant class.
Given a bordism $(M, \partial_-M,\partial_+M)$
we say that the system $(M,\partial_-M,\partial_+M, \Cal F)$ is of
determinant class if the three pairs
$\{(M,\partial_-M,\partial_+M),\Cal F\}$,
$\{(M,\emptyset ,\partial M), \Cal F\}$ and $\{\partial M ,
\Cal F\restriction_{\partial M}\}$
are of determinant class. As we will see in Proposition 2.5 (4) this
is equivalent with  $\{(M,\emptyset ,\partial M), \Cal F\}$ and
$\{\partial M ,
\Cal F\restriction_{\partial M}\}$ being of determinant class.
If $\Cal F$ is a parallel flat bundle over a compact manifold
$M$ induced from a $\Gamma-$principal covering where
$\Gamma$ is a residually finite group, one can derive from the
work of L\"uck (cf.[L\"u3]),
that for any bordism $(M,\partial_-M,\partial_+M)$ the system
$( M, \partial_-M, \partial_+M, \Cal F)$
is of determinant class (cf. Theorem A in Appendix). In particular,
one concludes that if $M$ is connected and has residually finite fundamental
group $\Gamma := \pi_1(M),$
the system $(M,\partial_-M,\partial_+M,\Cal F)$ is of determinant class
if $\Cal F$ is the Hilbert bundle induced by the
$(\Cal N(\Gamma),\Gamma^{op})-$Hilbert module
$\ell_2(\Gamma)$ (cf Appendix).

Given a pair $\{( M, \partial_-M, \partial_+M), \Cal F \}$ of determinant
class, denote by
\newline $\Cal R = \Cal R (M,\partial_-M, g, \tau, \Cal F)$ the
quotient of analytic and Reidemeister torsion,
$$
\log \Cal R = \log T_{an} (M, \partial_-M, g, \Cal F) -
\log T_{Re} (M, \partial_-M, g, \tau, \Cal F).$$
Extending our earlier results [BFK], [BFKM] and using similar
techniques, we prove
that $\log \Cal R$ depends only on the data on the boundary:

\proclaim{Theorem 2.10} Assume that for $j=1,2, (M_j,\partial_-M_j,
\partial_+M_j)$ is a bordism equipped with a Riemannian
metric $g_j$, a generalized triangulation $\tau_j$ and that
$\Cal F_j = ( \Cal E_j, \nabla_j)$ is a parallel flat bundle of
$\Cal A$-Hilbert modules of
finite type on $M_j$. Assume also that the systems
$(M_j,\partial_-M_j, \partial_+M_j,\Cal F_j)$ are of determinant class and
$$
(\partial_\pm M_1, g_1 \restriction_{\partial_\pm M_1}, \Cal F_1
  \restriction_{\partial_\pm M_1})
=(\partial_\pm M_2, g_2 \restriction_{\partial_\pm M_2}, \Cal F_2
  \restriction_{\partial_\pm M_2}). $$
Then
$$ \Cal R ( M_1, \partial_-M_1, g_1, \tau_1, \Cal F_1) =
\Cal R ( M_2, \partial_-M_2, g_2, \tau_2, \Cal F_2).$$
\endproclaim

This theorem allows us to extend the comparison theorem (cf. [BFKM] and
[BFK]) of analytic and Reidemeister torsions from closed manifolds to
bordisms.

\proclaim{Theorem 3.1}  Assume that the system
$(M, \partial_- M, \partial_+ M,
\Cal F)$
is of determinant class.
Then $\log \Cal R$ is given by
$$
\log \Cal R(M, \partial_-M, g, \tau, \Cal F) = \frac 1 4
\chi (\partial M; \Cal F)
\log 2$$
where $\chi (\partial M; \Cal F)$ is the Euler
characteristic of $\partial M$ with coefficients in $\Cal F$ and
is equal to
$\chi (\partial M; \Cal F)=
\sum_k\chi (\partial_k M)\cdot dim_N\Cal E
\restriction_{\partial_k M}$ where $\partial_k M$ are the connected
components of the boundary $\partial M$, $\chi(\partial_k M)$ is the
standard Euler characteristic of $\partial_k M$and
$dim_N\Cal E\restriction_{\partial_k M}$
is the von Neumann dimension of the fiber of $\Cal E$ above
$\partial_k M $.
\endproclaim

\remark{Remark } In the case $\Cal A$ is ${\Bbb R}$ or ${\Bbb C}$,
Theorem 3.1 is
due to W. L\"uck [L\"u1] and, independently, Vishik [Vi] (cf. also [Ch]).

The above result shows that $\Cal R$ does not depend on
the partition $\partial M_-, \partial M_+$ of the boundary
$\partial M.$  However we continue to use the
notation $\Cal R(M, \partial_-M, g, \tau, \Cal F)$ as the independence of
$\partial_-M$ will be proven only towards the end of the paper.
\endremark

Next we present a gluing formula for the analytic torsion.  For $j = 1, 2$,
let $M_j = (M_j, \partial_- M_j, \partial_+M_j)$, $g_j, \tau_j =
(h_j, g_j')$ and $\Cal F_j = (\Cal E_j, \triangledown_j)$ be as
above and suppose that there exist an isometry

$$
\omega : (\partial_+ M_1, g_1 \restriction_{\partial_+ M_1}) \rightarrow
(\partial_- M_2, g_2 \restriction_{\partial_- M_2})$$
and a connection preserving bundle isometry $\Phi$ above $\omega$ which
makes
the following diagram commutative

$$
\matrix
\Cal E_1 \restriction_{\partial_+ M_1} & \overset \Phi \to \longrightarrow
&\Cal E_2 \restriction_{\partial_- M_2}\\
\downarrow && \downarrow\\
\partial_+ M_1 & \overset \omega \to \longrightarrow & \partial_- M_2.
\endmatrix$$
Then one can form the bordism $(M = M_1 \cup_\omega M_2, \partial_- M_1,
\partial_+ M_2)$ by glueing $\partial_+ M_1$ to $\partial_- M_2$ by $\omega$
and
the parallel flat bundle $\Cal F$ by glueing $\Cal F_1$ and $\Cal F_2$ by
$(\omega, \Phi)$.  The metrics $g_1, g_1'$ and $g_2, g_2'$ determine
Riemannian
metrics $g$ and $g'$ on $M$, and the functions $h_1$ and $h_2$ determine the
$C^{\infty}$-function $h : M \rightarrow {\Bbb R}$ given by

$$
h(x):=\cases h_1(x)  &(x \in M_1)\\
b_1 - a_2 + h_2(x) & (x \in M_2)\endcases$$
where, for $j=1,2,$ $h_j(M_j) = [a_j,b_j].$
\definition{Definition} Generalized triangulations $\tau_1$ and $\tau_2$
as above are said to be
compatible if $\tau = (h, g')$ is a generalized triangulation
of $(M, \partial_-M_1 ,\partial_+ M_2).$
\enddefinition

Notice that by
an arbitrary small perturbation,
localized in a given neighborhood of $\partial_+ M_1$ one can modify
the metric $g_1'$ to $\overset \sim
\to g_1'$ so that
the triangulation $\overset \sim \to \tau_1:=(h_1, \overset \sim\to g_1')$
is compatible with $\tau_2$.

Finally the manifolds $M, M_1,M_2$ and the metric g induce a
cohomology sequence $\Cal H_{an}(g),$ which is a long weakly exact
sequence of $\Cal A-$Hilbert modules, and hence a cochain complex.
Different metrics induce isomorphic, but not necessarily isometric
complexes. If $\Cal H_{an}(g)$ is of determinant class we
denote by $T(\Cal H_{an}(g))$ its torsion. Similarly, the manifolds
$M,M_1,M_2$ and the generalized triangulation $\tau$ induce a
cohomology sequence $\Cal H_{comb}(\tau),$
which is a long weakly exact sequence of
$\Cal A-$Hilbert modules and thus, again, a cochain complex.
Integration theory (cf section 2) provides an isomorphism
(which typically is not an isometry) between the two sequences
$\Cal H_{an}(g)$ and $\Cal H_{comb}(\tau).$ Both sequences
are of the following form:

$$\eqalignno{
0&\to  H^0 (M_2, \partial_- M_2, \Cal F_2)  \to
H^0 (M, \partial_- M, \Cal F) \to
H^0 (M_1, \partial_- M_1, \Cal F_1) \to \cr
&\to  H^1 (M_2, \partial_- M_2, \Cal F_2 ) \to \dots
\to H^k (M_2,\partial_- M_2, \Cal F_2) \to
H^k (M, \partial_- M, \Cal F)
\to \cr &\to  H^k (M_1, \partial_- M_1, \Cal F_1)
\to H^{k +1} (M_2, \partial_- M_2, \Cal F_2) \to  \dots \cr}$$

\proclaim{Theorem 3.2}
Assume that for $i=1,2$ the system $(M_i,\partial_-M_i,\partial_+M_i,\Cal
F_i)$
is of determinant
class.  Then the following statements hold:
\newline (i) The system $(M,\partial_-M,\partial_+M, \Cal F)$ and the
complexes $\Cal H_{comb}, \Cal H_{an} $ are of
determinant class;
$$\leqalignno{
\log T_{Re} (M, \partial_-M, g, \tau, \Cal F) =
\overset 2 \to
{\underset {j = 1} \to \sum} &\log T_{Re} (M_j, \partial_-M_j,
g_j, \tau_j,\Cal F_j)
&{(ii)}\cr
&+ \log T (\Cal H_{an}).\cr}$$
$$
\leqalignno{
\log T_{an} (M, \partial_-M, g, \Cal F)
&= \sum^2_{j = 1} \log T_{an} (M_j, \partial_-M_j, g_j, \Cal F_j)
+ \log T(\Cal H_{an}) &{(iii)}\cr
&-\frac {\chi ( \partial_+ M_1; \Cal F_1)} {2}
\log 2.\cr}$$
\endproclaim

\remark{Remark} In the case where $\Cal A = {\Bbb R}$ this result is due to
Vishik
[Vi3]. Vishik's proof is, however, very
different from ours.
\endremark

A sightly more general form of this result is
contained in Theorem 3.2$'$ where only
some components of $\partial_+M_1$ and $\partial_-M_2$ are glued together.
The last result is Theorem 3.3 which compares the analytic torsions of
bordisms with the same underlying compact manifold $M$.

The paper is organized as follows:
\newline In section 1 we prove a number of auxiliary results about cochain
complexes of
determinant class and a generalization of Milnor's result concerning
the torsions of a short exact sequence of finite
dimensional complexes to infinite dimension (Proposition 1.13,
Theorem 1.14). This section is quite elementary
but it is included for the sake of completeness.
In section 2 we prove Theorem 2.9 on which all subsequent results depend on.
Theorems 3.1, 3.2, 3.2', 3.3 are proven in section 3.
In the appendix, we prove that for any bordism
$(M, \partial_- M, \partial_+M)$
with $M$ a compact connected manifold, $\tilde{M} \to M$ a principal
covering with residually finite group
and $\Cal F$ the parallel flat bundle on $M,$ induced
by this covering, the system
$(M, \partial_- M, \partial_+M ,\Cal F)$
is of determinant
class.  This result is implicit in L\"uck [L\"u3].  Unfortunately, in [L\"u3]
there are a number of misleading misprints and the definition of the
$L^2$-determinant is incorrect.  For the convenience of the reader we present
an
outline of the proof of this result.

By the same arguments as given in [BFKM, Proposition 5.11]
it suffices to prove
the above Theorems in the case where the fiber $\Cal W$
of the bundle $\Cal E \to M$ is a free $\Cal A$-module.
For the rest of the paper the bundle $\Cal E$ has as a fiber a free
$\Cal A$-module.

\vfill \eject

\subhead {Section 1 Linear homological algebra in the von Neumann sense}
\endsubhead

In this section we follow, if not stated otherwise, the notations and use the
definitions of Section 1 in [BFKM]. In particular, contrary to the
standard notations for cochain complexes, the indice for the
degree is denoted as a subscript and not as a superscript.

Throughout this section, let $\Cal A$ be a von Neumann algebra of finite
type.
By
$\Cal W, \Cal W_1, \Cal W_2, \Cal U, \dots$ we denote arbitrary $\Cal
A$-Hilbert
modules of finite type.

Recall that $f: \Cal W_1 \rightarrow \Cal W_2$ is said to be a weak
isomorphism
if $f$ is 1-1 and its range, $Im(f)$, is dense in $\Cal W_2$.  It admits a
polar decomposition $f = f_{iso} f_{wiso}$ where $f_{iso} := f (f^*f)^{-1/2}
:
\Cal W_1 \rightarrow \Cal W_2$ is an isometry and $f_{wiso} := (f^* f)^{1/2}
:
\Cal W_1 \rightarrow \Cal W_1$ is a weak isomorphism which is selfadjoint and
positive.  Therefore $f_{wiso}^2 +
\epsilon (\epsilon >0)$ is an isomorphism and one can define $\Vol(f)$ by

$$\leqalignno{
\log \Vol (f) := \lim_{\epsilon \rightarrow 0} \frac 1 2 \log \det
(f_{wiso}^2
+ \epsilon)&&(1.1)\cr}$$
where the determinant is the one given by Fuglede-Kadison [FK] and where
$\log$
denotes the branch of the logarithm with $\log \, 1 = 0$.  Notice that $\frac
1 2
\log \det (f_{wiso}^2) := \log \Vol (f)$ is in ${\Bbb R} \cup \{- \infty \}$
and

$$\leqalignno{
\log \Vol (f) = \lim_{\epsilon \rightarrow 0} \frac 1 2 \log \det (f^*f+
\epsilon).&&(1.2)\cr}$$
The weak isomorphism $f$ is said to be of \underbar{determinant class} iff

$$\leqalignno{
- \infty < \log \Vol (f).&&(1.3)\cr}$$
Using the spectral distribution function of $f,$
$$F_f(\lambda):=\sup \{ \,{\roman {dim}}_N \Cal L \,; \Cal L
 \, {\roman {subspace\, of }} \,\Cal W_1,
\|f(x))\| \leq \sqrt \lambda \|x\| (x \in \Cal L)\} ,$$
one can show
that
$$\leqalignno{
\log \Vol (f) = \int_{0+}^{\infty} \log(\lambda)dF_f (\lambda)
&&(1.2')\cr}$$
and that (1.3) is equivalent to
$$\leqalignno{
- \infty < \int_{0+}^{1} \log(\lambda)dF_f (\lambda).
&&(1.3')\cr}$$
In case $(1.3')$ holds, one verifies that
$$\int_0^1 \frac {F_f (\lambda) - F_f (0)} {\lambda} d\lambda =
- \int_{0^{+}}^1 \log (\lambda) dF_f(\lambda ),
$$
using that the boundary term obtained, when integrating by parts,
vanishes.
It will be convenient to extend the concept of determinant class to
an arbitrary morphism $f \in \Cal L_{\Cal A} (\Cal W_1, \Cal W_2).$
For this purpose we factor $f$ as $f=j\cdot f'\cdot p$ with
$p:\Cal W_1\to \Cal W_1'=\Cal W_1 /{Null f},$ the canonical
projection and
$j:\Cal W'_2 := \overline{Range f} \to \Cal W_2$ the inclusion.
 Then $f': \Cal W'_1 \to \Cal W'_2$ is a weak
isomorphism.
The morphism $f$ is said to be of determinant class if
$f'$
is of determinant class.

\proclaim{Proposition 1.1}

\noindent (A) Assume that $f \in \Cal L_{\Cal A} (\Cal W_1, \Cal W_2)$ and $g
\in
\Cal L_{\Cal A} (\Cal W_2, \Cal W_3)$ are
both weak isomorphisms.  Then $g \cdot f$ is a weak isomorphism.  Moreover

$$\leqalignno{
\log \Vol (g \cdot f) = \log \Vol (g) + \log \Vol (f).&&(1.4)\cr}$$
In particular $g \cdot f$ is of determinant class iff both $g$ and $f$ are
of determinant class.

\noindent (B) Assume that $f \in \Cal L_{\Cal A} (\Cal W_1, \Cal W_1')$ and
$g
\in \Cal L_{\Cal A} (\Cal W_2, \Cal W_2')$ are both weak isomorphisms and
$h \in \Cal L_{\Cal A} (\Cal W_2, \Cal W_2')$.  Then $f' = \pmatrix
f & h\\ 0 & g \endpmatrix : \Cal W_1 \oplus \Cal W_2 \rightarrow \Cal W_1'
\oplus \Cal W_2'$ is a weak isomorphism and
$$\leqalignno{
\log \Vol (f') = \log \Vol (f) + \log \Vol (g).&&(1.5)\cr}$$

\noindent (C) If $f_i:\Cal W \to \Cal V_i, i=1,2$ are two morphisms and
both are of determinant class then $f_1\oplus f_2:\Cal W
\to \Cal V_1\oplus\Cal V_2$ is of determinant class.
\endproclaim

\remark{Remark}: (a)  In the case where $\Cal W_1 = \Cal W_2$, statement (A)
is
verified in [FK].

\noindent (b)  In the case $f, g$ are isomorphisms and not only weak
isomorphisms,
statement (B) is proved in [BFKM, Proposition 1.9].
\endremark

\demo{Proof}  (A)  Consider the polar decomposition $f = f_{iso} f_{wiso}$
and
$g= g_{iso} g_{wiso}$.  Then $gf$ can be written as $gf = \overset \sim \to g
\overset \sim \to f$ where

$$\leqalignno{
\overset \sim \to g := g_{iso} f_{iso} : \Cal W_1 \rightarrow \Cal W_3; \,
\overset \sim \to f := f_{iso}^{-1} g_{wiso} f : \Cal W_1 \rightarrow
\Cal W_1. &&(1.6)\cr}$$
Note that $\overset \sim \to g$ is an isometry.  Therefore $\overset \sim \to
g^*
\overset \sim \to g= Id_{\Cal W_1}$ and $(\overset \sim \to g \overset \sim
\to f)^*
(\overset \sim \to g \overset \sim \to f) = \overset \sim \to f^*\overset
\sim
\to
f$.  From definition (1.1) it then follows that

$$\leqalignno{
\log \Vol (gf) = \log \Vol (\overset \sim \to f).&&(1.7)\cr}$$
The map $\overset \sim \to f$ can be decomposed, $\overset \sim \to f=
(f_{iso}^{-1} g_{wiso} f_{iso})f_{wiso}$ and we conclude from (1.4),
which is applied in the case
$\Cal W_1 = \Cal W_2,$ that

$$\leqalignno{
\log \Vol (\overset \sim \to f) = \log \Vol (f_{iso}^{-1} g_{wiso} f_{iso} )
+
\log \Vol (f_{wiso}).&&(1.8)\cr}$$
Taking into account (1.2) one concludes that

$$\leqalignno{
\log \Vol (f_{iso}^{-1} g_{wiso} f_{iso}) = \log \Vol(g_{wiso}).&&(1.9)\cr}$$
{}From definition (1.1) we know that $\log \Vol (f) = \log \Vol (f_{wiso})$
and similarly for $g$.  Combining (1.7) - (1.9) leads to

$$
\log \Vol (gf) = \log \Vol (g) + \log \Vol (f).$$

(B)  To verify (B) decompose $f'$ as follows

$$\leqalignno{
f'= \pmatrix Id_{\Cal W_i'}&0\\ 0&g\endpmatrix \pmatrix Id_{\Cal W_i'} & h \\
0 & Id_{\Cal W_2}\endpmatrix \pmatrix f & 0\\ 0 & Id_{\Cal W_2}\endpmatrix
&&(1.10)\cr}$$
and apply (A) together with Proposition 1.9 in [BFKM].

 (C)   Denote by $F_1(\lambda):= F_{f_1}(\lambda),
F_2(\lambda):= F_{f_2}(\lambda)$ and
$F(\lambda):= F_{f_1\oplus f_2}(\lambda).$ It is straightforward
to verify that $F(\lambda)-F(0)\le \sup \{ F_1(\lambda)-F_1(0),
F_2(\lambda)-F_2(0)\}.$  Integration by parts  in $(1.3')$ implies the
result.
\hfill $\square$
\enddemo

Consider a finite $\Cal A$-cochain complex $(\Cal C, d)$ of $\Cal A$-Hilbert
modules of finite type

$$\leqalignno{
0 \rightarrow \Cal C_0 \overset {d_0} \to \rightarrow \Cal C_1 \overset {d_1}
\to \rightarrow \dots \rightarrow 0.&&(1.11)\cr}$$
Define the Hodge decomposition $\Cal C_i, \Cal C_i = \Cal H_i \oplus
\Cal C_i^{ +} \oplus \Cal C_i^{-}$, where

$$\leqalignno{
\Cal C_i^{+} := \overline{d_{i-1} (\Cal C_{i-1})} \, ;\quad \Cal C_i^{-}:=
\overline {d_i^* (\Cal C_{i+1})} \, ;\quad \Cal H_i := Null (d_i) \cap Null
(d_{i-1}^*).&&(1.12)\cr}$$
Recall that $\Cal H_i$ is isomorphic to $Null d_i/\overline {Range(
d_{i-1})}=H^i(\Cal C).$
 With respect to this decomposition

$$\leqalignno{
d_i = \pmatrix 0 & 0 & 0 \\
0 & 0 & \underline d_i\\
0 & 0 & 0 \endpmatrix ; \, \underline d_i : \Cal C_i^{-} \rightarrow
\Cal C_{i+1}^+ .&&(1.13)\cr}$$
Note that $\underline d_i$ is a weak isomorphism.

\definition{Definition} The $\Cal A$-cochain complex $\Cal C$ is said to be
of determinant class if, for any $i$, $d_i$ (or, equivalently,
$\underline d_i$) is of determinant class.

If $\Cal C$ is of determinant class, its torsion $T(\Cal C)$ is a positive
real number given by

$$\leqalignno{
\log T(\Cal C) = \sum_j (-1)^j \log \Vol (\underline d_j).&&(1.14)\cr}$$
For $k = 1, 2$, let $\Cal C^k$ be $\Cal A_k$-cochain complexes.  Denote by
$\chi (\Cal C^k)$ the Euler-Poincar\'e characteristic of $\Cal C^k$.
Introduce
the $\Cal A_1 \hat \otimes \Cal A_2$-cochain complex $\Cal C$
given by $\Cal C_j = \oplus_k (\Cal C^1_k \hat \otimes \Cal C^2_{j-k})$.
Here the tensor product is taken in the category of Hilbert spaces.
\enddefinition

In [BFKM, Corollary 1.22] we have proved the following

\proclaim{Proposition 1.2}  Suppose $\Cal C^1$ and $\Cal C^2$ are both of
determinant class.  Then $\Cal C$ is of determinant class and

$$\leqalignno{
\log T(\Cal C) = \chi (\Cal C^2) \log T (\Cal C^1) + \chi ( \Cal C^1) \log
T (\Cal C^2).&&(1.15)\cr}$$
\endproclaim

Suppose $f : \Cal C^1 \rightarrow \Cal C^2$ is a morphism of chain complexes.
By
straightforward inspection (cf. [BFKM, section 1]) one sees that, with
respect
to the Hodge decompositions of $\Cal C^1_i = \Cal H^1_i
\oplus \Cal C_i^{1, +}
\oplus \Cal C_i^{1, -}$ and $\Cal C_i^2 = \Cal H^2_i
\oplus \Cal C_i^{2, +}
\oplus \Cal C_i^{2, -}$, the morphisms $f_i : \Cal C^1_i \rightarrow
\Cal C^2_i$
are of the form

$$\leqalignno{
f_i = \pmatrix f_{i, 11} & 0 & f_{i, 13}\\
f_{i, 21} & f_{i, 22} & f_{i, 23}\\
0 & 0 & f_{i, 33}\endpmatrix. &&(1.16)\cr}$$
Moreover

$$\leqalignno{
\underline d_i^2 f_{i, 33} = f_{i+1, 22} \underline d_i^1.&&(1.17)\cr}$$
Further it is straightforward to conclude from Proposition 1.1 that for any
$i
\ge 0$
the following statements hold:

$$\leqalignno{
f_i &\text { is a [weak] isomorphism iff } f_{i, 11}, f_{i, 22} \text
{ and }&(1.18)\cr
& f_{i, 33} \text { are all [weak] isomorphisms}.&\cr}$$

$$\leqalignno{
\log \Vol f_i = \log \Vol f_{i, 11} + \log \Vol f_{i, 22} + \log \Vol
f_{i, 33}.&&(1.19)\cr}$$

$$\leqalignno{
&\quad\text {Assume that } f_i \text{ is a weak isomorphism.  Then } f_i
\text { is of determinant class }&(1.20)\cr
&\qquad \text {iff }f_{i,11}, f_{i, 22}, f_{i, 33} \text { are all of
determinant class}.\cr}$$

In the sequel, we will occasionally write $H(f_i)$ for $f_{i, 11} :
\Cal H^1_i
\rightarrow \Cal H^2_i$.

\proclaim{Proposition 1.3} Let $f: \Cal C^1 \rightarrow \Cal C^2$ be
a morphism of cochain complexes.

(A) If, for any i, $f_i$ is a weak isomorphism and of determinant
class,
 then the following statements hold:
\newline (i) For any i, $H(f_i)$ is of determinant class;
\newline (ii) $\Cal C^1$ is of determinant class iff $\Cal C^2$ is;
\newline (iii) if $\Cal C^1$ is of determinant class (and thus, by
(ii), $\Cal C^2$ is of determinant class as well), then
$$\leqalignno{
\log T(\Cal C^2) & = \log T(\Cal C^1) - \sum_i (-1)^i \log \Vol (f_i)
&(1.21)\cr
&\quad + \sum_i (-1)^i \log \Vol \big(H(f_i)\big).\cr}$$
(B)  If $f$ is a homotopy equivalence (cf. [BFKM], definition 1.14
) then $\Cal C^1$ is of determinant class iff $\Cal C^2$is.
\endproclaim

\demo{Proof}  (A) (i) is obtained from (1.20) and (ii) follows from
(1.20), (1.17) and Proposition 1.1.
To check (iii), apply Proposition 1.1 to $\underline d_i^2 f_{i, 33} =
f_{i+1,
22}
\underline d_i^1$ to conclude that

$$\leqalignno{
(-1)^i &\log \Vol \,\underline d_2^i &(1.22)\cr
&\quad = (-1)^i \log \Vol \underline d_1^i + (-1)^{i+1}
\log \Vol f_{i, 33} + (-1)^i \log \Vol f_{i+1, 22}\cr}$$
which, after summing up and using (1.14), leads to

$$\leqalignno{
\log T(\Cal C^2) & - \log T(\Cal C^1) &(1.23)\cr
&= - \sum_i (-1)^i \log \Vol (f_{i, 22}) - \sum_i (-1)^i \log \Vol (f_{i,
33})\cr
&= - \sum_i (-1)^i \log \Vol (f_i) + \sum (-1)^i \log \Vol H(f_i),\cr}$$
where for the last equality, we used (1.19).

(B) This result is due to Gromov-Shubin [GS] ( cf. also Proposition 1.18 in
[BFKM]).
\enddemo\hfill $\square$

Let $\Cal C$ be a cochain complex of the form

$$\leqalignno{
0 \rightarrow \Cal C_0 \overset {d_0} \to \rightarrow \Cal C_1 \overset {d_1}
\to \rightarrow \Cal C_2 \rightarrow 0.&&(1.24)\cr}$$
Such a complex is called a three stage complex.

\definition{Definition}  A three stage complex $\Cal C$ is said to be
a weakly exact sequence if $d_0$ \
is injective, $ \overline {\text {Range }(d_0)} = Null (d_1)$ and $ \overline
{\text {Range }(d_1)} = \Cal C^2$.
\enddefinition

Note that in this situation $\underline d_0 : \Cal C_0 \rightarrow
\overline {\text {Range }(d_0)}$ and $\underline d_1 : Range (d_0)^\perp
\rightarrow \Cal C_2$ are weak isomorphisms where $\overline
{\text {Range }(d_0)}^\perp$ denotes the orthogonal complement of
$Range (d_0)$ in $\Cal C_1$.

In the case where $\Cal C$ is of determinant class, its torsion is given by

$$\leqalignno{
\log T (\Cal C) = \log  \Vol (\underline d_0) - \log \Vol (\underline d_1).
&&(1.25)\cr}$$

\proclaim{Lemma 1.4}  Suppose that the three stage cochain complex
$\Cal C$ is a weakly exact sequence of the form

$$\leqalignno{
0 \rightarrow \Cal C_0 \overset {d_0} \to \rightarrow \Cal C_1^{ +} \oplus
\Cal C_1^{-} \overset {d_1} \to \rightarrow \Cal C_2 \rightarrow 0
&&(1.26)\cr}$$
where $d_0 = {f\choose {f_1}}$ and $d_1 = (g_1 \, g)$.  Further assume that
$d_1$
is onto.  Then the following statements hold:

(A) If $f : \Cal C_0 \rightarrow \Cal C_1^{ +}$ and
$g : \Cal C_1^{-} \rightarrow
\Cal C_2$ are weak isomorphisms of determinant class, then the complex is of
determinant class and

$$\leqalignno{
\log T(\Cal C) = \log \Vol (f) - \log \Vol (g).&&(1.27)\cr}$$

(B)  Assume that (1.26) is of determinant class.  If
$f : \Cal C_0 \rightarrow
\Cal C_1^{+}$ [resp. $g: \Cal C_1^{-} \rightarrow \Cal C_2$]
is of determinant
class, so is $g: \Cal C_1^{-} \rightarrow \Cal C_2$ [resp. $f :
\Cal C_0 \rightarrow \Cal C_1^{ +}]$ and formula (1.27) holds.
\endproclaim

\demo{Proof}  First note that, due to the assumption that $d_1$ is onto,
$\underline
d_1: \big( Range (d_0)\big)^\perp \rightarrow \Cal C_2$ is an isomorphism
and therefore of determinant class.  Moreover, $\log \Vol (\underline d_0)
\in
{\Bbb R} \cup \{ - \infty\}.$  Therefore, we can define $\log T(\Cal C)= \log
\Vol (\underline d_0) - \log \Vol (\underline d_1) \in {\Bbb R} \cup
\{ -\infty \}$.  Denote the inverse $(\underline d_1)^{-1}: \Cal C_2
\rightarrow
Range (d_0)^\perp \subseteq \Cal C_1^{+} \oplus \Cal C_1^{-}$ by
$(\underline d_1)^{-1} = {{h_1}\choose h}$ and conclude that

$$\leqalignno{
Id_{\Cal C_2} = ( g_1 \, \, g) {{h_1}\choose h} = g_1 h_1 + gh.&&(1.28)\cr}$$
\enddemo

Consider the complex $\Cal C'$ given by

$$
0 \rightarrow \Cal C_0 \overset {{Id_{\Cal C_0}}\choose  0\quad} \to
\longrightarrow \Cal C_0 \oplus \Cal C_2
\overset {(0\,\, Id_{\Cal C_2})}
\to \longrightarrow \Cal C_2\rightarrow 0$$
which is of determinant class and satisfies

$$\leqalignno{
\log T(\Cal C') = 0. &&(1.29)\cr}$$

Further introduce the morphism  $f: \Cal C \rightarrow \Cal C'$ given by
the following diagram

$$
\matrix 0 &\rightarrow &\Cal C_0 &\rightarrow &\Cal C_0 \oplus \Cal C_2
&\rightarrow &\Cal C_2 &\rightarrow &0\\
&&\qquad\downarrow Id_{\Cal C_0}&&\quad\qquad\qquad\downarrow {\pmatrix f &
h_1\\
f_1 & h\endpmatrix}&&\qquad\downarrow Id_{\Cal C_2}\\
0 & \rightarrow & \Cal C_0 & \overset {f\choose {f_1}} \to \rightarrow&
\Cal C_1^{+} \oplus \Cal C_1^{-} &\overset {(g_1, g)} \to \rightarrow &
\Cal C_2 & \rightarrow & 0.\endmatrix$$

Using the same arguments as in the proof of Proposition 1.3 one shows that,
in
view
of (1.29),

$$\leqalignno{
\log T (\Cal C) = 0 - \log \Vol \pmatrix f & h_1\\ f_1& h\endpmatrix + 0.
&&(1.30)\cr}$$
In particular, one concludes that $\Cal C$ is of determinant class iff
$\pmatrix f & h_1\\ f_1 & h \endpmatrix$ is of determinant class.

To compute $\log \Vol \pmatrix f & h_1 \\ f_1 & h \endpmatrix$, note that the
weak exactness of $\Cal C$ implies that $g_1 f + gf_1 = 0$.  Therefore

$$\leqalignno{
\pmatrix Id_{\Cal C_1^{+}} & 0 \\ 0 & g \endpmatrix \pmatrix f & h_1\\ f_1 &
h
\endpmatrix &= \pmatrix f & h_1 \\ gf_1 & gh\endpmatrix = \pmatrix f & h_1\\
-g_1f& gh\endpmatrix  &(1.31)\cr
&=\pmatrix Id_{\Cal C_1^{+}} & h_1\\ -g_1 & gh\endpmatrix
\pmatrix f & 0 \\ 0 & Id_{\Cal C_2}\endpmatrix \cr}$$
and, by Proposition 1.1, one obtains (in ${\Bbb R} \cup \{ - \infty \} $)

$$
\log \Vol \pmatrix f & h_1\\ f_1 & h \endpmatrix + \log \Vol (g) = \log \Vol
(f)
+ \log \Vol \pmatrix Id_{\Cal C_1^+} & h_1\\ -g_1 & gh \endpmatrix.$$

According to (1.28), $gh = Id_{\Cal C_2} - g_1 h_1$ and therefore

$$\leqalignno{
\pmatrix Id_{\Cal C_1^+} & h_1 \\ -g_1 & gh \endpmatrix = \pmatrix Id_{\Cal
C_1^+} & 0\\
-g_1 & Id_{\Cal C_2} \endpmatrix \pmatrix Id_{\Cal C_1^+} & h_1 \\ 0 &
Id_{\Cal C_2} \endpmatrix . &&(1.32)\cr}$$

Formula (1.32) implies that $\pmatrix Id_{\Cal C_1^{+}} & h_1\\ -g_1 & gh
\endpmatrix $ is a weak isomorphism of determinant class and

$$\leqalignno{
\log \Vol \pmatrix Id_{\Cal C_1^{+}} & h_1 \\ -g_1 & gh \endpmatrix = 0 + 0.
&&(1.33)\cr}$$
Combining (1.30) - (1.33) statements (A) and (B) follow.\hfill $\square$

\proclaim{Lemma 1.5}  Let $\Cal C^1, \Cal C^2, \Cal C^3$ be three
stage complexes which are short
exact sequences, (and thus, in particular, of determinant class)

$$\eqalign{
\Cal C^1 &: 0 \rightarrow \Cal C_0^1 \overset {d_0^1}
\to \rightarrow \Cal C_1^1
\overset {d_1^1} \to \rightarrow \Cal C_2^1 \rightarrow 0\cr
\Cal C^2 &: 0 \rightarrow \Cal C_0^1 \oplus \Cal C^3_0 \overset {d_0^2} \to
\rightarrow
\Cal C_1^1 \oplus \Cal C^3_1 \overset {d_1^2}
\to \rightarrow \Cal C_2^1 \oplus
\Cal C_2^3 \rightarrow 0 \cr
\Cal C^3 &: 0 \rightarrow \Cal C_0^3 \overset {d_0^3}
\to \rightarrow \Cal C_1^3
\overset {d_1^3} \to \rightarrow \Cal C_2^3 \rightarrow 0\cr}$$
Further assume that the following diagram is commutative,

$$\leqalignno{
\qquad\matrix 0 & \rightarrow &\Cal C_0^1 & \rightarrow & \Cal C_1^1 &
\rightarrow&
\Cal C_2^1 & \rightarrow & 0 \\
&&&&&&\\
&& \qquad\downarrow {Id \choose 0} && \qquad\downarrow {Id \choose 0} &&
\qquad\downarrow {Id \choose 0} \\
&&&&&&\\
0 & \rightarrow &\Cal C_0^1 \oplus \Cal C_0^3 & \rightarrow &
\Cal C_1^1 \oplus \Cal C_1^3  & \rightarrow&
\Cal C_2^1  \oplus \Cal C_2^3 & \rightarrow & 0\\
&&&&&&\\
&& \qquad\quad\downarrow (0\, Id) && \qquad\quad\downarrow (0\, Id) &&
\qquad\quad\downarrow (0\, Id) \\
&&&&&&\\
0 & \rightarrow &\Cal C_0^3 & \rightarrow & \Cal C_1^3 & \rightarrow&
\Cal C_2^3 & \rightarrow & 0.\endmatrix &&(1.34)\cr}$$
Then

$$\leqalignno{
\log T (\Cal C^2) = \log T (\Cal C^1) + \log T ( \Cal C^3).&&(1.35)\cr}$$
\endproclaim

\demo{Proof}We first note that in the Hodge decomposition of
$\Cal C_1^1=\Cal C_1^{1,+}\oplus \Cal C_1^{1,-},$ and
$\Cal C_1^3=\Cal C_1^{3,+}\oplus \Cal C_1^{3,-}$   the Hilbert module
$\Cal C_1^{j,+}$ is isometric to $\Cal C_0^j$ and $\Cal C_1^{j,-}$
is isometric to $\Cal C_2^j$,  $(j \in \{ 1, 3 \}).$
Therefore one can write
$
d_0^j = {{\underline d^j_0} \choose 0} ; d_1^j = (0\,
\underline d_1^j),$
with $\underline d^j_0:\Cal C^j_0\to \Cal C^{j,+}_1$
and $\underline d^j_1:\Cal C^{j,-}_1 \to \Cal C^j_2$

Taking the above mentioned isometries into account
it follows that $\Cal C^2$ is isometric to the decomposition

$$
0 \rightarrow \Cal C^1_0 \oplus \Cal C^3_0 \overset {d_0^2} \to \rightarrow
\Cal C^1_0 \oplus \Cal C^3_0 \oplus \Cal C^1_2 \oplus \Cal C^3_2 \overset
{d_1^2}
\to \rightarrow \Cal C^1_2 \oplus \Cal C^3_2 \rightarrow 0 $$ and
the commutativity of the diagram (1.34) implies that $d_0^2$ must
be of the form

$$\leqalignno{
d_0^2 = \pmatrix \underline d_0^1 & 0\\0&0\\0&0\\0&0\endpmatrix +
\pmatrix 0&0\\0& \underline d_0^3\\0&0\\0&0\endpmatrix +
\pmatrix 0& h_1\\ 0&0\\ 0 &h_2\\0&0\endpmatrix &&(1.36)\cr}$$
where $h_1: \Cal C^3_0 \rightarrow \Cal C^1_0$ and $h_2: \Cal C^3_0
\rightarrow
\Cal C^1_2$ and $d_1^2$ must be of the form

$$\leqalignno{
d_1^2 = \pmatrix 0&0&\underline d_1^1 & 0\\ 0&0&0&0\endpmatrix +
\pmatrix 0&0&0&0\\ 0&0&0&\underline d_1^3 \endpmatrix +
\pmatrix 0& h_3 &0& h_4\\ 0&0&0&0 \endpmatrix &&(1.37)\cr}$$
where $h_3 : \Cal C^3_0 \rightarrow \Cal C^1_2$ and $h_4 : \Cal C^3_2
\rightarrow
\Cal C^1_2.$

The assumption that $\Cal C^2$ is an exact sequence implies that

$$\leqalignno{
h_3 \underline d_0^3 + \underline d_1^1 h_2 = 0.&&(1.38)\cr}$$

In order to compute the torsion $T(\Cal C^2)$ we want to apply Proposition
1.3.
First note that $\underline d_0^1, \underline d_0^3, \underline d_1^1,
\underline d_1^3$ are all isomorphism as $\Cal C^1, \Cal C^2, \Cal C^3$ are
exact sequences.  In view of (1.38),

$$\leqalignno{
\qquad\pmatrix \underline d_0^1 & 0\\ 0 & \underline d_0^3\\0&0\\0&0
\endpmatrix
\pmatrix Id_{\Cal C_0^1}& (\underline d_0^1)^{-1}h_1\\ 0& Id_{\Cal C_0^3}
\endpmatrix
= \pmatrix Id& 0 & 0 & 0\\0 & Id& 0 & 0\\ 0 & (\underline d_1^1)^{-1}h_3 & Id
& 0\\
0 & 0 & 0 & Id \endpmatrix
\pmatrix \underline d_0^1 & h_1\\ 0 &\underline d_0^3\\ 0& h_2\\
0&0\endpmatrix
&&(1.39)\cr}$$
and

$$\leqalignno{
&\pmatrix Id & -h_4 (\underline d_1^3)^{-1}\\0&Id \endpmatrix
\pmatrix 0&h_3 & \underline d_1^1 & h_4 \\ 0&0&0& \underline d_1^3\endpmatrix
&(1.40)\cr
&\qquad= \pmatrix 0&0&\underline d_1^1 &0\\ 0 & 0 & 0 &\underline d_1^3
\endpmatrix
\pmatrix Id_{\Cal C_0^1} & 0 & 0 & 0 \\ 0 & Id_{\Cal C_0^3} & 0 & 0\\
0 &( \underline d_1^1)^{-1}h_3& Id_{\Cal C_2^1} &0\\ 0 & 0 & 0 & Id_{\Cal
C_2^3}
\endpmatrix. \cr}$$
Consider the following isomorphic three stage complexes: $\Cal C^2,
\Cal C^{21}, \Cal C^{22}$ and the isomorphisms
$f': \Cal C^2 \rightarrow \Cal C^{21},$
$f'':\Cal C^{21}
\rightarrow \Cal C^{22}$ given by the following commutative diagram (use
equality (1.38))

{\eightpoint
$$\eqalign{
\Cal C^2:\quad\, 0 \rightarrow \quad \Cal C^1_0 \oplus \Cal C^3_0 \quad
\overset
{\left(\smallmatrix \underline d_0^1& h_1\\0 & \underline d_0^3\\ 0& h_2\\
0&0
\endsmallmatrix\right)}\to \longrightarrow
\quad\Cal C^1_0 \oplus \Cal C^3_0 &\oplus \Cal C^1_2 \oplus
\Cal C^3_2 \quad \overset {\left(\smallmatrix 0 & h_3 & \underline d_1^1 &
h_4\\
0 & 0 & 0 & \underline d_1^3\endsmallmatrix\right)} \to \longrightarrow
\qquad\Cal C^1_2 \oplus \Cal C^3_2 \rightarrow 0\cr
\left(\smallmatrix Id& (\underline d_0^1)^{-1} h_1\\ 0& Id
\endsmallmatrix \right)\downarrow
\qquad\qquad\qquad\qquad\qquad\quad\,&\downarrow Id
\qquad\qquad\qquad\qquad\qquad\qquad \quad\downarrow
\left(\smallmatrix Id & -h_4
(\underline d_1^3)^{-1}\\0 & Id\endsmallmatrix\right) \cr
\Cal C^{21}:\quad 0 \rightarrow \quad
\Cal C^1_0 \oplus \Cal C^3_0\quad
\overset {\left(\smallmatrix \underline d_0^1 & 0\\ 0 & \underline d_0^3\\
0& h_2\\0&0\endsmallmatrix\right)} \to \longrightarrow
\quad\Cal C^1_0 \oplus \Cal C^3_0 &\oplus\Cal C^1_2 \oplus \Cal C^3_2\qquad
\overset {\left(\smallmatrix 0 & h_3 & \underline d_1^1 & 0\\ 0 & 0 & 0 &
\underline d_1^3
\endsmallmatrix\right)} \to \longrightarrow \quad
\Cal C^1_2 \oplus \Cal C^3_2 \rightarrow 0\cr
Id \downarrow \qquad\qquad\qquad\qquad\qquad\quad\,
&\downarrow \left(\smallmatrix
Id&0&0&0\\0&Id&0&0\\ 0&(\underline d_1^1)^{-1}h_3& Id&0\\0&0&0&Id
\endsmallmatrix\right)\qquad\qquad\quad\downarrow Id \cr
\Cal C^{22}:\quad 0 \rightarrow \quad \Cal C^1_0 \oplus \Cal C^3_0
\quad\overset{\left(\smallmatrix \underline d_0^1 & 0 \\ 0 & \underline
d_0^3\\
0&0\\0&0\endsmallmatrix\right)}
\to \longrightarrow \quad
\Cal C^1_0 \oplus \Cal C^3_0
&\oplus \Cal C^1_2 \oplus \Cal C^3_2\qquad
\overset {\left(\smallmatrix 0&0&\underline d_1^1&0\\ 0&0&0&\underline d_1^3
\endsmallmatrix\right)}
\to \longrightarrow \qquad \Cal C^1_2 \oplus \Cal C^3_2 \rightarrow 0.
\cr}$$}

Now apply Proposition 1.3 and use that all three complexes $\Cal C^2,
\Cal C^{21},
\Cal C^{22}$ are exact sequences to conclude that

$$\leqalignno{
\log T(\Cal C^2) = \log T(\Cal C^{21}) = \log T(\Cal C^{22})&&(1.41)\cr}$$

Moreover,

$$\leqalignno{
\log T(\Cal C^{22}) = \log T(\Cal C^1) + \log T(\Cal C^3).&&(1.42)\cr}$$
Combining (1.41) and (1.42) leads to (1.35).  \enddemo \hfill $\square$

The proof of Lemma 1.5 actually leads to a slightly stronger result:

\proclaim{Lemma 1.5$'$ }  Assume that the three stage complexes
$\Cal C^1, \Cal C^2, \Cal C^3$ are weakly
exact sequences of determinant class and, in addition, assume that
the morphisms
$\underline d_0^1, \underline d_1^1$ and $\underline d_1^3$ (or $\underline
d_0^1,
\underline d_0^3$ and $\underline d_1^3$) are isomorphisms.  Then

$$
\log T(\Cal C^2) = \log T(\Cal C^1) + \log T(\Cal C^3).$$
\endproclaim

\remark{Remark}
To conclude the result of Lemma $1.5'$
in the case where $\underline{d}_0^1$,
$\underline{d}_0^3$ and $\underline{d}_1^3$ are isomorphisms
one replaces the entry $(\underline d_1^1)^{-1} h_3$
in the 4 x 4 matrix which appears in the diagram above
by $-h_2(\underline d_0^3)^{-1}.$

The following considerations are a preparation for the proof of Theorem 1.14.

Assume that

$$
0 \rightarrow \Cal C^1 \overset f \to \rightarrow \Cal C^2 \overset g \to
\rightarrow \Cal C^3 \rightarrow 0$$
is an exact sequence of $\Cal A$-cochain complexes.  Consider the Hodge
decomposition of each of the three complexes, $k= 1, 2, 3$,

$$
\Cal C_i^k = \Cal H^k_i \oplus \Cal C_i^{k, +} \oplus \Cal C_i^{k, -}.$$
With respect to these decompositions, $f_i: \Cal C_i^1
\rightarrow \Cal C_i^2$
and $g_i : \Cal C_i^2 \rightarrow\Cal C_i^3$ take the form

$$\leqalignno{
f_i = \pmatrix f_{i, 11} & 0 & f_{i, 13}\\ f_{i, 21} & f_{i, 22} & f_{i,
23}\\
0&0&f_{i, 33}\endpmatrix, \,\,
g_i = \pmatrix g_{i, 11} & 0 & g_{i, 13}\\ g_{i, 21} & g_{i, 22} & g_{i,
23}\\
0&0&g_{i, 33}\endpmatrix .&&(1.43)\cr}$$
Moreover

$$\leqalignno{
\underline d_i^2 f_{i, 33} = f_{i+1, 22} \underline d_i^1;\, \underline d_i^3
g_{i, 33} = g_{i+1, 22} \underline d_i^2.&&(1.44)\cr}$$

The exactness of the sequences $0 \rightarrow \Cal C_k^1 \rightarrow
\Cal C_k^2
\rightarrow \Cal C_k^3 \rightarrow 0$ imply that the
following sequences

$$\leqalignno{
0 \rightarrow \Cal C_i^{1,+} \overset {f_{i, 22}}
\to \rightarrow \Cal C_i^{2,+}
\overset {g_{i, 22}} \to \rightarrow \Cal C_i^{3,+}\rightarrow
0&&(1.45)\cr}$$

$$\leqalignno{
0 \rightarrow \Cal C_i^{1,-} \overset {f_{i, 33}}
\to \rightarrow \Cal C_i^{2,-}
\overset {g_{i, 33}} \to \rightarrow \Cal C_i^{3,-}\rightarrow
0&&(1.46)\cr}$$
are three stage complexes
with the property that $f_{i, 22}$ is 1-1 and $g_{i, 33}$ is onto.  Using
(1.44)
one verifies that $f_{i,22}$ has closed image, the range of $g_{i,22}$ is
dense in
$\Cal C_i^{3, +}$ and $f_{i, 33}$ is 1-1.  Therefore, the cochain complexes
(1.45) and
(1.46) have a Hodge decomposition

$$\leqalignno{
0 \rightarrow \Cal C_i^{1,+} \overset {f_{i, 22}}
\to \rightarrow \Cal C_i^{2,+0}
\oplus \Cal C_i^{2, ++} \oplus \Cal C_i^{2, +-} \overset {g_{i, 22}} \to
\rightarrow \Cal C_i^{3,+}\rightarrow 0&&(1.47)\cr}$$
and

$$\leqalignno{
0 \rightarrow \Cal C_i^{1,-} \overset {f_{i, 33}}
\to \rightarrow \Cal C_i^{2,-0}
\oplus \Cal C_i^{2, -+} \oplus \Cal C_i^{2, --} \overset {g_{i, 33}} \to
\rightarrow \Cal C_i^{3,-}\rightarrow 0.&&(1.48)\cr}$$

With respect to this decomposition,

$$\leqalignno{
f_{i, 22} = \pmatrix 0\\ \alpha_{i, 43} \\ 0 \endpmatrix;\,
g_{i, 22} = \pmatrix 0&0&\beta_{i, 35}\endpmatrix &&(1.50)\cr}$$

$$\leqalignno{
f_{i, 33} = \pmatrix 0\\ \alpha_{i, 74} \\ 0 \endpmatrix;\,
g_{i, 33} = \pmatrix 0&0&\beta_{i, 48}\endpmatrix. &&(1.51)\cr}$$
Notice that $\alpha_{i, 43}$ and $\beta_{i, 48}$ are isomorphisms and
$\beta_{i, 35}$ and $\alpha_{i, 74}$ are weak isomorphisms.

The exact sequence of $\Cal A$-cochain complexes of finite type

$$
0 \rightarrow \Cal C^1 \overset f \to \rightarrow \Cal C^2 \overset g \to
\rightarrow \Cal C^3 \rightarrow 0$$
induces a long weakly exact sequence $\Cal H$ in cohomology (cf. [CG, p. 10,
Thm 2.1])
$$\leqalignno{
\dots \rightarrow H^i(\Cal C^1) \overset H(f_i) \to \rightarrow H^i(\Cal C^2)
\overset H(g_i) \to \rightarrow H^i(\Cal C^3) \overset H(\delta_i) \to
\rightarrow H^{i+1}(\Cal C^1)
\rightarrow \dots &&(1.52)\cr}$$
or, equivalently, in terms of the harmonic spaces $\Cal H_i^k :=
Null(d_i^k) \cap Null(d_{i-1}^{k*})$

$$\leqalignno{
\dots \rightarrow \Cal H^1_i \overset {f_{i, 11}}
\to \rightarrow \Cal H^2_i
\overset {g_{i,11}} \to \rightarrow \Cal H^3_i \overset {\delta_i} \to
\rightarrow \Cal H^1_{i+1} \overset {f_{i+1, 11}}
\to \rightarrow \Cal H^2_{i+1}
\rightarrow \dots &&(1.52')\cr}$$
In particular the cochain complex $\Cal H$ is acyclic.  Consider its Hodge
decomposition

$$\leqalignno{
\dots \rightarrow \Cal H_i^{1,+} \oplus \Cal H_i^{1, -} \overset
{(\smallmatrix
0 &\alpha_{i, 12}\\ 0 & 0\endsmallmatrix)
}\to \longrightarrow \Cal H_i^{2, +}
\oplus \Cal H_i^{2, -} \overset{(\smallmatrix 0 &\beta_{i, 12}\\ 0 & 0
\endsmallmatrix)} \to \longrightarrow &&(1.53)\cr
\Cal H_i^{3, +} \oplus \Cal H_i^{3, -} \overset{(\smallmatrix 0 & \gamma_i\\
0&0\endsmallmatrix)} \to \longrightarrow \Cal H_{i+1}^{1, +} \oplus \Cal
H_{i+1}^{1, -} \rightarrow \dots\cr}$$
Notice that $\alpha_{i, 12}, \beta_{i, 12}$ and $\gamma_i$ are weak
isomorphisms.

In view of (1.43) - (1.53) $f_i$ and $g_i$ can be written as

$$\leqalignno{
\qquad f_i = \pmatrix f_{i, 11} & 0 & f_{i, 13}\\f_{i, 21}&f_{i, 22}& f_{i,
23}\\
0&0&f_{i, 33}\endpmatrix =
\pmatrix
{\pmatrix 0 & \alpha_{i, 12}\\ 0&0 \endpmatrix}
&{\pmatrix 0\\0\endpmatrix}
&{\pmatrix \alpha_{i, 14}\\\alpha_{i, 24}\endpmatrix}\\
{\pmatrix \alpha_{i, 31}& \alpha_{i, 32}\\
\alpha_{i, 41}& \alpha_{i, 42}\\ \alpha_{i, 51}& \alpha_{i, 52}\endpmatrix}
&{\pmatrix 0\\ \alpha_{i, 43}\\0\endpmatrix}
&{\pmatrix \alpha_{i, 34}\\ \alpha_{i, 44}\\ \alpha_{i, 54}\endpmatrix}\\
{\pmatrix 0&0\\0&0\\0&0\endpmatrix}
&{\pmatrix 0\\0\\0\endpmatrix}
&{\pmatrix 0\\ \alpha_{i, 74}\\0\endpmatrix}
\endpmatrix&&(1.54)\cr}$$
and

$$\leqalignno{
\quad g_i &= \pmatrix g_{i, 11} & 0 & g_{i, 13}\\g_{i, 21}&g_{i, 22}& g_{i,
23}\\
0&0&g_{i, 33}\endpmatrix &(1.55)\cr
&=
\pmatrix
{\pmatrix 0 & \beta_{i, 12}\\ 0&0 \endpmatrix}
&{\pmatrix 0\,&0\,&0\\0\,&0\,&0\endpmatrix}
&{\pmatrix \beta_{i, 16}& \beta_{i, 17}& \beta_{i, 18}\\
\beta_{i, 26} & \beta_{i, 27} & \beta_{i, 28}\endpmatrix}\\
{\pmatrix \beta_{i, 31} & \beta_{i, 32}\endpmatrix}
&\quad{\pmatrix 0&\, 0&\beta_{i, 35}\endpmatrix}
&{\pmatrix \beta_{i, 36} & \beta_{i, 37}& \beta_{i, 38}\endpmatrix}\\
{\pmatrix 0\, & 0\endpmatrix}
&{\pmatrix 0\,&0\,&0\endpmatrix}
&{\pmatrix 0\qquad&0\quad&\beta_{i, 48}\endpmatrix}
\endpmatrix \cr}$$

Recall that $f_i: \Cal H^{1, +}_i \oplus \Cal H_i^{1, -}
\oplus \Cal C_i^{1, +}
\oplus C_i^{2, -} \rightarrow \Cal H_i^{2, +} \oplus \Cal H_i^{2, -}
\oplus \Cal C_i^{2, +0} \oplus \Cal C_i^{2, ++} \oplus \Cal C_i^{2, +-}
\oplus \Cal C_i^{2, -0} \oplus \Cal C_i^{2, -+}\oplus \Cal C_i^{2, --}$ and
$g_i: \Cal H_i^{2, +} \oplus \Cal H_i^{2, -}
\oplus \Cal C_i^{2, +0} \oplus \Cal C_i^{2, ++} \oplus \Cal C_i^{2, +-}
\oplus \Cal C_i^{2, -0} \oplus \Cal C_i^{2, -+}\oplus \Cal C_i^{2, --}
\rightarrow \Cal H^{3, +}_i \oplus \Cal H_i^{3, -} \oplus \Cal C_i^{3, +}
\oplus C_i^{3, -}. $

We have already noticed that $\alpha_{i, 12}, \alpha_{i, 74}, \beta_{i, 12},
\beta_{i, 35}$ are weak isomorphisms, while $\alpha_{i, 43}$ and
$\beta_{i,48}$ are isomorphisms.  Using that $g_i \cdot f_i =0$ and the fact
that $\alpha_{i, 74}$ and $\beta_{i, 35}$ are weak isomorphisms one verifies
that

$$\leqalignno{
\alpha_{i, 51} = 0;\quad \beta_{i, 27} = 0&&(1.56)\cr}$$
$$\leqalignno{
&\beta_{i, 12} \alpha_{i, 24} + \beta_{i, 17} \alpha_{i, 74} = 0;\quad
\beta_{i, 31}\alpha_{i, 12} + \beta_{i, 35} \alpha_{i, 52} = 0 &(1.56')\cr
&\beta_{i, 31} \alpha_{i, 14} + \beta_{i, 32} \alpha_{i, 24} + \beta_{i, 35}
\alpha_{i, 54} + \beta_{i, 37} \alpha_{i, 74} = 0.\cr}$$

\proclaim{Lemma 1.6}  For any $i$,

\noindent (A) $\alpha_{i, 31} : \Cal H_i^{1, +}
\rightarrow \Cal C_i^{2, +0}$
is an isomorphism;

\noindent (B) $\beta_{i, 26} : \Cal C_i^{2, -0}
\rightarrow \Cal H_i^{3, -}$
is an isomorphism.
\endproclaim

\demo{Proof}  (A)  Inspecting (1.54) and in view of $\alpha_{i, 51} = 0$
(cf. 1.56) one concludes that $h_i := f_i
\restriction_{\Cal H_i^{1, +}\oplus \Cal C_i^{1,+}}$ has a range which is
contained in $\Cal C_i^{2, +0} \oplus \Cal C_i^{2, ++}$ and admits a
representation of the form

$$
h_i := \pmatrix \alpha_{i, 31} & 0\\  \alpha_{i, 41}& \alpha_{i, 43}
\endpmatrix : \Cal H_i^{1, +} \oplus \Cal C_i^{1,+}
\rightarrow \Cal C_i^{2, +0}
\oplus \Cal C_i^{2, ++}.$$

{}From (1.55) one sees that $\Cal C_i^{2, +0} \oplus \Cal C_i^{2, ++} \subseteq
Null \, g_i.$  Due to the exactness $0 \rightarrow \Cal C_i^1 \overset {f_i}
\to \rightarrow \Cal C_i^2 \overset {g_i} \to \rightarrow \Cal C_i^3
\rightarrow 0$, the map $h_i$ is an isomorphism.  Further we already know
that $\alpha_{i, 43}$ is an isomorphism and therefore conclude that
$\alpha_{i, 31}: \Cal H_i^{1, +} \rightarrow \Cal C_i^{2, +0}$ is an
isomorphism as well.

\noindent (B) From (1.54) one sees that $\Cal C_i^{2, -0}
\oplus \Cal C_i^{2, --}
\subseteq (Range \, f_i)^\perp$.  Therefore, using that $0 \rightarrow
\Cal C_i^1 \overset {f_i} \to \rightarrow \Cal C_i^2
\overset {g_i} \to \rightarrow \Cal C_i^3 \rightarrow 0$ is exact,
$h_i' := g_i \restriction_{\Cal C_i^{2, -0} \oplus \Cal C_i^{2, --}}$ is 1-1
and
its image is given by $\Cal H_i^{3, -} \oplus \Cal C_i^{3, -}$, i.e.
$h_i' := \pmatrix \beta_{i, 26} & \beta_{i, 28}\\ 0& \beta_{i, 48}
\endpmatrix
: \Cal C_i^{2, -0} \oplus \Cal C_i^{2, --} \rightarrow \Cal H_i^{3, -} \oplus
\Cal C_i^{3, -}$ is bijective.  We already know that $\beta_{i, 48}$ is an
isomorphism and (B) follows. \hfill $\square$
\enddemo

Further notice that $d_i^2$ takes the following form

$$\leqalignno{
d_i^2 = \pmatrix 0&0&0\\ 0&0&\underline d_i^2\\ 0&0&0\endpmatrix =
\pmatrix 0&0&0\\ &&&\\ 0&0& {\pmatrix d^2_{i, 36}&0&d^2_{i, 38}\\
d^2_{i, 46}& d^2_{i, 47} & d^2_{i, 48}\\ 0&0&d^2_{i, 58}\endpmatrix}\\
&&&\\ 0&0&0\endpmatrix &&(1.57)\cr}$$
where we used again (1.44) to verify that $d^2_{i, 56}, d^2_{i, 57}$ and
$d^2_{i, 37}$ all vanish and that, for any $i$,

$$\leqalignno{
\underline d_i^3 \beta_{i, 48} = \beta_{i+1, 35} d^2_{i, 58} \text
{ is a weak isomorphism}&&(1.58)\cr}$$
and

$$\leqalignno{
d^2_{i, 47} \alpha_{i, 74} = \alpha_{i+1, 43} \underline d_i^1 \text
{ is a weak isomorphism.}&&(1.59)\cr}$$
Further, considering $\Cal C_i^{2, -} \overset {g_{i, 13}} \to
\longrightarrow
\Cal H^3_i \overset {\delta_i} \to \longrightarrow \Cal H^1_{i+1} \overset
{f_{i+1, 21}} \to \longrightarrow \Cal C_{i+1}^{2, +}$ and taking into
account
the
definition of $\delta_i = \pmatrix 0& \gamma_i\\ 0&0\endpmatrix$ one obtains
$$\leqalignno{
d^2_{i, 36} = \alpha_{i+1, 31} \gamma_i \beta_{i, 26}.&&(1.60)\cr}$$

In a straightforward way one verifies the following

\proclaim{Lemma 1.7} Assume that $F : \Cal W_1 \oplus \Cal W_2
\rightarrow \Cal V_1 \oplus \Cal V_2$ is an $\Cal A$-linear map of the form
$F= \pmatrix f & h\\ 0&g\endpmatrix$.

\noindent (A) If $F$ is 1-1 and has closed range and $f$ is an isomorphism,
then
$g$ is 1-1 and has closed range.

\noindent (B) If $F$ is onto and $g$ is an isomorphism, then $f$ is onto.
\endproclaim

\proclaim{Lemma 1.8}  For any $i$, $\alpha_{i, 12}$ and $\alpha_{i, 74}$ are
of
determinant class iff $\beta_{i, 12}$ and $\beta_{i, 35}$ are of determinant
class.
\endproclaim

\demo{Proof}  Introduce $\varphi_i : \Cal H_i^{1, -}
\oplus \Cal C_i^{1, -}
\rightarrow \Cal H_i^{2, +} \oplus \Cal C_i^{2, -+}
\oplus \Cal H_i^{2, -}
\oplus \Cal C_i^{2, +-}$ and $\psi_i : \Cal H _i^{2, +}
\oplus \Cal C_i^{2, -+}
\oplus \Cal H_i^{2, -} \oplus \Cal C_i^{2, +-} \rightarrow
\Cal H_i^{3, +} \oplus \Cal C_i^{3, +}$ given as follows

$$\leqalignno{
\psi_i = \pmatrix 0 &\beta_{i, 17} & \beta_{i, 12} & 0\\ \beta_{i, 31}&
\beta_{i, 37} & \beta_{i, 32} & \beta_{i, 35} \endpmatrix ; \varphi_i =
\pmatrix \alpha_{i, 12}&\alpha_{i, 14}\\ 0& \alpha_{i, 74}\\
0 & \alpha_{i, 24}\\ \alpha_{i, 52} & \alpha_{i, 54}\endpmatrix .
&&(1.61)\cr}$$
Due to (1.56) - $(1.56')$ we have $\psi_i \cdot \varphi_i = 0$.  Moreover
$\alpha_{i, 12} , \alpha_{i, 74} , \beta_{i, 12}$ and $\beta_{i, 35}$ are
all weak isomorphisms.  Therefore, for any $i$,

$$\leqalignno{
\quad 0 \rightarrow \Cal H_i^{1, -} \oplus \Cal C_i^{1, -} \overset
{\varphi_i}
\to\longrightarrow \Cal H_i^{2, +} \oplus \Cal C_i^{2, -+}
\oplus \Cal H_i^{2, -}
\oplus \Cal C_i^{2, +-} \overset {\psi_i}
\to \longrightarrow \Cal H_i^{3,+}
\oplus \Cal C_i^{3, +} \rightarrow 0&&(1.62)\cr}$$
is a weak exact sequence.

Using that $\pmatrix \alpha_{i, 31} & 0\\ \alpha_{i, 41} & \alpha_{i, 43}
\endpmatrix$ is an isomorphism and applying Lemma 1.7 (A) twice
to $f_i$ one concludes
that $\varphi_i$ has closed image.  Similarly, using that $\pmatrix
\beta_{i, 26} & \beta_{i, 28}\\ 0 & \beta_{i, 48}\endpmatrix$ is an
isomorphism
and applying Lemma 1.7 (B) one concludes
that $\psi_i$ is onto.  Therefore, the sequence (1.62) is, in fact, an
exact sequence and thus of determinant class.  The statement then follows
from Lemma 1.4(B) and Proposition 1.1(B). \enddemo \hfill $\square$

\proclaim{Lemma 1.9}  For any $i, d^2_{i, 36}, d^2_{i, 47}, d^2_{i, 58}$ are
weak
isomorphisms and satisfy:

\noindent (A) $d^2_{i, 36}$ is of determinant class iff $\gamma_i$ is of
determinant class.

\noindent (B) $d^2_{i, 47}$ and $\alpha_{i, 74}$ are of determinant class iff
$\underline{d_i^1}$ is of determinant class.

\noindent (C) $d^2_{i, 58}$ and $\beta_{i+1, 35}$ are of determinant class
iff
$\underline{d_i^3}$ is of determinant class.
\endproclaim

\demo{Proof}  The fact that $d^2_{i, 36}, d^2_{i, 47}$ and $d^2_{i, 58}$ are
weak isomorphisms follows from inspection of (1.57).  To prove (A), (B) and
(C)
one applies Proposition 1.1 and uses formula (1.60) (in case of (A)), formula
(1.59) (in case of (B)) and formula (1.58) (in case of (C)) together with
the properties that $\alpha_{i+1,31}$ and $\beta_{i, 26}$ are isomorphisms
(Lemma 1.6) and therefore of determinant class (in case (A)),
$\alpha_{i+1, 43}$ is an isomorphism (in case (B)) and $\beta_{i, 48}$ is
an isomorphism (in case (C)). \enddemo \hfill $\square$

\proclaim{Lemma 1.10} Assume that $\Cal C^1$ and $\Cal C^3$ are of
determinant class.  Then, for any $i$, the following statements hold:

\noindent (A) the three stage complexes $\Cal D^i_{\pm}$ given by $0
\rightarrow
\Cal C_i^{1, \pm} \rightarrow \Cal C_i^{2, \pm} \rightarrow
\Cal C_i^{3, \pm} \rightarrow 0$ are of determinant class and

$$\leqalignno{
\log T ( \Cal D^i_{+}) = \log \Vol \,\alpha_{i, 43} - \log \Vol \,\beta_{i,
35}
&&(1.63)\cr}$$

$$\leqalignno{
\log T ( \Cal D^i_{-}) = \log \Vol \,\alpha_{i, 74} - \log \Vol \,\beta_{i,
48}.
&&(1.64)\cr}$$

\noindent (B) Assume, in addition, that the cochain complex
$\Cal H$ (cf (1.52)) is of determinant class.
Then, for any $i$, the three stage complex $\Cal D^i,$ given by $0
\rightarrow
\Cal C^1_i \overset {f_i} \to \longrightarrow \Cal C^2_i \overset {g_i} \to
\longrightarrow \Cal C^3_i \rightarrow 0,$ is of determinant class and

$$\leqalignno{
\log T ( \Cal D^i)& = \log \Vol \,\alpha_{i, 12} + \log \Vol \,\alpha_{i, 31}
+ \log \Vol \,\alpha_{i, 43} + \log \Vol \,\alpha_{i, 74} &(1.65)\cr
&\quad - \log \Vol \,\beta_{i, 12} - \log \Vol \,\beta_{i, 26} - \log \Vol\,
\beta_{i, 35} - \log \Vol \,\beta_{i, 48}.\cr}$$
\endproclaim

\demo{Proof} (A) Recall that $\Cal D^i_{ +}$ is given by

$$
0 \rightarrow \Cal C_i^{1, +} \overset {f_{i, 22}} \to \longrightarrow
\Cal C_i^{2, +} \overset {g_{i, 22}} \to \longrightarrow \Cal C_i^{3, +}
\rightarrow 0$$
where, according to (1.50), $f_{i, 22} = (\smallmatrix 0\\ \alpha_{i, 43} \\
0 \endsmallmatrix)$ and $g_{i, 22} = ( 0 \quad 0 \quad \beta_i,_{35})$.  The
map $\alpha_{i, 43}$ is an isomorphism and therefore of determinant class.
Due
to Lemma 1.9, the weak isomorphism $\beta_{i, 35}$ is of determinant class
as well and (1.63) follows.

Similarly one argues for $\Cal D^i_{-}$ given by

$$
0 \rightarrow \Cal C_i^{1, -} \overset {f_{i, 33}} \to \longrightarrow
\Cal C_i^{2, -} \overset {g_{i, 33}} \to \longrightarrow \Cal C_i^{3, -}
\rightarrow 0$$
where, according to (1.50),

$$f_{i, 33} = \pmatrix 0\\ \alpha_{i, 74} \\ 0 \endpmatrix, \,\,
g_{i, 33} = \pmatrix 0&0& \beta_{i, 48}\endpmatrix.$$
Notice
that $\beta_{i, 48}$ is an isomorphism and therefore of determinant class.
According to Lemma 1.9, the weak isomorphism $\alpha_{i, 74}$ is of
determinant class and therefore (1.64) is proved.

\noindent (B) The complexes $\Cal D^i$, given by $0
\rightarrow \Cal C^1_i
\overset {f_i} \to \longrightarrow \Cal C^2_i \overset {g_i} \to
\longrightarrow \Cal C^3_i \rightarrow 0$ are exact and therefore of
determinant
class.  To verify formula (1.65) we want to apply Lemma 1.4.

Decompose $\Cal C^2_i = \Cal W^2_i \oplus \Cal V^2_i$ where

$$
\Cal W_i^2 := \Cal H_i^{2, +} \oplus \Cal C_i^{2, +0}
\oplus \Cal C_i^{2, ++}
\oplus \Cal C_i^{2, -+}$$
and

$$\Cal V_i^2 := (\Cal W_i^2)^\perp = \Cal H_i^{2, -}
\oplus \Cal C_i^{2, +-}
\oplus \Cal C_i^{2, -0} \oplus \Cal C_i^{2, --}.$$
Then $f_i = \pmatrix f_i'\\ f_i '' \endpmatrix$ and $g_i = (g_i' \,\,
g_i'')$ have the following representations

$$\leqalignno{
f_i' = \pmatrix 0& \alpha_{i, 12} & 0 & \alpha_{i, 14}\\
\alpha_{i, 31} & \alpha_{i, 32} & 0 & \alpha_{i, 34} \\
\alpha_{i, 41} & \alpha_{i, 42} & \alpha_{i, 43} & \alpha_{i, 44} \\
0& 0& 0& \alpha_{i, 74} \endpmatrix &&(1.66)\cr}$$

$$\leqalignno{
g_i'' = \pmatrix \beta_{i, 12} & 0 & \beta_{i, 16} & \beta_{i, 18} \\
0 & 0 & \beta_{i, 26} & \beta_{1, 28}\\
\beta_{i, 32}& \beta_{i, 35} & \beta_{i, 36} & \beta_{i, 38}\\
0&0&0& \beta_{i, 48}\endpmatrix . &&(1.67)\cr}$$

Due to the assumption that $\Cal H$ is of determinant class and due to
Lemma 1.6 and Lemma 1.9, the maps $\alpha_{i, 31}, \alpha_{i, 12}, \alpha_{i,
43},
\alpha_{i, 74}, \beta_{i, 12}, \beta_{i, 35} , \beta_{i, 26} \text { and }
\beta_{i, 48}$
are weak isomorphisms of determinant class.  Therefore $f_i'$ and $g_i''$ are
weak isomorphisms of determinant class.  As $g_i$ is onto, we can apply
Proposition 1.4 to conclude that

$$\leqalignno{
\log T(0 \rightarrow \Cal C_i^1 \rightarrow \Cal C_i^2 \rightarrow
\Cal C_i^3 \rightarrow 0) = \log \Vol f_i' - \log \Vol g_i''.&&(1.68)\cr}$$

Applying Proposition 1.1, we conclude from (1.66) - (1.67)

$$\leqalignno{
\log \Vol f_i' &= \log \Vol \alpha_{i, 31} + \log \Vol \alpha_{i, 12} +
\log \Vol \alpha_{i, 43} + \log \Vol \alpha_{i, 74} &(1.69)\cr
\log \Vol g_i '' &= \log \Vol \beta_{i, 12} + \log \Vol \beta_{i, 35} +
\log \Vol \beta_{i, 26} + \log \Vol \beta_{i, 48}.\cr}$$
Substituting (1.69) into (1.68) leads to the claimed result (1.65). \hfill
$\square$
\enddemo

Let $f:\Cal C^1\to \Cal C^2$ be a morphism of cochain complexes. Define the
mapping cone of $f$ to be the cochain complex $\Cal C(f):=
(\Cal C_i(f),
d_i(f))$ where $\Cal C_i(f)=\Cal C_i^2\oplus \Cal C_{i+1}^1,$
$d_i(f) = \pmatrix d^2_i& f_{i+1}\\ 0&-d^1_{i+1}\endpmatrix.$
Denote by $ S\Cal C=(S\Cal C_i, Sd_i)$ the
suspension of a cochain complex $\Cal C$,
defined by
$S\Cal C_i= \Cal C_{i+1}, S d_i=-d_{i+1}.$
\newline Consider the suspension $S\Cal C^1$ of $\Cal C^1$ and the morphisms
$p(f):\Cal C(f)\to S\Cal C^1$ and
$j(f):\Cal C^2\to \Cal C(f)$ given by
$p_i(f)=(0\, Id)$ and $j_i(f)=\pmatrix Id \\ 0 \endpmatrix.$
Then
$$
0 \rightarrow \Cal C^2 \overset j(f) \to \rightarrow \Cal C(f)
\overset p(f) \to
\rightarrow S\Cal C^1 \rightarrow 0$$
is a short exact sequence of cochain complexes. The induced
long weakly exact sequence in cohomology (cf. 1.52) becomes in this case
$$\leqalignno{
\dots \rightarrow H^{i}(\Cal C^2)  \rightarrow H^{i}(\Cal C(f))
\rightarrow H^{i}(S\Cal C^1) \overset {H(\delta_i)}
\to \rightarrow H^{i+1}(\Cal C^2)
\rightarrow \dots &&(1.70)\cr}$$
Notice that $H^{i}(S\Cal C^1)$ is equal to $H^{i+1} (\Cal C^1)$ and
the connecting homomorphism
$H(\delta_i): H^i(S\Cal C^1) \to H^{i+1}(\Cal C^2)$ is given by the
map
$H(f_{i+1}):H^{i+1}(\Cal C^1)\to H^{i+1}(\Cal C^2),$
induced by f in cohomology.
\proclaim{ Lemma 1.11 (A)} Let $f: \Cal C^1 \to \Cal C^2$
be a morphism of cochain complexes. Assume that
 the cochain complexes $\Cal C^k,
(k=1,2)$ are of determinant class as well as the morphisms $H(f_i)
\,(i \ge 0).$   Then the
 cochain complex $\Cal C(f)$ is of determinant
class.
\endproclaim

\demo{Proof} The claimed result follows from applying Proposition 1.1 and
Lemma 1.9 (cf the proof of Proposition 1.13(i)).
\hfill $\square$
\enddemo

Let
$$
0 \rightarrow \Cal C^1 \overset f \to \rightarrow \Cal C^2
\overset g \to
\rightarrow \Cal C^3 \rightarrow 0$$
be a short exact sequence of cochain complexes. Let
$\pi:\Cal C(f)\to \Cal C^3$, $\rho :S\Cal C^1\to \Cal C(g)$ and
$S(f): S\Cal C^1 \to S\Cal C^2$ be the morphisms defined by
$\pi_i := (g_i \, 0),$  $\rho_i := \pmatrix 0 \\  f_{i+1} \endpmatrix,$
$S_i(f) := f_{i+1.}$ Clearly
$\pi_i\cdot j_i(f) = g_i$ and $p_i(g)\cdot \rho_i= S_i(f).$

\proclaim{ Lemma 1.11 (B)} Both $\pi$ and $\rho$ are homotopy
equivalences.
\endproclaim

\demo{Proof} Without loss of generality one can assume that
$\Cal C_i^2=\Cal C_i^1\oplus \Cal C_i^3$, $f_i = \pmatrix Id \\ 0
\endpmatrix$,
$g_i = (0 \, Id)$  and
$d_i^2 = \pmatrix d_i^1& \theta_i\\ 0&d_i^3 \endpmatrix$ with
$\theta_i:\Cal C_i^3\to \Cal C_{i+1}^1,$  satisfying

$$\leqalignno{
d^1_{i+1}\theta_i+\theta_{i+1} d^3_i=0. &&(1.71)\cr}
$$

Then $\rho_i = \pmatrix 0\\ Id \\ 0 \endpmatrix $ and
$\pi_i =(0 \, Id \, 0
),$ as well as $\Cal C_i(f)=C_i^1\oplus \Cal C_i^3\oplus \Cal C_{i+1}^1$
and
$\Cal C_i(g)=\Cal C_i^3\oplus \Cal C_{i+1}^1\oplus \Cal C_{i+1}^3.$
Notice that $d_i(f),$ resp. $d_i(g),$ are given by
$$\leqalignno{
d_i(f)= \pmatrix d_i^1&\theta_i&
Id\\ 0&d^3_{i}&0\\ 0&0&-d^1_{i+1}\endpmatrix
. &&(1.72)\cr}$$
$$\leqalignno{
d_i(g)= \pmatrix d_i^3&0&Id\\
0&-d_{i+1}^1&-\theta_{i+1}\\ 0&0&-d^3_{i+1}\endpmatrix
. &&(1.72')\cr}$$

Define the morphisms of cochain complexes
$\sigma:\Cal C^3\to \Cal C(f)$  and
$\omega:\Cal C(g)\to S\Cal C^1$
by $\sigma_i = \pmatrix 0\\ Id \\ -\theta_i \endpmatrix $
respectively
$\omega_i=(-\theta_i \, -Id \,\, 0)$ and notice that $\pi_i
\cdot \sigma_i =Id$ and
$\omega_i \cdot \rho_i =-Id.$  It remains to check that $\sigma\cdot \pi$
and $-\rho \cdot \omega$ are both homotopic  to $Id.$
Since the cochain complexes $Null (\pi)$
and $Null (\omega)$ are, as is verified easily, exact sequences
this will follow from
Lemma 1.12 below
(applied to $0 \to Null(\pi) \to \Cal C(f) \to \Cal C^3 \to 0,$
respectively $0 \to Null(\omega) \to \Cal C(g) \to S\Cal C^1 \to 0).$
 \hfill $\square$
\enddemo

\proclaim{ Lemma 1.12} Let
$$
0 \rightarrow \Cal C^1 \overset f \to \rightarrow \Cal C^2
\overset g \to
\rightarrow \Cal C^3 \rightarrow 0$$
be a short exact sequence of cochain complexes and
let $s:\Cal C^3\to \Cal C^2$ be a morphism so that $g\cdot s=Id$. If
$\Cal C^1$ is an exact sequence, then $s\cdot g$ and $Id$ are
homotopic.
\endproclaim

\demo{Proof} It is to prove (cf. Definition 1.14 in [BFKM]) that
there exist morphisms $t_i:\Cal C_i^2 \to \Cal C_{i-1}^2$ such that
$Id - s_ig_i = d^2_{i-1}t_i + t_{i+1}d^2_i.$
Consider the Hodge decomposition of $\Cal C^1,
\Cal C^1_i=\Cal C_i^{1,+}\oplus \Cal C_i^{1,-}$ with $d_i^1 = \pmatrix 0
&\underline d_i^1 \\ 0&
0\endpmatrix.$ In view of the exactness of $\Cal C^1,$
$\underline d_i^1$
are isomorphisms. Without loss of generality one can assume that
$\Cal C_i^2 = \Cal C^{1,+}_i\oplus \Cal C^{1,-}_i\oplus \Cal C^3_i,$
$f_i =
\pmatrix Id&0 \\ 0& Id \\
0 & 0\\ \endpmatrix,$ $g_i=(0 \,0 \, Id)$, $s_i = \pmatrix s^{+}_i\\
s^{-}_i \\ Id\endpmatrix $ and  $d_i^2 =
\pmatrix 0&\underline d_i^1&
\theta_i^{+}\\
0&0&\theta_i^{-}\\ 0&0&d^3_i\endpmatrix.$
Note that $d_{i+1}^2\cdot d_i^2= 0$ implies that
$$\leqalignno{
\underline d^1_{i+1}\theta_i^- +\theta_{i+1}^+d_i^3=0 ,
\,\,\theta^-_{i+1} d^3_i=0. &&(1.73)\cr}
$$
The condition that s is a morphism of cochain
complexes implies that
$$\leqalignno{
\underline d^1_{i} s^-_i +\theta_{i}^+ =s^+_{i+1}d^3_i, \,\, \theta^-_i=
s^-_{i+1} d^3_i. &&(1.74)\cr}$$

Define $t_i:\Cal C_i^{1,+}\oplus \Cal C_i^{1,-}\oplus \Cal C^3_i\to
\Cal C_{i-1}^{1,+}\oplus \Cal C_{i-1}^{1,-}\oplus \Cal C^3_{i-1} $
by

$$\leqalignno{
t_{i+1}= \pmatrix 0&0&
0\\ (\underline d^1_i)^{-1}&0&-(\underline d^1_i)^{-1} \cdot
s^+_{i+1}\\ 0&0&0\endpmatrix
.&&(1.75)\cr}$$
One verifies that
$Id-s_i g_i=t_{i+1} d^2_i +d^2_{i-1}t_i.$
\hfill $\square$
\enddemo
Let $0 \rightarrow \Cal C^1 \overset f \to
\rightarrow \Cal C^2 \overset g \to \rightarrow \Cal C^3 \rightarrow 0$
be an exact sequence of $\Cal A$-cochain complexes of finite type.
Recall that it
induces the long weakly exact sequence in cohomology (1.52)
$$\leqalignno{
\dots \rightarrow H^i(\Cal C^1) \overset H(f_i) \to \rightarrow H^i(\Cal C^2)
\overset H(g_i) \to \rightarrow H^i(\Cal C^3) \overset H(\delta_i) \to
\rightarrow H^{i+1}(\Cal C^1)
\rightarrow \dots &&(1.52)\cr}$$
which will be also viewed as a cochain complex and denoted by
$\Cal H.$

\proclaim{Proposition 1.13}  Let $0 \rightarrow \Cal C^1 \overset f \to
\rightarrow \Cal C^2 \overset g \to \rightarrow \Cal C^3 \rightarrow 0$
be an exact sequence of $\Cal A$-cochain complexes of finite type.
\newline (i) If $\Cal C^1$, $\Cal C^3$ is of determinant class and
$H(\delta_i):
H^i(\Cal C^3)\to H^{i+1}(\Cal C^1)$ is of determinant class for any
i, then $\Cal C^2$ is
of determinant class.
\newline (ii) If $\Cal C^1$, $\Cal C^2$ are of determinant class and
$H(f_i):
H^i(\Cal C^1)\to H^i(\Cal C^2)$ is of determinant class for any i,
then $\Cal C^3$ is of determinant class.
\newline (iii) If $\Cal C^2$, $\Cal C^3$ are of determinant class and
$H(g_i):
H^i(\Cal C^2)\to H^i(\Cal C^3)$ is of determinant class for any i,
then $\Cal C^1$ is of determinant class.
\endproclaim

\demo{Proof} (i) Recall that
$$\underline d_i^2 = \pmatrix
{\pmatrix d_{i, 36}^2 & 0 \\ d_{i, 46} & d_{i, 47} \endpmatrix} &
{\pmatrix d^2_{i, 38} \\ d^2_{i, 48} \endpmatrix} \\
&&&\\
{\pmatrix 0 \qquad& 0 \endpmatrix}
& d^2_{i, 58} \endpmatrix.$$
{}From Proposition 1.1 we conclude that $\underline d_i^2$ is
of determinant class iff $d^2_{i, 36}, d^2_{i, 47}$ and $d^2_{i, 58}$ are of
determinant class.  In view of Lemma 1.9 it then follows that
$\Cal C^2$
is of determinant class, using that $\underline d_i^1, \underline d_i^3$
and $\gamma_i$ are of determinant class for any $i$.
\newline(ii) By Lemma 1.11(A) (or statement (i) above), the mapping cone
$\Cal C(f)$ is of determinant class. As $\Cal C(f)$
is homotopy equivalent to $\Cal C^3$ (Lemma 1.11(A)) it follows therefore
from
Proposition 1.3(B) that $\Cal C^3$ is of determinant class.
\newline (iii) Applying again Lemma 1.11(A) and (B)
one concludes that the mapping cone $\Cal C(g)$ is of determinant class
and homotopy equivalent to $S\Cal C^1.$ Thus Proposition 1.3(B)
implies that $\Cal C^1$ is of determinant class.
\hfill $\square$
\enddemo

\proclaim{Theorem 1.14}  If three out of the four cochain complexes $\Cal
C^1,
\Cal C^2, \Cal C^3, \Cal H$ are of determinant class then so is the
fourth and one has the following equality
$$\leqalignno{
\log T(\Cal C^2) &=\log T(\Cal C^1) + \log T(\Cal C^3) + \log T(\Cal H)
&(1.75)\cr
&\quad - \sum_i (-1)^i \log T (0 \rightarrow \Cal C_i^1 \rightarrow
\Cal C_i^2 \rightarrow \Cal C_i^3 \rightarrow 0). \cr}$$
\endproclaim

\demo{Proof} In view of Proposition 1.13 it remains to prove, apart
from formula (1.75), that $\Cal H$ is of determinant class if
$\Cal C^1, \Cal C^2$ and $\Cal C^3$ are of determinant class.
Notice that the following statements hold:
\newline (i) If $\Cal C^2$ is of determinant class then
$d^2_{i,36}$ is of determinant class which implies that  $H(\delta_i)$
is of
determinant class for any $i$ (cf Lemma 1.9(A)).
\newline (ii) If $\Cal C^3$ is of determinant class then
$\Cal C(f)$ is of determinant class (Lemma 1.11(B), Proposition 1.3(B)).
Moreover, by Proposition 1.13(i),
the morphisms
$H^i(S\Cal C^1) \to H^{i+1} (\Cal C^2)$ are of determinant class which
is equivalent to $H(f_{i+1})$ being of determinant class (cf (1.70)).
\newline (iii) If $\Cal C^1$ is of determinant class then $\Cal C(g)$ is of
determinant class (Lemma 1.11(B), Proposition 1.3(B))
and, by Proposition 1.13(i), $H(g_i): H^i(\Cal C^2)
\to H^i(\Cal C^3)$ is of determinant class for any $i.$
\newline Combining (i)-(iii) one concludes that
 if $\Cal C^1,\Cal C^2, \Cal C^3$ are of determinant class,
then $\Cal H$  is of determinant class.

To prove
formula (1.75) consider the commutative diagram

$$\leqalignno{
&&(1.76)\cr
0\rightarrow \qquad \Cal C_i^{1, -} \qquad\overset {\left(\smallmatrix 0 \\
\alpha_{i, 74} \\ 0 \endsmallmatrix \right)} \to \rightarrow \qquad
\Cal C_i^{2, -0} \oplus \Cal C_i^{2, -+} \oplus \Cal C_i^{2, --}\qquad
\overset {( 0 \, 0 \, \beta_{i, 48})}
\to \longrightarrow \qquad \Cal C_i^{3, -}
\rightarrow \qquad 0 \cr\cr
\qquad\downarrow \underline {d_i^1} \qquad\qquad\qquad\qquad\qquad \downarrow
\underline {d_i^2} \qquad\qquad\qquad\qquad\qquad\qquad\quad
\downarrow \underline {d_i^3}\qquad\qquad \cr
&&(1.77)\cr
0 \rightarrow \, \Cal C_{i+1}^{1, +} \quad\overset {\left(\smallmatrix 0\\
\alpha_{i+1, 43} \\ 0 \endsmallmatrix \right)} \to \rightarrow \quad
\Cal C_{i+1}^{2,+0} \oplus \Cal C_{i+1}^{2, ++}
\oplus \Cal C_{i+1}^{2,+-} \quad
\overset {( 0 \, 0 \, \beta_{i+1, 35})} \to \longrightarrow
\quad \Cal C_{i+1}^{3, +}
\rightarrow 0.\cr}$$
Notice that $\underline d_i^1, \underline d_i^2$ and $\underline d_i^3$
establish a weak isomorphism of determinant class between the complexes
(1.76) and (1.77).  By Lemma 1.9 both complexes (1.76) and (1.77) are of
determinant class.  By Proposition 1.3, for any $i$

$$\leqalignno{
\log &T \pmatrix 0 \rightarrow & \Cal C_{i+1}^{1, +} \rightarrow
& \Cal C_{i+1}^{2, +}
\rightarrow & \Cal C_{i+1}^{3, +} \rightarrow 0 \endpmatrix &(1.78)\cr
&=\log T \pmatrix 0 \rightarrow & \Cal C_i^{1, -} \rightarrow &
\Cal C_i^{2, -}
\rightarrow & \Cal C_i^{3, -} \rightarrow 0 \endpmatrix\cr
&\quad - \log \Vol (\underline d_i^1) + \log \Vol (\underline d_i^2)
- \log \Vol (\underline d_i^3) - \log \Vol (d^2_{i, 36}),\cr}$$
where we used that $H(\underline {d_i^1}) = 0 = H(\underline d_i^3)$ and
$H(\underline d_i^2) = d^2_{i, 36}$.  Rearrange
terms in (1.78), multiply by $(-1)^{i+1}$
and sum over $i$ to obtain

$$\leqalignno{
\sum_i &(-1)^i \log \Vol (\underline{d_i^1}) - \sum_i (-1)^i \log \Vol
(\underline {d_i^2}) + \sum_i (-1)^i
\log \Vol (\underline {d_i^3})&(1.79)\cr
&= - \sum_i (-1)^i \log T( 0 \rightarrow \Cal C_{i+1}^{1, +} \rightarrow
\Cal C_{i+1}^{2, +} \rightarrow \Cal C_{i+1}^{3, +} \rightarrow 0)\cr
& \quad + \sum_i (-1)^i \log T ( 0 \rightarrow \Cal C_i^{1, -} \rightarrow
\Cal C_i^{2, -} \rightarrow \Cal C_i^{3, -} \rightarrow 0)\cr
& \quad - \sum_i (-1)^i \log \Vol (d^2_{i, 36}).\cr}$$

Recall formula (1.60), $d^2_{i, 36} = \alpha_{i +1, 31} \gamma_i \beta_{i,
26}$.
Using that both $d^2_{i, 36}$ and $\gamma_i$ are weak isomorphisms of
determinant class and, according to Lemma 1.6 which states
 that $\alpha_{i+1, 31}$ and $ \beta_{i, 26}$
are isomorphisms we conclude from (1.4) that

$$\leqalignno{
\log \Vol d^2_{i, 36} = \log \Vol \alpha_{i+1, 31} + \log \Vol \gamma_i +
\log \Vol \beta_{i, 26}.&&(1.80)\cr}$$
Substituting (1.80) into (1.79) and taking into account that

$$
\log T(\Cal H) = \underset i \to \sum (-1)^i
\log \Vol (\alpha_{i, 12}) - \sum (-1)^i
\log \Vol (\beta_{i, 12}) + \sum (-1)^i \log \Vol (\gamma_i),$$
one obtains
$$\leqalignno{
\log &T(\Cal C^1) - \log T(\Cal C^2) + \log T(\Cal C^3)&(1.81)\cr
&= -\log T(\Cal H) - \sum_i (-1)^i \log \Vol \alpha_{i+1, 31} -
\sum_i (-1)^i \log \Vol \beta_{i, 26} \cr
& \quad - \sum_i (-1)^i \log \Vol \beta_{i, 12} + \sum_i (-1)^i
\log \Vol \alpha_{i, 12} \cr
& \quad + \sum_i (-1)^i \log T ( 0 \rightarrow \Cal C_i^{1, -} \rightarrow
\Cal C_i^{2, -} \rightarrow \Cal C_i^{3, -} \rightarrow 0)\cr
& \quad - \sum_i (-1)^i \log T ( 0 \rightarrow \Cal C_{i+1}^{1, +}
\rightarrow
\Cal C_{i+1}^{2, +} \rightarrow \Cal C_{i+1}^{3, +} \rightarrow 0 ). \cr}$$

By Lemma 1.10,

$$\leqalignno{
\log &T( 0 \rightarrow \Cal C_{i+1}^{1, +}
\rightarrow \Cal C_{i+1}^{2, +}
\rightarrow \Cal C_{i+1}^{3, +} \rightarrow 0) &(1.82)\cr
&= \log \Vol (\alpha_{i+1, 43})
- \log \Vol (\beta_{i+1, 35})\cr}$$
and

$$\leqalignno{
\log T( 0 \rightarrow \Cal C_i^{1, -} \rightarrow \Cal C_i^{2, -}
\rightarrow \Cal C_i^{3, -} \rightarrow 0) = \log \Vol \alpha_{i, 74}
- \log \Vol \beta_{i, 48}.&&(1.83)\cr}$$
Changing the index of summation from $i + 1$ to $i$ where necessary, one
obtains

$$\eqalign{
\log &T( \Cal C^2) = \log T(\Cal C^1) + \log T (\Cal C^3) + \log T
( \Cal H)\cr
& -  \sum_i (-1)^i \log \Vol \alpha_{i, 31} - \sum (-1)^i \log \Vol
\alpha_{i, 12}\cr
&- \sum (-1)^i \log \Vol \alpha_{i, 74} - \sum (-1)^i
\log \Vol \alpha_{i, 43} \cr
&+ \sum_i (-1)^i \log \Vol \beta_{i, 12} + \sum_i (-1)^i \log \Vol
\beta_{i, 26}\cr
&+ \sum_i (-1)^i \log \Vol \beta_{i, 35} + \sum_i (-1)^i
\log \Vol \beta_{i, 48} \cr
&= \log T (\Cal C^1) + \log T (\Cal C^3) + \log T (\Cal H) - \sum
(-1)^i \log T (0 \rightarrow \Cal C^1 \rightarrow \Cal C^2 \rightarrow
\Cal C^3 \rightarrow 0)\cr}$$
where, for the last equality, we used again Lemma 1.10.  This proves
(1.75).
\hfill $\square$

\vfill \eject

\proclaim{2.  Torsions for compact manifolds with boundary}\endproclaim

\proclaim{2.1  Reidemeister and analytic torsion}\endproclaim

Let $(M, g)$ be a compact Riemannian manifold with boundary and
\newline $(M, \partial_- M, \partial_+ M)$ a bordism. Further
let $\Cal A$ be a unital von Neumann algebra and
$\Cal E \overset p \to \rightarrow M$ be a bundle of $\Cal A$-Hilbert
modules of finite type
equipped with
a flat connection $\nabla$
with the property that the inner products
$(\cdot, \cdot)$ of $\Cal E_y = p^{-1}
(y)$ are parallel. Such a pair $( \Cal E, \nabla)$ will be denoted by
$\Cal F$ and refered to as a parallel flat bundle of $\Cal A$-Hilbert
modules.

Given a generalized triangulation $\tau = ( h, g')$ of $M$, the collection
of the
unstable manifolds in $M\backslash \partial M$,
$W_y^-$, associated to the critical points
$y \in Cr(h)$, and the flow, induced by $grad_{g'} h,$ provides a
\underbar{relative} CW-complex structure on the pair of spaces
$(S,\partial_{-}M)$,
$S=\partial_-M\cup_{y\in Cr(h)} W_y^-$ whose open cells are given by $W_y^-.$
The inclusion of pairs $(S,\partial_-M)\subset (M,\partial_-M)$ is a
homotopy equivalence.
Notice
that if $\tau = ( h, g')$ is a generalized triangulation for $M= (M,
\partial_- M, \partial_+ M)$, then $\tau_{D} =
(- h, g')$ is a triangulation for $-M$.
Further if $M,$ resp. $N,$ are two closed manifolds  and
$\tau'=(h',g'),$ resp. $\tau''=(h'',g''),$ are generalized
triangulations,
then $\tau_0:=(h'+h'', g'\oplus g'')$
provides a generalized
triangulation on $M\times N.$
The CW-complex structure induced
by $\tau_0$ is the product of the CW-complex structures induced by
$\tau'$ and $\tau''.$ Notice that if $(M,\partial_-M,\partial_+M)$
is a bordism
and $\partial M\ne \emptyset,$ then $h_0:=h'+h''$ is not
constant on the boundary components $\partial_{\pm}M\times N,$
and therefore $\tau_0$ is not a generalized triangulation as defined
in the introduction. However it is
possible to
modify $h_0$ in an arbitrary small
neighborhood of $\partial M\times N,$ so that the new function h has
the same critical points as $h_0,$
$\tau := (h, g'\oplus g'')$ is a
generalized triangulation and the relative
CW-complex structure induced by $\tau$ is the product of
the relative CW complex structure induced by $\tau' $ and
the CW complex structure induced by
$\tau''.$ The generalized
triangulation $\tau$ will be refered to as a product
triangulation $\tau=\tau'\times \tau''.$

Denote by $M = \cup M_j$ the disjoint union of the connected
components $M_j$ of a manifold $M.$ Let $\overset \sim \to M_j$ be
the universal cover of $M_j$ and
introduce $\overset \sim \to M = \cup\overset \sim \to M_j.$
Let $\overset \sim \to h$ be the lift of $h$ to
$\overset \sim \to M,$ and $\overset \sim \to g'$
be the pull back of the Riemannian
metric $g'$ on $\overset \sim \to M.$
Then the gradient vector field $grad_{\overset
\sim \to g'} \overset \sim \to h$ satisfies the Morse-Smale condition as
well.

For an arbitrary critical point $\overset \sim \to y \in Cr(\overset \sim \to
h),$ choose an orientation $\Cal O_{\overset \sim \to y}$ of the descending
manifold $W^{-}_{\overset \sim \to y}$ and let $\Cal O_h := \{
\Cal O_{\overset \sim \to y}; \overset \sim \to y \in Cr(h)\}$.  Given this
choice one constructs, as in [BFKM, section 4], the cochain complex
$\big( \Cal C^q (M, \partial_-M, \tau, \Cal O_h, \Cal F) , \delta_q)$ of
$\Cal
A$-Hilbert
modules of finite type.

Denote by $\Delta_q^{comb}$ the combinatorial Laplacians

$$
\Delta_q^{comb}:= \delta_q^* \delta_q + \delta_{q-1} \delta_{q-1}^*.$$
Here $\delta_q^*$ denotes the adjoint of $\delta_q.$
Observe that $\Delta_q^{comb}$ is a bounded, nonnegative selfadjoint
$\Cal A$-linear operator on $\Cal C^q$ and for any $\epsilon > 0$ one
can define its regularized determinant in the von Neumann sense, by the
following Stieltjes integral

$$
\log \, det_N ( \Delta_q^{comb} + \epsilon ) := \int_{0^+}^\infty \log
( \lambda + \epsilon) d N_{\Delta_q^{comb}} (\lambda)$$
where $N_{\Delta_q^{comb}}(\lambda)$ is the spectral
distribution function associated
to $\Delta_q^{comb}$.

The function $\log \, det_N(\Delta_q^{comb}+ \cdot)$ is
an element in the vector space $\Bbb D$ consisting of equivalence classes
$[f]$ of real analytic functions $f: ( 0, \infty) \rightarrow \Bbb R$ with
$f \sim g$ iff $\underset {\lambda \rightarrow 0} \to \lim \big( f (\lambda)
- g (\lambda) \big) = 0$.  (The elements of $\Bbb D$ represented by the
constant
functions form a subspace of $\Bbb D$ which can be identified with $\Bbb R$,
the space of real numbers.)  We then define $T_{comb} = T_{comb} (M,
\partial_-M,
\tau,
\Cal F)$ as the following element in $\Bbb D$

$$
\log T_{comb} = \frac 1 2 \underset q \to \sum (-1)^{q+1} q \log \, det
(\Delta_q^{comb} + \cdot ).$$
One can show that $\log T_{comb}$ is independent of the choice of
orientations $\Cal O_h$.

To define the analytic torsion we first have to introduce some more notation.

Denote by $\Lambda^q(M; \Cal E)$ the $\Cal A$-module of all smooth $q$-forms
and let $d_q: \Lambda^q(M; \Cal E) \rightarrow \Lambda^{q+1} (M; \Cal E)$
be the exterior differential induced by the exterior differential on
scalar valued $q$-forms and the covariant differentiation given by the
connection $\nabla$ on $\Cal E$.  Then $\big( \Lambda^* (M; \Cal E),
d_*)$ is a cochain complex of $\Cal A$-modules.  The Riemannian metric $g$
together with the connection $\nabla$ induce the Hodge $*$-operators
$J_q: \Lambda^q(M;  \Cal E) \rightarrow \Lambda^{d-q} (M; \Cal E)$, denoted
by $*$.  The operator $J_q$ induces a scalar product
$\langle \cdot, \cdot \rangle$
on $\Lambda^q (M; \Cal E)$, with $\langle \omega', \omega''\rangle$ given by
$\big(
\omega', \omega'' \text { in } \Lambda^q(M; \Cal E) \big)$

$$
\leqalignno{
\langle \omega', \omega''\rangle =
\underset M \to\int \omega ' \wedge_{\Cal E} * \omega''&&(2.1)\cr}$$
where $\wedge \equiv \wedge_{\Cal E}$

is defined by the composition

$$
\Lambda^q (M; \Cal E) \times \Lambda^k (M; \Cal E) \rightarrow \Lambda^{q+k}
( M; \Cal E \otimes \Cal E) \rightarrow \Lambda^{q+k} (M).$$

Denote by $d_{q-1}^* : \Lambda^q (M; \Cal E) \rightarrow \Lambda^{q-1}
(M; \Cal E)$ the formal adjoint of $d_{q-1}$ with respect to this scalar
product.  Notice that $d_{q-1}^*$ is $\Cal A$-linear and

$$
d_{q-1}^* = (-1)^{d(q-1)+1} * d_{d-q} * .$$

Define the submodules $\Lambda_1^q(M; \Cal E)$ and $\Lambda_2^q(M; \Cal E)$
of $\Lambda^q(M; \Cal E)$

$$\leqalignno{
\Lambda_1^q(M; \Cal E) := \{ \omega \in \Lambda^q(M; \Cal E);
\,i^\#_{\partial_-M} \omega = 0 ;\,
i^\#_{\partial_+M} (* \omega) = 0\}
&&(2.2)\cr}$$
and
$$\leqalignno{
\Lambda_2^q(M; \Cal E) &:= \{ \omega \in \Lambda^q(M; \Cal E);
i^\#_{\partial_-M} \omega = 0;\,
i^\#_{\partial_+M} (* \omega)
= 0; &(2.3)\cr
& i^\#_{\partial_-M}
(d_{q-1}^* \omega) = 0; i^\#_{\partial_+M}
\big( d^*_{d-q-1} ( *\omega) \big) = 0
\}\cr}$$
where $i^\#_{\partial_\pm M} : \Lambda^* (M; \Cal E) \rightarrow
\Lambda^* (\partial M_\pm; \Cal E)$ are the pullbacks induced by the
embeddings $i_{\partial_\pm M} : \partial_\pm M \hookrightarrow M$.

Consider the Laplacians acting on $q$-forms in $\Lambda_2^q (M; \Cal E)$,

$$
\Delta_q:= d_q^* d_q + d_{q-1} d^*_{q-1}.$$

The operators $\Delta_q$ are essentially self-adjoint, nonnegative, elliptic
and $\Cal A$-linear.  For any $\epsilon > 0$ one can define its regularized
determinant in the von Neumann sense,

$$
\log \, det_N (\Delta_q + \epsilon ) =
-\frac{d}{ds}\zeta_{\Delta_q+\epsilon}(s)\vert_{s=0}
-\epsilon \log (tr_N P_q)
$$
where $P_q$ denotes the orthogonal projection onto the null space of
$\Delta_q$
and $\zeta_{\Delta_q+\epsilon}(s)$ is the meromorphic continuation of
the function defined for $\Cal Re(s)> d/2$ by the formula (cf. [BFKM])
$$
\zeta_{\Delta_q+\epsilon}(s) = \frac {1}{\Gamma(s)}
\int_0^\infty t^{s-1}e^{-t\epsilon}
tr_N e^{-t\Delta_q}dt.
$$
Then $\log \, det_N (\Delta_q + \cdot)$ is an element in the
vector space $\Bbb D$ introduced above and the analytic torsion
 $T_{an} =
T_{an} (M,\partial_-M, g, \Cal F)$
is defined as the following element in $\Bbb D$

$$\leqalignno{
\log \, T_{an} = \frac 1 2 \sum_q (-1)^{q+1} q \log \, det
(\Delta_q + \cdot ).&&\cr}$$

Denote by $L_2 \big(
\Lambda^q (M; \Cal E) \big )$ the
completion of $\Lambda^q (M; \Cal E)$ with respect to the scalar product
$\langle \cdot, \cdot \rangle$. $L_2 \big(
\Lambda^q (M; \Cal E) \big )$ is an
$\Cal A$-Hilbert module. Notice that the completions of $ \big( \Lambda^q_j
(M; \Cal E)
\big)$ \, $(j = 1, 2)$ give the same space $L_2 \big(
\Lambda^q (M; \Cal E) \big ).$

Introduce the space of harmonic forms,

$$
\Cal H_q (M; \Cal E) := \{ \omega \in \Lambda_2^q(M; \Cal E) ;\, d_q \omega
= 0;\, d^*_{q-1} \omega = 0\}.$$
It is an immediate consequence of the definition of $\Lambda_j^q (M; \Cal E),
j = 1, 2$, that the Hodge $*$ operators induce isometries (Poincar\'e
duality)

$$\eqalignno{
\Lambda_j^q (M; \Cal E) &\rightarrow \Lambda_j^{d-q} (-M; \Cal E)\cr
\Cal H_q (M; \Cal E) &\rightarrow \Cal H_{d-q} (-M; \Cal E).\cr}$$
(Recall that $-M$ is obtained from the bordism $(M, \partial_- M,
\partial_+ M$) by interchanging the role of $\partial_+M$ and $\partial_-M$.)

By the theory of elliptic differential $\Cal A$-operators one concludes that,
for any $q$, $\Cal H_q(M; \Cal E)$ is an $\Cal A$-Hilbert module of finite
von Neumann dimension (in fact it is of finite type).
Denote by $P_q(\lambda) : \Lambda^q (M; \Cal E) \rightarrow
\Lambda_2^q (M; \Cal E)$ the spectral projections corresponding to the
Laplacian $\Delta_q$ and the boundary operator
{\eightpoint
$$\leqalignno{
\quad B_q: \Lambda^q (M; \Cal E)
&\rightarrow \Lambda^q (\partial_- M; \Cal E)
\times  \Lambda^{q-1} (\partial_- M; \Cal E) \times  \Lambda^{d-q}
(\partial_+ M; \Cal E) \times  \Lambda^{d-q-1} (\partial_+ M;
\Cal E)&(2.4)\cr
\omega &\mapsto \big( i^\#_{\partial_-M} (\omega),i^\#_{\partial_-M}
(d^*_{q-1} \omega), i^\#_{\partial_+M} (*\omega), i^\#_{\partial_+M}
\big(d^*_{d-q-1}(*\omega)\big).\cr}$$}

\proclaim{Proposition 2.1}  (Hodge decomposition)  For any $0 \le q \le d$,
there exist orthogonal decompositions of $\Cal A$-pre Hilbert modules

{\eightpoint
$$\leqalignno{
\quad\Lambda^q (M; \Cal E ) = \Cal H_q (M; \Cal E) \oplus
\text { closure }
d_{q-1} \big( \Lambda_1^{q-1} (M; \Cal E) \big) \oplus \text { closure }
d_q^*\big( \Lambda_1^{q+1} (M; \Cal E) \big) &&{(HD)}\cr}$$}
where the word closure refers to the closure
with respect to the usual $C^{\infty}-$topology.
{\eightpoint
$$\leqalignno{
\quad\Lambda_1^q (M; \Cal E ) = \Cal H_q (M; \Cal E) \oplus
\text { closure }
d_{q-1} \big( \Lambda_2^{q-1} (M; \Cal E) \big) \oplus \text { closure }
d_q^*\big( \Lambda_2^{q+1} (M; \Cal E) \big) &&{(HD)_1}\cr}$$}
and

{\eightpoint
$$\leqalignno{
L_2 &\big(\Lambda_1^q(M; \Cal E)\big )
= \Cal H_q (M; \Cal E) \oplus
L_2 \Big(d_{q-1} \big( \Lambda_2^{q-1} (M; \Cal E) \big)\Big) \oplus
L_2 \Big(d_q^* \big( \Lambda_2^{q+1} (M; \Cal E) \big)\Big).
&{\overline {(HD)}}\cr}$$}
\endproclaim

\demo{Proof}  We begin by noticing that
 $\Cal H_q(M; \Cal E),d_{q-1}\big( \Lambda_2^{q-1} (M; \Cal E)
\big)$ and $d_q^* \big( \Lambda_2^{q+1} (M; \Cal E) \big)$ are
$\Cal A$-submodules of $\Lambda_1^q(M; \Cal E)$
and $\Cal H_q(M; \Cal E),\allowmathbreak d_{q-1}
\big( \Lambda_1^{q-1} (M; \Cal E)
\big)$ and
\newline $d_q^* \big( \Lambda_1^{q+1} (M; \Cal E) \big)$ are
$\Cal A$-submodules of $\Lambda^q(M; \Cal E).$
In each of the two cases the three spaces are, due to Stokes' theorem,
pairwise orthogonal with respect to
$\langle\cdot, \cdot \rangle.$ Moreover
$\overline {(HD)}$ follows from $(HD)$ or $(HD)_1.$

Thus it remains to  check that an
arbitrary element $\omega$ in $\Lambda^q(M; \Cal E) \text{ or }
\Lambda_1^q(M; \Cal E)$
can be decomposed as stated. With this in mind, we introduce
for $\epsilon >0$ the operator
$ G_{q; \epsilon}: \Lambda^q(M; \Cal E) \rightarrow \Lambda_2^q
(M; \Cal E),$
which vanishes on the subspace $P_q (\epsilon) \Lambda^q
(M; \Cal E)$ and is equal to $\Delta_q^{-1}$ (considered with boundary
conditions defined by the boundary operator (2.4)) on $\big( Id-P_q
(\epsilon ) \big) \Lambda^q(M; \Cal E).$

Notice that
$\Delta_q G_{q; \epsilon} = Id - P_q(\epsilon ).$
Thus, for $\omega \in \Lambda^q (M; \Cal E)$, one has
$$ \omega = P_q (\epsilon )\omega + d_{q-1} d^*_{q-1} G_{q; \epsilon}
\omega +  d_q^* d_q G_{q; \epsilon} \omega$$
and therefore
$$\leqalignno{
\omega-\big(P_q (\epsilon) - P_q (0)\big) \omega = P_q (0) \omega +
d_{q-1} d_{q-1}^* G_{q; \epsilon} \omega + d_q^* d_q G_{q; \epsilon}
\omega. &&(2.5)\cr}$$
Since  $G_{q; \epsilon} \omega \in \Lambda_2^q(M;\Cal E)$
$$\leqalignno{
d_{q-1}^* G_{q; \epsilon } \omega \in \Lambda_1^{q-1} (M; \Cal E)
\ \text {and}\
d_{q} G_{q; \epsilon} \omega \in \Lambda_1^{q+1} (M; \Cal E).
&&(2.5')\cr}$$
One concludes that (2.5) is a decomposition of type $(HD)$ for
\newline $\omega_{\epsilon}= \omega - ( P_q(\epsilon ) -
P_q (0) ) \omega .$
Let us show that if $\omega\in \Lambda_1^q(M;\Cal E)$ then
$$d_{q-1}^* G_{q; \epsilon } \omega \in \Lambda_2^{q-1} (M; \Cal E)
\ \text {and}\  d_q G_{q; \epsilon } \omega \in \Lambda_2^{q+1} (M; \Cal E)
.$$
In view of (2.3),(2.5') and $d\cdot d=0$ it remains to check that

$$\leqalignno{
i^\#_{\partial_+M} \big(*d_{q-1} d_{q-1}^* G_{q; \epsilon}
\omega \big)=0\ \text {and}\
i^\#_{\partial_-M} \big(d_{q}^* d_q G_{q; \epsilon}
\omega \big)=0.&&(2.6)\cr}
$$

Observe that
$$\leqalignno{
*d_{q-1}d_{q-1}^*
G_{q;\epsilon}\omega =*\big (&\Delta_q-d_q^*d_q\big ) G_{q;\epsilon} \omega=
*\omega-*P_q(\epsilon)\omega \pm d_{d-q-1}\ast d_{q} G_{q;\epsilon}\omega\cr}
$$
and
$$\leqalignno{
d^*_{q}d_q G_{q;\epsilon}\omega = \big (\Delta_q-d_{q-1}
d_{q-1}^*\big ) G_{q;\epsilon} \omega=
\omega-P_q(\epsilon)\omega - d_{q-1}
d^*_{q-1}G_{q;\epsilon}\omega.\cr}
$$
But $i^\#_{\partial_+M} *\omega=0$ and $i^\#_{\partial_-M} \omega=0$
(as $\omega\in \Lambda_1^q(M;\Cal E)),$  $i^\#_{\partial_+M}
*P_q(\epsilon) \omega=0$ and $i^\#_{\partial_-M} P_q(\epsilon )\omega=0$
(as $P_q(\epsilon) \omega\in \Lambda_2^q(M;\Cal E))$ and
$i^\#_{\partial_+M}d_{d-q-1}*d_q G_{q;\epsilon} \omega=
d_{d-q-1} i^\#_{\partial_+M}*d_q G_{q;\epsilon} \omega=0$ and
$i^\#_{\partial_-M}d_{q-1} d^*_{q-1} G_{q;\epsilon} \omega=
d_{q-1} i^\#_{\partial_-M}d_{q-1}^* G_{q;\epsilon} *\omega=0$
(as $G_{q;\epsilon}\omega \in \Lambda_2^q(M;\Cal E)).$

Consequently the
decompositions $(HD)$ and $(HD)_1$ hold for
$\omega_{\epsilon}=\omega -(P_q(\epsilon)- P_q(0))\omega.$ Note that
$\lim_{\epsilon\to 0}\omega_{\epsilon}=\omega $ in
$\Lambda^q(M;\Cal E)$ with respect to the $C^{\infty}-$topology
and all spaces in $(HD)$ and $(HD)_1$ are
closed in this topology.  \hfill $\square$
\medskip

Define
$$\leqalignno{
\Lambda^q(M,\partial_-M; \Cal E) := \{ \omega \in \Lambda^q(M; \Cal E);
\,i^\#_{\partial_-M} \omega = 0\}.
&&(2.7)\cr}$$

Notice that $(\Lambda^q(M, \partial_-M; \Cal E),d_q)$ is a subcomplex.
The integration $Int^{(q)}$ on the $q$-cells of a generalized triangulation
$\tau = (h, g')$, which is given by the unstable manifolds of $grad_{g'} h$,
defines an $\Cal A$-linear map

$$
Int^{(q)}: \Lambda^q (M, \partial_-M; \Cal E) \rightarrow
\Cal C^q (M,\partial_- M, \tau, \Cal O_h;
\Cal F)$$
so that $\delta_q Int^{(q)} = Int^{(q+1)} d_q$.  Denote by $\pi_q$ the
canonical projection $\pi_q: \Cal C^q \rightarrow Null \, (
\Delta_q^{comb})$.
By a theorem of Dodziuk [Do] of de Rham type, generalized for manifolds
with boundary, the map $\pi_q \, Int^{(q)}$, restricted to $\Cal H_q$,
is an isomorphism of Hilbert
modules.  Denote its inverse by $\theta_q$.  Define the metric part
$T_{met} = T_{met} (M,\partial_-M, g, \tau, \Cal F)$
of the Reidemeister torsion by

$$
\log T_{met}:= \frac 1 2 \sum_q (-1)^q \log \, det_N (\theta_q^* \theta_q).$$
$\log T_{met}$  is independent of the choice of the orientations
$\Cal O_h.$
The Reidemeister torsion
$T_{Re} = T_{Re} (M,\partial_-M, g, \tau, \Cal F)$ is defined as an
element in $\Bbb D$ by

$$
\log T_{Re} = \log \, T_{comb} + \log \, T_{met}.$$
In the remaining of this subsection we will discuss the behavior of the
torsions with respect to the Poincar\'e duality and the product formula
for torsions.

Suppose $M$ is oriented. Then the orientations $\Cal O_h$ of the unstable
manifolds corresponding to the gradient vector field
$grad_{g'} h$ induce orientations on the stable manifolds of
$ grad_{g'}h.$ These manifolds can be identified with
the stable manifolds of $grad_{g'}(-h)$ with
 the orientations $\Cal O_{-h}.$
Since the critical points of index $k$ of $h$ are
the critical points of index $d-k$ of
$-h$ one obtains isometries of Hilbert modules
$J_q^{comb}:\Cal C^q(M,\partial_-M,\tau;\Cal F) \to
\Cal C^{d-q}(M,\partial_+M,\tau_{D};\Cal F)$ ( which induce Poincar\'e
duality)
and verifies that

$$\leqalignno{
J_{q+1}^{comb}\cdot \delta_{q;\tau}
= (-1)^{q(d-q)}\delta^*_{d-q-1;\tau_{D}}
\cdot J_q^{comb} &&(2.8)\cr}$$
where $\delta_{q;\tau}$ resp. $\delta_{q;\tau_{D}}$ are the
differentials in the
cochain complexes $\Cal C(M, \partial_-M, \tau; \Cal F)$ resp.
$\Cal C(M, \partial_+M, \tau_{D}; \Cal F)$ and
$\delta^*_{q;\tau}$
denotes the adjoint of $\delta_{q;\tau}.$ Further it is easy to check
that the restrictions of $Int^{(q)} $ and $Int^{(d-q)} ,$
$$
Int^{(q)}: \Cal H_q (M; \Cal E) \rightarrow
\Cal C^q (M,  \partial_-M, \tau, O_h;\Cal F),$$
$$
Int^{(d-q)}: \Cal H_{d-q} (-M; \Cal E)
\rightarrow \Cal C^{d-q} (M, \partial_+M, \tau_{D}, \Cal O_h;
\Cal F)$$
intertwines $J_q:\Cal H_q (M; \Cal E)
\to \Cal H_{d-q} (-M; \Cal E)$ with $J_q^{comb}.$
Let us also recall that $J_{d-q}\cdot J_q = (-1)^{q(d-q)}Id$
and therefore
$$\leqalignno{
J_{d-q} \Delta_{d-q} J_q = (-1)^{q(d-q)}\Delta_q \cdot
J_{q+1}^{comb}\cdot \delta_{q;\tau} =
(-1)^{q(d-q)}\delta^*_{d-q;\tau_{D}}
\cdot J_q^{comb}. &&(2.9)\cr}$$
By using 2.8 and 2.9 and the intertwining mentioned above
one obtains

\proclaim {Proposition 2.2}
The following identities hold in $\Bbb D:$

(i) $\log T_{comb}(M,\partial_-M, \tau ;\Cal F)= (-1)^{d+1}
\log T_{comb}(M,\partial_{+}M, \tau_{D} ;\Cal F);$

(ii) $\log T_{an}(M,\partial_-M, g ;\Cal F)= (-1)^{d+1}
\log T_{an}(M,\partial_+M, g ;\Cal F);$

(iii) $\log T_{met}(M,\partial_-M, g, \tau ;\Cal F)= (-1)^{d+1}
\log T_{met}(M,\partial_+M, g, \tau_{D} ;\Cal F).$
\endproclaim

\proclaim {Proposition 2.3}  Let $(M, \partial_-M,\partial_+ M)$ be a
bordism equipped with a Riemannian metric $g$, a generalized
triangulation $\tau$ and a parallel flat bundle $\Cal F$ of
$\Cal A$-Hilbert modules of finite type. Let $N$ be a closed
manifold equipped with a Riemannian metric $g_N$, a generalized
triangulation $\tau_N$ and a parallel flat bundle $\Cal F_N$ of
$\Cal A_{N}-$Hilbert modules of finite type. Denote by
$\underline \tau$ a product triangulation and by
$\underline {\Cal A}$ the von
Neumann algebra $\underline{\Cal A}=\Cal A\hat{\otimes}\Cal A_N.$
Then $\Cal F\hat{\otimes}\Cal F_N$ is a parallel flat bundle
of $\underline{\Cal A}-$Hilbert modules of finite type.
\newline If the pairs $\{(M,\partial_-M,\partial_+M),\Cal F\}$
and $\{N, \Cal F_N\}$ are of determinant class, so is the pair
$\{(M \times N, \partial_-M\times N,\partial_+M \times N),
\Cal F \hat\otimes \Cal F_N\}$
and the following equalities hold:

$$\leqalignno{
\log T_{an} & ( M\times N,(\partial_-M)\times N,  g\oplus g_N,
\Cal F\hat\otimes \Cal F_N) =
\chi (M,\partial_-M; \Cal F)\cr \log
&T_{an} (N,g;\Cal F_N) +
\chi (N;\Cal F_N) \log T_{an} (M,\partial_-M, g; \Cal F)&(2.10)\cr}$$

$$\leqalignno{
\log T_{Re} & (M\times N,(\partial_-M)\times N, g\oplus g_N,\underline \tau,
\Cal F\hat\otimes \Cal F_N) =
\chi (M,\partial_-M; \Cal F)\cr &\log
T_{Re} (N,g,\tau_N;\Cal F_N) + \chi (N;\Cal F_N)
\log T_{Re} (M,\partial_-M, g,\tau;  \Cal F).&(2.10')\cr}$$
\endproclaim

The proof of Proposition 2.3 is identical to the one of Proposition 4.1 in
[BFKM].

\proclaim {2.2 Determinant class}
\endproclaim
Following [BFKM, section 4.1] we introduce the following

\definition{Definition }
(1) The triple $\{(M, \partial_- M, \partial _+ M), \tau, \Cal F\}$ is said
to
be of
$c$-determinant class iff for $0 \le k \le d$

$$
- \infty < \int^1_{0^+} \log \lambda d N_{\Delta_k^{comb}}(\lambda).$$

(2)  The triple $\{(M, \partial_- M, \partial_+ M), g, \Cal F\}$ is said to
be
of $a$-determinant class iff for $0 \le k \le d$

$$
- \infty < \int^1_{0^+} \log \lambda d N_{\Delta_k}(\lambda).$$
\enddefinition

As in [BFKM], one can conclude from work of Gromov-Shubin [GS] (cf. also
[Ef1],[Ef2]) that the following result holds:

\proclaim{Proposition 2.4} ([GS])
The triple $\{(M, \partial_- M, \partial_+ M), g,\Cal F\}$ is of
$a$-determinant
class iff
$\{(M, \partial_- M, \partial_+ M), \tau, \Cal F\}$ is of
$c$-determinant
class.
\endproclaim

Proposition 2.4 allows us to introduce the following

\definition{Definition } (i) The pair $\{( M, \partial_- M, \partial_+ M),
\Cal F
\}$ is said to be of determinant class if there exists a generalized
triangulation $\tau$ such that $\{(M, \partial_- M, \partial_+ M), \tau,
\Cal F \}$ is of $c$-determinant class.
\newline (ii) The system $( M, \partial_- M, \partial_+ M, \Cal F
)$ is said to be of determinant class if the three pairs
$\{( M, \partial_- M, \partial_+ M), \Cal F\},
\{(M,\emptyset, \partial M), \Cal F\},
\{\partial M,\Cal F \restriction_{\partial M}\}$
are of determinant class.
\enddefinition
Notice that a system $(M,\partial_- M, \partial_+ M, \Cal F)$ being
of
determinant class implies
that the two pairs $\{\partial_{\pm}M, \Cal F\restriction_{\partial_{\pm}M}
\}$ are of determinant class as well as,
 by Poincar\'e duality, that the two pairs
$\{( M, \partial_+ M, \partial_- M), \Cal F\}$ and
$\{(M,\partial M,\emptyset), \Cal F\}$ are of determinant class.

\proclaim {Proposition 2.5 }
(1) If the system $(M, \partial_- M, \partial_+ M, \Cal F)$ is
of determinant class then the
long weakly
exact sequences in cohomology with coefficients in $\Cal F$
associated to $(M,\partial M)$ and
$(M,\partial_{\pm} M)$ are of determinant class.
\newline (2)
Suppose that the systems $( M_i, \partial_- M_i,
\partial_+ M_i, \Cal F_i) \, (i=1,2)$ have the property that
$(\partial_+M_1, \Cal F_1\restriction_{\partial _+ M_{1}})
= (\partial_-M_2, \Cal F_2 \restriction_{\partial _- M_{2}}).$
Let $( M, \partial_- M,
\partial_+ M, \Cal F)$ be the system where $M$ is obtained by glueing
$M_1$ to $M_2$ along $\partial_+M_1=\partial_-M_2$ and $\Cal F$
defined by $\Cal F_i$ on $M_i$.
If the systems
$( M_i, \partial_- M_i,
\partial_+ M_i, \Cal F_i)$ are of determinant class then so is the
system
$( M, \partial_- M,
\partial_+ M, \Cal F)$.
\newline (3) Suppose that each connected component of
a closed Riemannian manifold $N$ is simply connected
and let $\tilde \Cal F$ denote the pullback of
$\Cal F$ by the canonical projection of $M\times N$ on $M.$
Then the pair
$\{(M,\partial_-M,\partial_+ M),\Cal F\}$ is of determinant class iff
the pair $\{(M\times N,\partial_-M\times N,\partial_+M\times N),
\tilde \Cal F\}$ is.
\newline (4) The system $(M,\partial_-M, \partial_+M , \Cal F)$ is
of determinant class iff the system
\newline $(M,\emptyset,\partial M, \Cal F)$
is of determinant class.
\endproclaim

\demo{Proof} (1) We verify the statement for
$(M,\partial_-M).$ Similar arguments can be used for $(M,\partial M)$ and
$(M,\partial_+M).$
Choose generalized triangulations
$\tau_-=(h_-,g'_-)$ of $\partial_-M$ and
$\tau=(h,g')$ of $(M,  \partial_-M,\partial_+ M).$
The
descending manifolds of $grad_{g'_-}h_-$ in $\partial_-M$ and of
$grad_{g'}h$ in $ M\backslash \partial_-M$ provide a
cell structure for the space $X=S\cup \partial_-M$,
where $S$ is the union of all open cells of $\tau$ (i.e. of all
descending manifolds of
$grad_{g'}h$ in $M\backslash\partial_-M$). By an arbitrary small
change of the metric $g'_-$ which keeps the CW-complex structure
of $\partial_-M$ isomorphic to the unperturbed one,
one can arrange that the cell structure
is in fact a CW-complex structure, denoted by
$\tilde \tau .$
Choose orientations
$\Cal O_h$ resp. $\Cal O_{h_-}$ for the cells of $\tau$ resp.
$\tau_-.$
These orientations define orientations $\Cal O$ for the
cells of $\tilde \tau$. Consider the cochain complex
$ (\Cal C^q(X,\tilde {\tau},\Cal O; \Cal F), \delta_q)$ of $\Cal A$
-Hilbert
modules of finite type associated to the CW-complex
$(X,\tilde \tau)$ and $\Cal F$. The short exact sequence of cochain
complexes
$$
0\to \Cal C^*(M,\partial_-M, \tau, \Cal O_h; \Cal F)
\to \Cal C^*(X,\tilde \tau,\Cal O; \Cal F)
\to
\Cal C^*(\partial_-M, \tau_-, \Cal O_{h_-}; \Cal F \restriction_{\partial_-M}
)\to 0
$$
induces a long cohomology sequence  which is weakly exact,
hence a cochain complex. This complex will be called
the cohomology sequence of the pair $(M, \partial_-M).$
Different generalized triangulations give rise to isomorphic (but not
necessarily isometric) cochain complexes.
Using subdivisions, one can show that any
generalized triangulation $\tau'=(h',g'')$ of the bordism
$(M, \emptyset,\partial M)$ produces a cochain complex
$\Cal C^*(M, \tau',\Cal O_{h'};\Cal F)$ of
$\Cal A$-Hilbert modules which is homotopy equivalent to
$\Cal C^*(X,\tilde \tau,\Cal O; \Cal F).$
By hypotheses,
$\Cal C^*(M, \tau',\Cal O_{h'};\Cal F)$ is of determinant class
and therefore, by Proposition 1.3(B),
$\Cal C^*(X,\tilde \tau, \Cal O; \Cal F)$ is of determinant class as well.
Thus, by Theorem 1.14, the cohomology sequence of the pair
$(M,\partial_-M)$ is of determinant class.

(2) We verify that the pair
$\{(M,\emptyset,\partial M),\Cal F\}$ is of determinant class.
Similar arguments can be applied to the pair
$\{(M,\partial_-M,\partial_+ M), \Cal F\}.$
(Notice that, by assumption, $\{\partial M,
\Cal F \restriction_{\partial M}\}$ is of determinant class.)
Choose generalized triangulations $\tau_1=(h_1,g'_1)$ for
$(M_1,\partial_+M_1,\partial_-M_1)$, $\tau_2=(h_2,g'_2)$ for
$(M_2,\partial_-M_2,\partial_+M_2)$  and $\tau_0 =(h_0,g'_0)$ for $\partial_+
M_1= \partial_-M_2.$
The open cells of $\tau_0,\tau_1,\tau_2$ provide a
CW-complex structure $\tilde \tau $
(if necessary, modify $g_0'$ slightly) for the space
$ Y=\partial_+M_1\cup S_1\cup S_2$ where $S_1$ resp. $S_2$ are the
union of the open cells of $\tau_1$ resp. $\tau_2.$
Choose orientations $\Cal O_{h_0}, \Cal O_{h_1}, \Cal O_{h_2}$ for
the cells of $\tau_0, \tau_1, \tau_2$
which define orientations $\Cal O$ for the cells of $\tilde \tau,$
and consider the associated cochain complex
$(\Cal C^q(Y,\tilde \tau,\Cal O; \Cal F), \delta_q)$  of $\Cal A-$Hilbert
modules of finite type.
Consider the short exact sequence of cochain complexes
$$
\leqalignno{
0\to \Cal C^q(M_1,\partial_+M_{1},\tau_1,\Cal O_{h_1};\Cal F_1)\oplus
&\Cal C^q(M_2,\partial_-M_{2}, \tau_2,\Cal O_{h_2};\Cal F_2)
\to \Cal C^q(Y,\tilde  \tau,\Cal O; \Cal F)
\to \cr \to &\Cal C^q (\partial_+M_1,\tau_0,\Cal O_{h_0};
\Cal F\restriction_{\partial_+M_{1}})\to 0.
\cr}$$
Notice that the connecting homomorphisms in the induced
cohomology sequence, $H^q(\partial_+M_1,\tau_0;
\Cal F\restriction_{\partial _+M_{1}})
\to H^{q+1}(M_1,\partial_+M_1,\tau_1; \Cal F_1)\oplus
H^{q+1}(M_2,\partial_-M_2, \tau_2; \Cal F_2),$
are the sum of the connecting homomorphisms in the cohomology
sequences
of
the pairs
$(M_1,\partial_+M_1)$ and $(M_2,\partial_-M_2)$ and hence of determinant
class by (1) and Proposition 1.1 (C).
Therefore, Proposition 1.13 (i) together with the made assumptions implies
that  the complex $\Cal C^*(Y,\tilde
\tau,\Cal O; \Cal F)$ is of determinant class.
By Propostion 1.3(B), the homotopy
equivalence of $\Cal C^*(Y,\tilde
\tau,\Cal O; \Cal F)$
with the cochain complex associated to an arbitrary generalized triangulation
of the bordism $(M,\emptyset,\partial M)$ implies
that the pair
$\{ (M,\partial_- M,\partial_+ M),\Cal F\}$ is of determinant class.
\newline (3) follows from Proposition 1.2.
\newline (4) In view of the definition of a system being of determinant
class it suffices to prove that the pair
 $\{(M,\partial_-M,\partial_+M);\Cal F\}$ is of determinant class if
the system $(M,\emptyset ,\partial M, \Cal F)$ is of determinant
class. First notice the following two facts:
\newline (i) The disjoint union of two systems of determinant class
is a system of determinant class and
\newline (ii) if $\{N,\Cal F\}$ is a pair of determinant class ($N$
a closed manifold) and $\Cal F_{[a,b]}$ denotes the pullback of
$\Cal F$ by the projection $p:N\times [a,b]\to N $ ($[a,b]$ a compact
interval on $\Bbb R$) then
$(N\times [a,b], N\times a, N\times b ;\Cal F_{[a,b]})$
is a system of determinant  class.

For $i=1,2$ define bordisms $(N_i,\partial_-N_i,\partial_+N_i)$
$$
(N_1,\partial_-N_1, \partial _+N_1)
:= (M,\emptyset,\partial M)\, \sqcup \,
(\partial_-M\times [-1,0], \partial_-M\times \{-1\},\partial_-M\times \{0\})
$$
$$
(N_2,\partial_-N_2, \partial_+N_2):=
(\partial_-M\times [0,1],\partial_-M \times \{0,1\},\emptyset)$$
$$\sqcup \, (\partial_+M\times [0,1],
\partial_+M\times \{0\},\partial_+M\times \{1\})$$
and whose parallel flat bundle $\Cal G_i$
of $\Cal A-$Hilbert modules are equal to $\Cal F$  on $M$ and to
$({\Cal F}\restriction_{\partial_{\pm}M})_{[a,b]}$
on the other components. (Recall that $A\sqcup B$ denotes the
disjoint union of two sets A and B.)
\newline The statements (i) and (ii) above imply that the
systems $(N_i,\partial_-N_i,\partial_+N_i; \Cal G_i)$ are of determinant
class. Since $(\partial_+N_1, {\Cal G_1}\restriction_{\partial_+N_1})=
(\partial_-N_2, {\Cal G_2}\restriction_{\partial_-N_2})$ one can glue
these two systems to obtain a system $(N,\partial_-N,\partial_+N,\Cal G)$
which by (2) is of determinant class. The result follows once we
observe that  $(N,\partial_-N,\partial_+N,\Cal G)$ is isomorphic
to $(M,\partial_-M,\partial_+M,\Cal F).$
\hfill $\square$
\enddemo

\proclaim{2.3 Witten deformation of the de Rham complex for manifolds
with boundary}\endproclaim

Given a triangulation $\tau = (h, g')$ of $M$, consider the Witten
deformation
$d_q(t)$ of $d_q$,

$$\leqalignno{
d_q(t) = e^{-th} d_q e^{th} = d_q + tdh \wedge.&&(2.8)\cr}$$

Noticing that $d^*_{q-1} (t) = (-1)^{d(q-1) +1} * d_{d-q} (-t) *$
we define the
spaces
\newline $\Lambda^*_{j;t} (M; \Cal E) \equiv \Lambda_j^* (M; \Cal E) (t) \,
(j = 1, 2)$ and $\Cal H_{*; t} (M; \Cal E)$ by replacing in the definition
(2.1) - (2.2) the exterior derivative $d_q$ and its adjoint $d_q^*$ with
$d_q(t)$ respectively $d_q^*(t)$.

Using the same proof, the above Proposition 2.1 can be generalized as
follows:

\proclaim{Proposition 2.1'}  For any $t \in \Bbb R$ and $0 \le q \le d$,
there
exists an orthogonal decomposition of $\Cal A$-pre Hilbert modules

{\eightpoint
$$\leqalignno{
\quad\Lambda^q &(M; \Cal E)= \Cal H_{q; t} (M; \Cal E) \oplus\,
closure \, d_{q-1}(t) \big( \Lambda_{1;t}^{q-1} (M; \Cal E) \big) \oplus
\, closure \,
d_q^* (t) \big( \Lambda_{1; t}^{q+1} (M; \Cal E)
\big),&{(HD)^t}\cr}$$}

{\eightpoint
$$\leqalignno{
\quad\Lambda_1^q &(M; \Cal E)= \Cal H_{q; t} (M; \Cal E) \oplus\,
closure \, d_{q-1}(t) \big( \Lambda_{2;t}^{q-1} (M; \Cal E) \big) \oplus
\, closure \,
d_q^* (t) \big( \Lambda_{2; t}^{q+1} (M; \Cal E) \big)&(HD)_{1}^{t}\cr}$$}

(where the word closure refers to closure with respect to the
$C^{\infty}-$topology) and

{\eightpoint
$$\leqalignno{
\qquad L_2 &\big(\Lambda^q (M; \Cal E)\big)
= \Cal H_{q; t} (M; \Cal E) \oplus
L_2 \Big( d_{q-1}(t) \big( \Lambda^{q-1}_{2; t} (M; \Cal E) \big) \Big)
\oplus L_2 \Big (d_q^* (t) \big( \Lambda_{2; t}^{q+1} (M; \Cal E) \big)
\Big).&{(\overline{HD})^t}\cr}$$}
\endproclaim

Consider the differential boundary operator
{\eightpoint
$$\leqalignno{
\quad B_q(t): \Lambda^q (M; \Cal E)
&\rightarrow \Lambda^q (\partial_- M; \Cal E)
\times  \Lambda^{q-1} (\partial_- M; \Cal E) \times  \Lambda^{d-q}
(\partial_+ M; \Cal E) \times  \Lambda^{d-q-1} (\partial_+ M;
\Cal E)&(2.4)\cr
\omega &\mapsto \big( i^\#_{\partial_-M} (\omega),i^\#_{\partial_-M}
(d^*_{q-1}(t) \omega), i^\#_{\partial_+M} (*\omega), i^\#_{\partial_+M}
\big(d^*_{d-q-1}(t)(*\omega))\big).\cr}$$}

One verifies that, for any $t \in \Bbb R, \big( \Delta_q (t), B_q (t) \big)$
is an elliptic boundary value problem which is essentially selfadjoint and
nonnegative.  Hence the spectrum $\sigma_q (t):= Spec \, \big( \Delta_q
(t), B_q(t) \big)$ of $\big( \Delta_q (t), B_q (t) \big)$ is contained in
$[0, \infty )$.  Following the proof of Proposition 5.2 in [BFKM] one
concludes that there exist constants $C>0, C'>0$ and $t_0 >0$ so that for
$t \ge t_0$ and $0 \le q \le d,$

$$
\sigma_q (t) \cap (e^{-tC} , C't) = \emptyset .$$

Choose $t_0$ so large that

$$
e^{-t_0C} < 1 < C't_0$$
and introduce, for $t \ge t_0$, the orthogonal projections $Q_q (t):
\Lambda^q (M; \Cal E) \rightarrow \Lambda^q_{2; t} (M; \Cal E),$

$$
Q_q(t) := \frac 1 {2 \pi i} \int_{S^1} \big( \lambda - \Delta_q (t)
\big)^{-1}
d \lambda$$
where $\big ( \lambda - \Delta_q (t) \big)^{-1}$ is viewed as a map

$$
\big( \lambda - \Delta_q (t) \big)^{-1} : \Lambda^q (M; \Cal E)
\rightarrow \Lambda_{2; t}^q (M; \Cal E)$$
and $S^1$ is the circle in the complex $\lambda$-plane, centered at $0$,
with radius 1.

Introduce the following $\Cal A$-submodules of $\Lambda_{1; t}^q (M; \Cal
E)$,

$$
\Lambda_t^q (M; \Cal E)_{sm} := Q_q (t) \big( \Lambda^q (M; \Cal E) \big)$$
and

$$
\Lambda_t^q (M; \Cal E)_{la} := \big(Id-Q_q (t)\big) \big( \Lambda^q
(M; \Cal E) \big).$$

Notice that $\Lambda_t^q (M; \Cal E)_{sm} \subseteq \Lambda_{2;t}^q
(M; \Cal E)$ is a $\Cal A$-Hilbert submodule of finite von Neumann dimension.
Further $\big( \Lambda_t^q (M; \Cal E)_{sm}, d_q (t) \big)$ and $\big(
\Lambda_t^q (M; \Cal E)_{la}, d_q (t) \big)$ is a subcomplex of
$\big( \Lambda^q (M, \partial_-M; \Cal E), d_q (t) \big)$.

Let us now describe $\big( \Lambda_t^q (M; \Cal E)_{sm}, d_q (t) \big)$ in
more detail.

Recall that the cochain complex $\big ( \Cal C^q =
\Cal C^q (M, \partial_-M, \tau,
\Cal O_h; \Cal F), \delta_q \big)$ is given by

$$
\Cal C^q = \underset {x \in Cr_q (h)} \to \oplus \Cal E_x$$
where $Cr_q (h) = \{ x_{q; j} ; 1 \le j \le m_q\}$ denotes the set of
critical points of $h$ which are of index $q$.  Given $q$ and $1 \le j \le
m_q
:= \# Cr_q (h)$, denote by $e_{q; j, i} (1 \le i \le \ell)$ an orthonormal
basis of $\Cal E_{x_{q; j}} \,\big( x_{q; j} \in Cr_q (h) \big)$.
(Recall that $\Cal E_{x}$ is assumed to be a free $\Cal A-$Hilbert
module (cf introduction).) Define
$E_{q; j, i} \in \Cal C^q$ by

$$
E_{q; j, i} ( x_{q; j'}) = \delta_{j'j} e_{q; j, i}$$
where $\delta_{j'j}$ denotes the Kronecker delta.  With respect to this
basis, $\delta_q$ is given by

$$
\delta_q (E_{q; i, j}) = \underset {1 \le j' \le \ell} \to {\underset
{1 \le i' \le m_{q+1}} \to \sum} \gamma_{q; ij, i'j'} E_{q+1; i'j'}.$$
Generalizing results of Helffer-Sj\"ostrand, one constructs for $t \ge t_0$
an orthonormal base $\varphi_{q; i,j} (t)$ of $\Lambda_t^q(M; \Cal E)_{sm}$
to conclude that $\Lambda_t^q (M; \Cal E)_{sm}$ is a free $\Cal A$-Hilbert
module of finite type.

Let us recall from [BFKM] the notion of a $H$-neighborhood $U_{qj} \equiv
U_{x_{q; j}}$ of a critical point $x_{q; j} \in Cr_q (h)$:

\definition{Definition }  $U_{qj} \subseteq M \backslash \partial M$ is
said to be an $H$-neighborhood of $x_{q; j}$ if there exist a disc
$B_{2 \alpha}
:= \{ x \in \Bbb R^d ; |x| < 2 \alpha \}$ and diffeomorphisms $\varphi :
B_{2 \alpha} \rightarrow U_{qj}$ and $\Phi : B_{2 \alpha} \times
\Cal W \rightarrow \Cal E\restriction_{U_{qj}}$ with the following
properties:

(i) $\varphi (0) = x_{q;j}$;

(ii) when expressed in the coordinates induced by $\varphi,$ $h$ is of the
form

$$
h (x) = h (x_{q;j}) - \frac 1 2 \sum_{k = 1}^q x_k^2 + \frac 1 2
\sum_{k = q+1}^d x_k^2;$$

(iii) the pullback $\varphi^*(g')$ of the Riemannian metric $g'$ is the
Euclidean metric;

(iv) $\Phi$ is a trivialization of $\Cal E\restriction_{U_{qj}}$.
\enddefinition

\noindent For later use, we call the coordinates provided by $\varphi$
$H$-coordinates and define $U'_{qj} := \varphi (B_\alpha)$.

A collection $(U_x)_{x \in Cr(h)}$ of $H$-neighborhoods is called a system
of $H$-neighborhoods if, in addition, they are pairwise disjoint.

Taking into account that $h$ is not necessarily a self indexing Morse
function one obtains the following version of Theorem 5.7 in [BFKM].  (We
recall that we assume throughout this section that the $\Cal W_j$'s are free
$\Cal A$-Hilbert modules.)

\proclaim {Theorem 2.6}  Assume $(U_x)_{x \in Cr(h)}$ is a system of
$H$-neighborhoods such that, for any $q$ and $j \neq j'$

$$
U_{x_{q; j}} \cap W^-_{q; j'} = \emptyset $$
where $W^-_{q;j'}$ denotes the descending manifold associated to the
critical point $x_{q;j'}$ and the gradient flow $grad_{g'} h.$
\newline Then there exists $t_1 \ge t_0$ such that for $t > t_1$, the
elements
$\varphi_{q; j, i} (t)$, constructed in [BFKM, (5.31)] with $1 \le j \le
m_q, 1 \le i \le \ell$ form an orthonormal basis of the $\Cal A$-Hilbert
module $\Lambda_t^q (M; \Cal E)_{sm}$ with the following properties:

(i)  There exist $C>0, \eta >0$ so that for $t \ge t_1, 1 \le r \le \ell,$

$$
\underset {x \in M \backslash U_{x_{q,j}}} \to \sup \Vert \varphi_{q; j,r}
(t) (x) \Vert \le C e^{-t \eta} .$$

(ii)  When expressed in $H$-coordinates on $U_{x_{q; j}} \cap W^-_{q,
j}, $
the $q$-forms $\varphi_{q; j, r} (t)$ satisfy the following estimate

$$
\varphi_{q; j, r} (t) (x) = \left ( \frac t \pi \right )^{d/4} e^{-t|x|^2/2}
\left( dx_1 \wedge \dots \wedge dx_q \otimes e_{q; j, r} + 0 \big (\frac 1 t
\big)\right).$$

(iii)  With respect to this basis,

$$
d_q(t) \varphi_{q; j, r} (t) = \underset {1 \le r' \le \ell} \to
{\underset {1 \le j' \le m_{q+1}} \to \sum} \eta_{q; jr, j' r'} (t)
\varphi_{q+1; j', r'} (t)$$
and the coefficients $\eta_{q; jr, j' r'}$ satisfy

$$
\eta_{q; jr, j' r'} (t) = \left ( \big ( \frac t \pi \big )^{1/2}
\gamma_{q; jr, j' r'} + 0 (1) \right) e^{-t \big ( h (x_{q+1, i'}) - h
(x_{q; i}) \big)}.$$
\endproclaim

We need the following application of the above results (cf. [BZ]):

\proclaim{Corollary 2.7}

$$
Int^{(q)} \big( e^{th} \varphi_{q; j, r} (t) \big) = \big( \frac t \pi
\big)^{\frac {d-2q} 4} e^{th(x_{q; j})} \big( E_{q; j, r} + 0 \big( \frac 1
t)
\big).$$
\endproclaim

\demo{Proof} As in [BFKM] we must show that for any cell $W^-_{q; j'}$

$$
\int_{W^-_{q; j'}} \varphi_{q; j, r} (t) e^{ht} = \big( \frac t \pi
\big)^{(d - 2q)/4} e^{th(x_{q; j})} \big( \delta_{j j'} e_{q; j, r}
+ 0 ( \frac 1 t) \big).$$

Estimates of $\varphi_{q; j, r}$ in terms of the Agmon distance ([HS],
[BZ], [BFKM])imply that it
suffices to consider the case where $j = j'$.  Moreover it suffices to
estimate

$$
\int_{W^-_{q; j} \cap U_{qj}} \varphi_{q; j, r} (t) e^{ht}.$$

Note that on $W^-_{q;j} \cap U_{qj}$, the function $e^{ht}$ is of the
form

$$
e^{ht} = e^{h (x_{q; j})t} e^{-t ( \underset 1 \to {\overset q \to\sum}
x_k^2 /2)}$$
and the statement follows therefore from Theorem 2.6(ii).
\hfill $\square$
\enddemo
Consider the isomorphism of $\Cal A$-Hilbert modules

$$
f_k (t) : \Lambda_t^k (M; \Cal E)_{sm} \rightarrow
\Cal C^k (M,\partial_-M,  \tau, \Cal O_h; \Cal F)$$
given by

$$\leqalignno{
f_k (t) \big ( \varphi_{k; jr} (t) \big) := \Big( (\frac \pi t )^{\frac
{d-2k} 4} e^{-th(x_{q; j})} \Big) Int^{(k)}
( e^{th}\varphi_{k;jr}(t)).&&(2.9)\cr}$$
Then the maps $f_k (t)$ provide an isomorphism between
$ \big ( \Lambda_t^k (M; \Cal E)_{sm}, \overset \sim \to d_k (t) \big)$
 and
$\big( \Cal C^k (M,\partial_-M, \tau, \Cal O_h; \Cal F), \delta_k \big)$
where $\overset \sim \to d_k (t)$, when expressed with respect to the
basis $\big( \varphi_{k; ir} (t) \big)$, is given by the matrix

$$
\Big( (\frac \pi t)^{1/2} \cdot e^{t \big( h (x_{q+1; i'}) - h(x_{q;i})
\big)} \cdot \eta_{q; ir, i' r'} (t) \Big)_{ir, i' r'}.$$
As a consequence of Corollary 2.7, $f_k (t)$, when expressed with respect to
the basis $(\varphi_{k; ir})_{i, r}$ of $\Lambda_t^k (M; \Cal E)_{sm}$ and
the
basis $(E_{k; ir})_{ir}$ of
$\Cal C^k (M,\partial_-M, \tau, \Cal O_h; \Cal F),$ is
of the form

$$f_k (t) = Id + 0 ( \frac 1 t ) \quad (t \rightarrow \infty ).$$

\proclaim{2.4  Asymptotic expansions and a comparison theorem}\endproclaim

Suppose that $(M,\partial_-M,\partial_+ M )$ is a bordism, $g$ a Riemannian
metric, $\tau=(h,g')$ a generalized triangulation and $\Cal F $ a
parallel flat bundle.
\newline Denote by
$\log T(M,\partial_-M, g, h, \Cal F) (t)$
the analytic torsion defined
by the Laplace
operators $\Delta_q (t)$ associated to $d_q (t)$ and the metric $g$.
For $t$ large enough
let $\log T_{sm} (M,\partial_-M, g, h, \Cal F) (t)$ and
$\log T_{la} (M,\partial_-M, g, h, \Cal F) (t)$
be the torsion of the subcomplex $\big( \Lambda_t^q (M; \Cal E)_{sm}, d_q
(t) \big)$ respectively $\big( \Lambda_t^q (M; \Cal E)_{la}, d_q (t) \big).$
\newline Both,
\newline $\log T (M,\partial_-M, g, h, \Cal F) (t)$
and $\log T_{sm} (M,\partial_-M, g, h, \Cal F) (t)$ are
elements in the vector space $\Bbb D$ and
$\log T_{la} (M,\partial_-M, g, h, \Cal F) (t)$
is a real number. In the case the pair
$\{( M,\partial_- M, \partial_+ M), \Cal F \}$ is of determinant class
the Laplace operators $\Delta_q(t)$ associated to
$d_q(t)$ and the metric $g$ are of determinant class.
As in [BFKM] one proves the following two theorems:

\proclaim {Theorem A} Let $M= (M, \partial_- M, \partial_+ M)$ be a
bordism, $g$ a Riemannian metric, $\tau =
(h, g)$ a generalized triangulation and $\Cal F $
a parallel flat bundle of $\Cal A-$Hilbert modules
so that the pair $\{(M, \partial_- M, \partial_+ M),
\Cal F \}$ is of determinant class.
Then the following statements hold:

(i) The functions $\log T (M,\partial_-M, g, h, \Cal F) (t)$,
$\log T_{sm} (M,\partial_-M, g, h, \Cal F)
(t)$ and
\newline $\log T_{la} (M,\partial_-M, g, h, \Cal F) (t)$
admit asymptotic expansions for
$t \rightarrow \infty$ of the form

$$
a_0 + \sum_{j=1}^{d+1} a_j t^j + b \log t + o(1).$$

(ii)  The asymptotic expansion of
$\log T (M,\partial_-M, g, h, \Cal F)(t)$ is of the form

$$\eqalign{
\log T &(M,\partial_-M, g, h, \Cal F)(0) -
\log T_{met} (M,\partial_-M, \tau, g, \Cal F)\cr
& + ( \log \pi) \big( \sum_{q = 0}^d (-1)^q \frac {d - 2q} 4 \beta_q \big)
+ \big( \sum_{q = 0}^d (-1)^{q+1} \frac {d-2q} 4 \beta_q \big) \log t \cr
& + \big( \sum_{q = 0}^d (-1)^{q+1} q \beta_q \big) t + \sum_{j=1}^{d+1}
\big( \sum_{q = 0}^d (-1)^q p_{q, j} \big) t^j + o(1)\cr}$$
where $p_{q, j}$ can be written as a sum, $p_{q, j} = p_{q, j}^I +
p_{q, j}^{II}$ with $p_{q, j}^I$ being a local term on $M$ and
$p_{q, j}^{II}$  being a local term on $\partial M$.

(iii)  The asymptotic expansion of
$\log T_{sm} (M,\partial_-M, g, h, \Cal F)(t)$ is
of the form

$$
\log T_{comb}(M,\partial_-M, \tau, \Cal F) + \frac 1 2 \sum_{q = 0}^d (-1)^q
q(m_q \ell - \beta_q) (2t - \log t + \log \pi)+ o (1).$$
\endproclaim

\proclaim{Theorem B} Assume that the pairs
$\{(M_j, \partial_- M_j, \partial_+ M_j),
\Cal F_j\}\, (j=1,2)$ are of determinant class and $\tau_j = (h_j, g_j)$
are generalized triangulations
so that there exist neighborhoods $\Cal U_j$ of $\partial M_j \cup Cr(h_j)$
in
$M_j$ and a diffeomorphism $\Psi: \Cal U_1 \rightarrow \Cal U_2$ with $\Psi
( \partial_\pm M_1) = \partial_\pm M_2$,
$\Psi \big( Cr_q (h_1) \big) =  Cr_q (h_2) $ for $0 \le q \le
d$, as well as $\Psi^* g_2 = g_1$ and $h_2 \circ \Psi = h_1$ on $\Cal U_1$.

Then the free term

$$
FT \big ( \log T_{la} (M_1,\partial_-M_1, g_1, h_1, \Cal F_1) (t) -
\log T_{la} (M_2,\partial_-M_2, g_2,  h_2, \Cal F_2) (t) \big)$$
of the asymptotic expansion of
$$\log T_{la} (M_1,\partial_-M_1, g_1, h_1, \Cal F_1) (t) -
\log T_{la} (M_2,\partial_-M_2, g_2, h_2, \Cal F_2) (t)$$
for $t \rightarrow \infty$ is given by

$$
\int_{M_1 \backslash Cr(h_1)} a_0 (h_1, \epsilon = 0, x_1) -
\int_{M_2 \backslash Cr(h_2)} a_0 (h_2, \epsilon = 0, x_2)$$
where $a_0 (h_1, \epsilon, x_1)$ and $a_0 (h_2, \epsilon, x_2)$ are densities
(forms of degree $d$) and are given by explicit local formulas;  the
difference of the two integrals has to be taken in the same way as in
[BFK, section 0, Remarks after Theorem 3].
\endproclaim
\medskip
Introduce

$$\eqalign{
A &( M_1, M_2, \tau_1, \tau_2, \Cal F_1,
\Cal F_2) \cr
&:= \int_{M_1 \backslash Cr(h_1)} \, a_0 (h_1; \epsilon = 0, x_1) -
\int_{M_2 \backslash Cr(h_2)} \, a_0 (h_2; \epsilon = 0, x_2).\cr}$$

\proclaim{Corollary C} Assume in addition to the hypotheses
in Theorem B that there exists a system
$(M_3,\partial_-M_3,\partial_+M_3, \Cal F_3)$
with the following properties $(j=1,2)$:
\newline (i) $\partial_{\pm}M_3=\partial_{\pm}M_j$ and
$\Cal F_3 \restriction_{ \partial_{\pm}M_3} =
\Cal F_j \restriction _{\partial_{\pm}M_j};$
\newline (ii) the pair $\{N_j,\tilde \Cal F_j\}$
is of determinant class where $N_j$ denotes the closed
manifold
obtained by gluing $M_j$
to $M_3$ and $\tilde \Cal F_j $ is given by $\Cal F_j $ on $M_j$
and by $\Cal F_3$ on $M_3.$
\newline Then
$$
A(M_1, M_2, \tau_1, \tau_2, \Cal F_1, \Cal F_2)=0$$
\endproclaim

\demo{Proof} Let a be the minimum of $h_1$ (which is also the minimum
of $h_2$) and let b be its maximum.
Choose a smooth function $h_3:M_3 \to \Bbb R$ with
$h_3(M_3) \subset [a-1,b+1]$
as well as $ h_3\restriction _{\partial_-M_3}=a,
h_3\restriction _{\partial_+M_3}=b$  (cf. Figure 1), and a metric $g_3$
on $M_3$
so that, for $i=1,2,$ $h_i$ together with $h_3$ defines
a Morse function $\tilde h_i$
on $N_i$ and
$g_i$ together with $g_3$ defines a metric $\tilde g_i$ on $N_i$.
Further choose $g_3$ (if necessary,
 modify  $g_3$ in an arbitrary small neighborhood of
$\partial M_3)$ so that $(\tilde h_i, \tilde g_i)$ is a generalized
triangulation for $N_i$. Since $\{N_i,\Cal F_i\}, (i=1,2)$ are, by
assumption,
of determinant class, one concludes from Theorem A
$$\eqalign{
FT &\big( \log T (N_1,\tilde g_1,\tilde h_1, \tilde{\Cal F_1}) (t)
- \log T (N_2,\tilde g_2,\tilde h_2,\tilde{ \Cal F_2})
(t) \big)\cr
&=\log T_{an} (N_1,\tilde g_1,\tilde h_1, \tilde{\Cal F_1})
- \log T_{an} (N_2, \tilde g_2,\tilde h_2, \tilde{\Cal F_2})\cr
& \quad - \log T_{met} (N_1,\tilde \tau_1,\tilde g_1,\tilde{
\Cal F_1})
+ \log T_{met}
(N_2,\tilde \tau_2, \tilde g_2, \tilde{\Cal F_2})\cr}$$
and

$$\eqalign{
FT &\big( \log T_{sm} (N_1,\tilde g_1,\tilde h_1, \tilde{\Cal F_1}) (t)
- \log T_{sm} (N_2,\tilde g_2,\tilde h_2,
\tilde{\Cal F_2}) (t) \big)\cr
&=\log T_{comb} (N_1,\tilde g_1, \tilde{\Cal F_1})
- \log T_{comb} (N_2,
\tilde \tau_2, \tilde{\Cal F_2}).\cr}$$

As $N_1$ and $N_2$ are both closed manifolds
of determinant class we conclude from [BFKM] that

$$\eqalign{
\log T_{an} (N_j, \tilde g_j, \tilde{\Cal F_j})
&= \log T_{comb} (N_j,\tilde \tau_j,
\tilde{\Cal F_j}) \cr
& \, + \log T_{met} (N_j,\tilde \tau_j,\tilde g_j, \tilde{\Cal F_j}).\cr}$$

Combining the above equalities with the identity

$$\eqalign{
FT\big( \log T(N_j,\tilde g_j,\tilde h_j, \tilde{\Cal F_j}) (t)\big)
&= FT \big(\log T_{sm} (N_j,\tilde g_j,\tilde h_j,
\tilde{\Cal F_j})(t) \big) \cr
&\quad + FT\big(\log T_{la} (N_j,\tilde  g_j,\tilde h_j,
\tilde{\Cal F_j}) (t) \big) \cr}$$
one obtains, with $A$ defined as above,
$$
A(N_1, N_2, \tilde \tau_1,\tilde \tau_2, \tilde{\Cal F_1},
\tilde{\Cal F_2}) = 0.
$$
In view of the fact that A is local, one concludes that
$$
A (M_1, M_2, \tau_1, \tau_2, \Cal F_1,
\Cal F_2) = 0.$$

\hfill $\square$

Consider the pair
$\{( M, \partial_- M, \partial_+ M), \Cal F\} $
and assume that it is of determinant class.
Define $\Cal R(M,\partial_-M,g,\tau,\Cal F)$ by
$$
\log \Cal R (M,\partial_-M, g, \tau, \Cal F) :=
\log T_{an} (M,\partial_-M, g, \Cal F)
- \log T_{Re} (M,\partial_-M, \tau, g, \Cal F).$$
Our aim is to show that the ratio $\Cal R$ depends only on the data on
the boundary at least in the case when the system $( M, \partial_-
M, \partial_+ M, \Cal F) $ is
of determinant class.

First observe that the result concerning the metric anomaly of the analytic
torsion [BFKM, Lemma 6.11] extends to manifolds with boundary:

\proclaim {Proposition 2.8}  Assume that the pair $\{( M,
\partial_- M, \partial_+ M), \Cal F \}$ is of determinant class
and that $g(s), a \le
s \le b$, is a smooth 1-parameter family of Riemannian metrics on $M$ so
that $g(s)$ is independent of $s$ in a collar neighborhood of $\partial M$.
Then
$$
\frac d {ds} \log T_{an} \big( M,\partial_-M, g(s), \Cal F \big)
= \frac d {ds}
\log T_{met} \big ( M,\partial_-M, g(s), \tau, \Cal F \big).$$
\endproclaim

\demo {Proof}  Arguing as in [RS] one has

$$\eqalign{
\frac d {ds} &\log T_{an} \big(M,\partial_-M, g(s), \Cal F)\cr
&=\frac d {ds} \log T_{met} \big( M,\partial_-M, g(s), \tau,
\Cal F \big) + \sum_{q=0}^d (-1)^q \hat {q\rho_q}\cr}$$
where $\hat {\rho_q} = \underset M \to \int \rho_q (x, s)$ and
$\rho_q (x, s)$
is a density (i.e. a d-form) on $M,$  which can be calculated in
terms of the symbol of $\Delta_q (s)$ acting on $\Lambda^q (M; \Cal E).$
Due to our assumption that $g(s)$ does not depend on s near the boundary,
$\rho_q(x,s)$ vanishes for $x$ near the boundary.
 If $d = dim \, M$ is odd one can then use a parity argument to conclude that
$\rho_q(x, s) \equiv 0$ for $0 \le q \le d$.  If $d = dim \, M$ is even,
$\underset {q = 0} \to {\overset d \to \sum} (-1)^q q \hat{\rho_q}= 0.$
Indeed take the double of $M$ and use that $\log T_{an}$ and $\log T_{met}$
are 0 for closed manifolds of even dimension.
\enddemo

\proclaim {Corollary 2.9}  Suppose that the pair $\{( M, \partial_- M,
\partial_+ M), \Cal F \})$ is of determinant class and that $\tau = (h,
g') $ is a generalized triangulation. If $g_1$ and $g_2$ are two
Riemannian metrics
on $M$ which agree on $\partial M$, then

$$
\Cal R (M,\partial_-M, g_1, \tau, \Cal F)
= \Cal R (M,\partial_-M, g_2, \tau, \Cal F).$$
\endproclaim

\proclaim{Theorem 2.10} Assume that, for $j=1,2,$ $M_j = (M_j,\partial_-M_j,
\partial_+M_j)$ is a bordism  equipped with a Riemannian
metric $g_j$, a generalized triangulation $\tau_j$ and
a parallel flat bundle
$\Cal F_j = ( \Cal E_j, \nabla_j)$ of
$\Cal A$-Hilbert modules of
finite type on $M_j$.Further  assume that the systems
$(M_j,\partial_-M_j, \partial_+M_j,\Cal F_j)$ are of determinant class and
$$
(\partial_\pm M_1, g_1 \restriction_{\partial_\pm M_1}, \Cal F_1
  \restriction_{\partial_\pm M_1})
=(\partial_\pm M_2, g_2 \restriction_{\partial_\pm M_2}, \Cal F_2
  \restriction_{\partial_\pm M_2}). $$
Then
$$ \Cal R ( M_1, \partial_-M_1, g_1, \tau_1, \Cal F_1) =
\Cal R ( M_2,\partial_-M_2, g_2, \tau_2, \Cal F_2).$$

\endproclaim

\demo{Proof}  The proof proceeds in several steps.

(I)  In the case where the hypotheses of Corollary C are satisfied the
claim follows from Theorem A.

(II)  Next assume that only the hypotheses of Theorem B are
satisfied. Consider $S^2$ equipped with a generalized triangulation
$\tau_{S^2}=(h_0,g_0)$ and denote by
$\underline \tau_j =(\underline h_j,\underline g'_j)$ a product
triangulation on $M_j\times S^2$. We claim that, for $j=1,2,$
$(M_j \times S^2, g_j \times g_0, \underline \tau_j,
\underline {\Cal F_j}),$ where $\underline {\Cal F_j} $ is the pullback of
$\Cal F_j$
by the projection $M_j\times S^2 \to M_j,$
satisfies the assumptions of Corollary C. Indeed
denote by $D^3$ the unit disc in $\Bbb R^3$ and choose for
$M_3$ the disjoint union of $\partial_-M_1\times D^3$ with
$\partial_+M_1\times D^3.$ Then $\partial_{\pm}M_3=
\partial_{\pm}M_1\times S^2.$
Let $\Cal F_3$
be the pull back
of $\Cal F_1\restriction_{\partial M_1}$ by the
projection $\partial M_1\times D^3 \to \partial M_1$. Since the pairs
$\{\partial_{\pm}M_i,\Cal F_i\restriction_{\partial_{\pm}M_{i}}\}$
are of determinant class,
the system $(M_3,\partial_-M_3, \partial_+M_3,\Cal F_3)$ is of determinant
class and therefore, by Proposition 2.3, the pairs
$\{N_i, \tilde{\Cal F_i}\}$ (as defined in Corollary C) are
of determinant class.
One concludes then by (I) that
$$
\Cal R (M_1 \times S^2, \partial_-M_1\times S^2, g_1 \times g_0,
\underline \tau_1,
\underline {\Cal F_1}) = \Cal R (M_2 \times S^2,
\partial_-M_2 \times S^2 , g_2 \times g_0,
\underline \tau_2, \underline {\Cal F_2}).$$
By applying Proposition 2.5 (3) and Proposition 2.3 (formulas (2.10)-
(2.10')), one obtains

$$
\Cal R (M_1, \partial_-M_1, g_1, \tau_1, \Cal F_1)=\Cal R (M_2,
\partial_-M_2, g_2, \tau_2,
\Cal F_2).$$

(III) Applying Corollary 2.9 it suffices to prove the statement in the case
where $g_j$ and $\tau_j = (h_j, g_j')$ $(j= 1,2)$ have the additional
property that $g_j = g_j'$.

(IV)
(A) It suffices to prove the result under the additional
hypothesis that $\chi(M_1,\partial_-M_1) = \chi(M_2,\partial_-M_2).$
Indeed, if $dim M_1$ is odd,
\newline $\chi (M_1,\partial_-M_1)=\chi
(M_2,\partial_-M_2).$ If $dim M_1$ is even, then
$$\chi(M_1,\partial_-M_1)= \chi(M_2,\partial_-M_2) + 2k.$$
If, in addition, $k > 0$ we replace, without loss of generality,
$M_1$ by the disjoint union of $M_2$ with k copies of $(S^d,\Cal F),$
where $\Cal F$ is the trivial parallel flat bundle with the same
fiber as $\Cal F_1.$ (Observe that $\Cal R =0$ for $(S^d,\Cal F).$ )
In the case when $k<0$ interchange the role of $M_1$ and $M_2.$

(B) Under the additional hypotheses
$\chi(M_1,\partial_-M_1)=\chi(M_2,\partial_-M_2)$ the result can be
proved as follows:
By the
invariance of the Reidemeister torsion under a subdivision of a
generalized triangulation (cf. Lemma 6.12 in [BFKM] for the case of a
closed manifold) it suffices to prove that there exist
subdivisions $\hat \tau_j$ of
$\tau_j$ so that the assumptions of Theorem B are satisfied.
But such
subdivisions exist
if $\chi (M_1,\partial_-M_1) = \chi(M_2,\partial_- M_2)$ (cf. [BFKM]).
\hfill $\square$

\newpage

\input amstex.tex
\documentstyle{amsppt}
\magnification=1200
\baselineskip 16pt plus 2pt
\parskip 2pt
\NoBlackBoxes

\def\Vol{\text {Vol}}

\document

\proclaim{3.  Applications to torsions}  \endproclaim
\proclaim{3.1 Comparison of analytic and Reidemeister torsion}
\endproclaim

In this subsection we provide a formula for the ratio of analytic and
Reidemeister torsion of a compact Riemannian
manifold with boundary as introduced at the end of
section 2.

We begin by some auxiliary considerations.
Assume that $(M,g)$ is a closed Riemannian manifold, $\tau =(h,g')$
a generalized triangulation,  and $\Cal F =
( \Cal E, \nabla)$ a parallel flat
bundle on $M.$
Consider the cylinder
$M_I = M \times I$ where $I$ is a compact interval $[a, b]$ in $\Bbb R$ with
$a,b \in 4\Bbb Z.$
Denote by $g_0$ the
standard metric on $[a, b]$ and by $\tau_i = (h_i, g_0),\,(i=1,2),$
the generalized triangulations of the bordisms
$([a,b],\{a\},\{b\})$ resp.
$([a,b], \emptyset,  \{a,b\})$
with $h_1(x)= x$ resp.
$h_2= \frac {1}{2}(x- (b+a)/2)^2.$ ( To satisfy condition (T1) in the
definition of a generalized triangulation, one can perturb $h_2$ slightly
so that it is linear near the boundary of the interval.)
Notice that the generalized triangulation $\tau_1$ has no cells.
Denote by $\Cal T$ the trivial 1-dimensional complex
vector bundle with
trivial connection. $\Cal T$ is a parallel flat bundle over $[a,b].$
An easy calculation gives

$$\leqalignno{
\log T_{comb} &(I,\{a\},\tau_1, \Cal T) = 0,\
\log T_{an} (I, \{a\}, g_0, \Cal T) = \frac 1 2 \log 2,&(3.1)\cr
&\log T_{met} (I ,\{a\},\tau_1, g_0,\Cal T) = 0,\cr}$$

$$\leqalignno{
\log T_{comb} &(I,\emptyset,\tau_2, \Cal T) = 0,\
\log T_{an} (I,\emptyset, g_0, \Cal T) =
\frac 1 2 (\log 2 +\log (b-a)),&(3.1')\cr
&\log T_{met} (I,\emptyset,\tau_2, g_0,\Cal T) =
\frac 1 2 \log (b-a).\cr}$$

Assume that the pair $\{M,\Cal F\}$ is of
determinant class. Then, by Proposition 2.3, the systems
$(M_I,\emptyset,M\times \partial I, \Cal F_I)$
and $(M_I,M\times \{a\},M\times \{b\}, \Cal F_I)$ are of determinant class
(where $\Cal F_I$ is the pull back of $\Cal F$ by the projection
$M\times I \to M)$ and

$$\leqalignno{
\log T_{comb}&( M_I, M\times \{a\}, \tau\oplus \tau_1, \Cal F_I)=0,
&(3.2)\cr}
$$
$$
\log T_{an} ( M_I, M\times \{a\}, g\oplus g_0, \Cal F_I) =
\chi (M;\Cal F)\frac {\log2}{2} $$
$$
\log T_{met} ( M_I, M\times \{a\}, g\oplus g_0,
\tau\oplus \tau_2,\Cal F_I) =0,$$

$$\leqalignno{
\log T_{comb}&( M_I, \emptyset, \tau\oplus \tau_2,
\Cal F_I)=\log T_{comb}(M, \tau, \Cal F),& (3.2')\cr
\log T_{an} ( M_I, \emptyset&, g\oplus g_0, \Cal F_I) =
\frac {\chi (M;\Cal F)}{2} (\log 2 +
\log(b-a))+\log T_{an}(M,g,\Cal F),\cr
\log T_{met} ( M_I,& \emptyset, g\oplus g_0, \tau\oplus \tau_2,\Cal F_I)
=\frac {\chi (M;\Cal F)}{2} \log(b-a) +\log T_{met}(M,g,\tau, \Cal F)
\cr}$$
(where we have used that $\chi(I,\emptyset)=1$ and $\chi(I,\{a\})=0$).

Since M is closed and the pair $\{M,\Cal F\},$ is of determinant class
one obtains, by [BFKM],
$$\log T_{an}(M,g,\Cal F) = \log T_{Re}(M,g,\tau,\Cal F).$$
Combining
with the formulas  (3.2)-(3.2'), one obtains
$$\leqalignno{
\log \Cal R (M_I, M\times \{a\}, g\otimes g_0,\tau\oplus
\tau_1, \Cal F_I) = \frac {\chi ( M;\Cal F)
\log 2} 2 .&&(3.3)\cr}$$

$$\leqalignno{
\log \Cal R (M_I,\emptyset, g\otimes g_0,\tau\oplus \tau_2, \Cal F_I)
= \frac {\chi ( M;\Cal F)
\log 2} 2 .&&(3.3')\cr}$$

\medskip

\proclaim {Theorem 3.1}  Assume
that the system $(M, \partial_-M, \partial_+M,\Cal F)$ is of determinant
class.  Then
$$\leqalignno{
\log \Cal R (M, \partial_-M, g,\tau, \Cal F) = \frac 1 4 \chi (\partial
M;\Cal
F) \log 2.
&&(3.4)\cr}$$
\endproclaim

\remark{Remark}  In the case where $ \Cal A = \Bbb R$, the result, as it
stands,
is due independently to L\"uck ([L\"u1]) and Vishik [Vi]
(cf. also [Ch]). However, even in the case $\Cal A = \Bbb R$, the proof
presented here
is new, elementary and short. \endremark

\demo {Proof}  The proof proceeds in two steps.

(A)  We first consider the case where $M = (M, \emptyset, \partial M)$.  Take
a \underbar {disjoint} union $M \sqcup M$ of two copies
of $M$ and consider the
bordism $\big(M \sqcup M, \emptyset, (\partial M) \sqcup
(\partial M)\big)$.  Further introduce $(\partial M)_I :=
\partial M \times I$ and consider the bordism
$\big( ( \partial M)_I, \emptyset, \partial (
(\partial M)_I) \big)$.  Notice that $\partial((\partial M)_I)
= ( \partial M
\times \{ 0 \} ) \sqcup (\partial M \times \{ 1 \} ) = (\partial M) \sqcup
(\partial M) .$ Since by hypothesis $(\partial M, \Cal F
\restriction_{\partial M})$
is of determinant class so is  $( \partial M \times I,\emptyset, (\Cal F
\restriction_{\partial M})_I )$ (cf Proposition 2.5(3)). Therefore, by
Theorem 2.10,

$$
\log \Cal R( \partial M \times I, \emptyset, g_I,
\tau_I, (\Cal F \restriction_{\partial M})_I ) =
\log \Cal R ( M \sqcup M, \emptyset, g\sqcup g,
\tau \sqcup \tau, \Cal F \sqcup \Cal F)$$

where $(\Cal F\restriction_{\partial M})_I$ is
defined as above, $g_I:= g\restriction_{\partial M} \oplus g_{0}$
and $\tau_I$ is a product triangulation.
Notice that

$$\leqalignno{
\log \Cal R( M \sqcup M, \emptyset, g\sqcup g,
\tau \sqcup \tau, \Cal F \sqcup \Cal F)=
2 \log \Cal R (M, \emptyset, g, \tau, \Cal F),
&&\cr}$$
and, by (3.3'),

$$\leqalignno{
\log \Cal R ( \partial M \times I, \emptyset, g_I,
\tau_I, \Cal F_I \restriction_{\partial M} ) =
\frac {\chi (\partial M;\Cal F)} 2  \log 2.&&\cr}$$
Combining the three equalities above, we obtain

$$
\log \Cal R (M, \emptyset, g,\tau, \Cal F ) =
\frac {\chi (\partial M;\Cal F)} 4 \log 2$$
and statement (3.4) is proved in the case where $\partial_-M = \emptyset.$

(B)  In view of Proposition 2.3 it suffices to check (3.4) for the
system
$(M\times S^2, \partial_-M\times S^2, \partial_+M\times S^2,\tilde \Cal F)$
which is, by Proposition 2.5 (3), of determinant class because
$(M,\partial_-M, \partial_+M, \Cal F)$ is. Here
$\tilde{\Cal F}$ denotes the pullback of $\Cal F$ by the
projection $M\times S^2 \to M.$
In view of Theorem 2.10
it suffices to prove the statement (3.4) for the two systems
$(\partial_{+}M\times D^3,\emptyset, \partial_{+}M\times S^2,
(\Cal F \restriction_{\partial_{+} M})\tilde{ })$ and
$(\partial_- M \times D^3, \partial_-M\times S^2, \emptyset,
(\Cal F \restriction _{\partial_-M}) \tilde{ })$ where
$(\Cal F \restriction_{\partial_{\pm} M} )\tilde{ }$ denotes the pullback
of $\Cal F\restriction_{\partial_{\pm} M}$ by the projections
$\partial_{\pm}M \times D^3 \to \partial_{\pm}M.$ For the first system it
suffices to apply (A) and for the second system, (3.4) follows from
Proposition 2.2 (Poincar\'e duality) and (A).
\hfill $\square$

\proclaim{3.2. Glueing formulas (Part I)}\endproclaim

In this subsection we present a gluing formula for the
analytic and Reidemeister torsions.  Let
$(M_j, \partial_- M_j, \partial_+M_j)$, be two bordisms equipped with
Riemannian metrics $g_j,$ generalized triangulations $ \tau_j =
(h_j, g_j')$ and parallel flat bundles
$\Cal F_j = (\Cal E_j, \nabla _j).$
Suppose that there exists an isometry

$$\leqalignno{
&\omega : (\partial_+ M_1, g_1 \restriction_{\partial_- M_1}) \rightarrow
(\partial_- M_2, g_2 \restriction_{\partial_+ M_2}) &(3.5)\cr}$$
and a connection preserving bundle isomorphism $\Phi$ above $\omega$ which
makes
the following diagram commutative

$$\leqalignno{
&\matrix
\Cal E_1 \restriction_{\partial_- M_1} & \overset \Phi \to \longrightarrow
&\Cal E_2 \restriction_{\partial_+ M_2}\\
\downarrow && \downarrow\\
\partial_- M_1 & \overset \omega \to \longrightarrow & \partial_+ M_2.
\endmatrix &(3.6)\cr}$$
Then one can form a bordism $(M := M_1 \cup_\omega M_2, \partial_- M_1,
\partial_+ M_2)$ by gluing $\partial_+ M_1$ to $\partial_- M_2$ by $\omega$
and a parallel flat bundle $\Cal F$ by gluing $\Cal F_1$ and $\Cal F_2$ by
$(\omega, \Phi)$.  The metrics $g_1, g_1'$ and $g_2, g_2'$ determine
Riemannian
metrics $g$, $g'$ on $M$, and the functions $h_1$ and $h_2$ determine the
$C^{\infty}$-function $h : M \rightarrow {\Bbb R}$ given by
$$\leqalignno{
&h(x):=\cases h_1(x)  &(x \in M_1)\\
b_1 - a_2 + h_2(x) & (x \in M_2)\endcases &(3.7)\cr}$$
where, for j=1,2, $h_j(M_j)=[a_j,b_j].$
\newline
Notice that if $\tau = (h,g')$ is not a generalized triangulation
(i.e. violates (T4)) one can modify, by
an arbitrary small perturbation and
localized in a given neighborhood of $\partial_+ M_1,$
the metric $g_1'$ to $\overset \sim
\to g_1'$ so that
the triangulation $\overset \sim \to \tau_1:=(h_1, \overset \sim\to g_1')$
is compatible with $\tau_2$ and $\overset \sim \to \tau_1 $ and
$\tau_1$ provide the same relative CW-complex structure.
Therefore, without loss of generality, we may assume that
$\tau = (h, g')$ is a generalized triangulation
of $(M, \partial_-M_1 ,\partial_+ M_2).$

Choose orientations $\Cal O_h$
for $M.$ They induce orientations $\Cal O_{h_i}$ for $M_i.$
Now consider the short exact sequence of cochain complexes
$$\leqalignno{
0 \rightarrow \Cal C^* (M_2,\partial_-M_2, \tau_2,\Cal O_{h_2};\Cal F_2)
&\overset {i_*} \to \rightarrow
\Cal C^*(M,\partial_-M_1, \tau, \Cal O_h ;\Cal F) \cr
&\overset {r_*} \to \rightarrow
\Cal C^*(M_1,\partial_-M_1, \tau_1, \Cal O_{h_1};\Cal F_1)
\rightarrow 0,&(3.8)\cr}$$
where $i_* $ is the map which extends the cochains, defined on the cells of
$M_2,$ to cochains defined on all cells of $M$ by assigning the value
zero on the cells of $M_1$ and where $r_*$
is the map which restricts the cochains, defined on the
cells of $M,$ to the cells of $M_1.$  The sequence (3.8) induces a
long weakly exact sequence $\Cal H_{comb}(\tau)$ in cohomology
of $\Cal A$- Hilbert modules (cf (1.52))

$$\leqalignno{
\cdots \rightarrow  H^q (M_2,\partial_-M_2, \tau_2 ;\Cal F_2)
&\rightarrow H^q(M, \partial_- M_1, \tau ;\Cal F)
\rightarrow  H^q (M_1,\partial_- M_1, \tau_1 ; \Cal F_1)
\rightarrow \cr & H^{q+1} (M_2,\partial_-M_2, \tau_2 ;\Cal F_2)\rightarrow
\cdots
&(3.9)\cr}$$
A similar sequence (depending on the Riemannian metric $g$),
 denoted by $\Cal H_{an}(g),$ can be obtained
using de Rham cohomology $H^q(M,\partial_-M, g ;\Cal F)$
instead of the combinatorial one,
$$\leqalignno{
\cdots \rightarrow  H^q(M_2,\partial_-M_2, g_2;\Cal F_2)
&\rightarrow H^q(M, \partial_- M, g;\Cal F)
\rightarrow  H^q (M_1,\partial_- M_1, g_1; \Cal F_1)
\rightarrow \cr & H^{q+1}(M_2,\partial_-M_2, g_2 ;\Cal F_2)\rightarrow \cdots
&(3.10)\cr}$$
As we have seen in section 2, the integration theory
provides isomorphisms $(\theta_j)^{-1}$ of $\Cal A$-Hilbert modules
from de Rham cohomology to the combinatorial cohomology. Moreover the
$(\theta_j)$ define an isomorphism of cochain complexes from
$\Cal H_{comb}(\tau)$ to $\Cal H_{an}(g).$

\proclaim{Theorem 3.2}
Assume that, for i=1,2, the system
$(M_i,\partial_-M_i,\partial_+M_i,\Cal F_i)$
is of determinant
class.  Then the following statements hold:
$$\leqalignno{
&\text{The system} (M,\partial_-M,\partial_+M, \Cal F)\ \text{and the
complexes}\ \Cal H_{comb}, \Cal H_{an}\cr &\text{ are of
determinant class;}&(i)\cr}$$
$$\leqalignno{
\log T_{Re} (M, \partial_-M, g, \tau, \Cal F) =
&\overset 2 \to
{\underset {j = 1} \to \sum} \log T_{Re} (M_j, \partial_-M_j, g_j,
\tau_j,\Cal
F_j)
&{(ii)}\cr
&+ \log T (\Cal H_{an}).\cr}$$
$$\leqalignno{
\log T_{an} (M, \partial_-M, g, \Cal F)
&= \sum^2_{j = 1} \log T_{an} (M_j, \partial_-M_j, g_j, \Cal F_j)
+ \log T(\Cal H_{an}) &{(iii)}\cr
&-\frac {\chi ( \partial_+ M_1;\Cal F_1\restriction_{\partial_+M_{1}})} 2
\log 2.\cr}$$
\endproclaim

\remark{Remark}  In the case where $\Cal A =\Bbb R$, the result (ii)
is due
to Milnor
[Mi].  Work related to (i) can be found in [LL].
In the case where $\Cal A = \Bbb R$, (iii) is due to
Vishik [Vi] .  Vishik's proof,
however, is completely different than ours.
\endremark

\demo{Proof}  (i) Proposition 2.5 implies that the system
$(M,\partial_-M,\partial_+M,\Cal F)$ is of determinant class. Theorem
1.14 applied to the short exact sequence of cochain complexes (3.8)
implies that $\Cal H_{comb}$ is of determinant class.
In view of the fact that $\theta$ defines an isomorphism of
cochain complexes, Proposition 1.3 implies that $\Cal H_{an}$ is of
determinant class as well.

(ii)  Concerning the torsion of the sequence (3.8), notice that
$$\leqalignno{
\log T (0 \rightarrow
\Cal C^q (M_2,\partial_-M_2,&\tau_2, \Cal O_{h_{2}}; \Cal F)
\overset {i_q} \to \rightarrow
\Cal C^q (M,\partial_- M, \tau, \Cal O_h ;\Cal F)  \cr
&\overset {r_q} \to \rightarrow
\Cal C^q (M_1,\partial_-M_1,\tau_1, \Cal O_{h_{1}}; \Cal F)
\rightarrow 0)=0 \cr}$$
since $\log Vol(i_q)=0$ and
$\log Vol(r_q\restriction (Null(r_q)^{\bot})=0.$ Theorem 1.14 implies that
$$\leqalignno{
\log T_{comb} (M,\partial_-M, & \tau, \Cal F) =
&{(3.11)} \cr
&\overset 2 \to
{\underset {j = 1} \to \sum}
\log T_{comb} (M_j,\partial_-M_{j}, \tau_j,\Cal F_j)
+ \log T (\Cal H_{comb}). \cr}$$

{}From Proposition 1.3 (C) and in view of the fact that
$\theta$ is an isomorphism of cochain complexes of $\Cal A-$Hilbert
modules one obtains
$$\leqalignno{
\log T_{met} (M,\partial_- M, g ,\tau, \Cal F) =
&\overset 2 \to{\underset {j = 1} \to \sum}
\log T_{met} (M_j,\partial_-M_{j}, g_j,\tau_j,
\Cal F_j)&{(3.12)}\cr
&+ \log T (\Cal H_{an}) - \log T( \Cal H_{comb}). \cr}$$
Combining (3.11) and (3.12) one obtains (ii).

(iii) follows from (ii) and Theorem 3.1.
\hfill $\square$
\enddemo

\proclaim {3.3. Glueing formulas (Part II)}\endproclaim

In this subsection we extend Theorem 3.2 for partial glueing, i.e.
to a situation where
not necessarily all of the components of
$\partial_+M_1$ and $\partial_-M_2$ are glued together.

Suppose that $(M_i,\partial_-M_i,\partial_+M_i ;\Cal F_i), i=1,2$ are two
systems with $\partial_+M_1$ resp. $\partial_-M_2$ consisting of two
disjoint components $\partial_+M_1=V_0\sqcup V_1,$ resp.
$\partial_-M_2=W_0\sqcup W_1.$
Assume that there exist a diffeomorphism
$\omega_0 : V_0\to W_0$
and a connection preserving  isomorphism above $\omega_0,$
$\Phi_0 : \Cal F_1\restriction_{V_0} \to \Cal F_2 \restriction_{W_0}.$
Then $M_1$ and $M_2$ can be glued together to give rise to the
bordism $(M, \partial_-M,\partial_+M)$ with $M:= M_1\cup_{\omega}
M_2$ and
boundary components $\partial_-M=\partial_-M_1\sqcup W_1$
and $\partial_+M= V_1 \sqcup \partial_+M_2.$ Using $\Phi_0,$
 $\Cal F_1$ can
be glued with $\Cal F_2$ to obtain a parallel flat bundle $\Cal F$
on $M.$
Suppose, in addition, that Riemannian metrics $g_1$ on $M_1$ and
$g_2$ on
$M_2$ are given so that
$\omega^{\#}(g_2\restriction_{W_{0}})=g_1\restriction_{V_{0}}.$
The metrics $g_1$ and $g_2$ define a
Riemannian metric $g$ on $M$. One then obtains the following (weakly exact)
sequence $\Cal H$ in de Rham cohomology
$$\leqalignno{
\cdots \rightarrow & H^q (M_2,\partial_-M_2, g_2;\Cal F_2)
\rightarrow H^q(M, \partial_- M,g;\Cal F)
\rightarrow\cr &\rightarrow H^q (M_1,\partial_- M_1, g_1; \Cal F_1)
\rightarrow  H^{q+1} (M_2,\partial_-M_2, g_2;\Cal F_2)
\rightarrow \cdots &(3.19)\cr}$$
This sequence is the cohomology sequence induced by the exact sequence of
cochain complexes (cf (2.7))
$$\leqalignno{
0\rightarrow \Lambda^{\ast}(M_2,\partial_-M_2; \Cal F_2)\rightarrow
\Lambda^{\ast}(M,\partial_-M; \Cal F)\rightarrow
\Lambda^{\ast}(M_1,\partial_-M_1; \Cal F_1)\rightarrow 0. &&(3.20)\cr}$$
(To be completely correct, $\Lambda^q(M_2,\partial_-M_2;\Cal F_2)$ should
be the smaller space of smooth $q-$forms which are the restriction
to $M_2$ of $q-$forms which are defined on all of M, but vanish on $M_1.)$

\proclaim{Theorem 3.2$'$} If the systems
$(M_i,\partial_-M_i,\partial_+M_i, \Cal F_i)$ (i=1,2) are of determinant
class then so are the system $(M, \partial_-M, \partial_+M, \Cal F)$
and the cohomology sequence $\Cal H,$ given by (3.19). Moreover
$$\leqalignno{
\log T_{an} (M,\partial_-M, g, \Cal F)
&= \sum^2_{j = 1} \log T_{an} (M_j,\partial_-M_j, g_j, \Cal F_j)
+ \log T(\Cal H) \cr
&-\chi (V_0;\Cal F_1\restriction_{V_{0}} )\frac {\log 2} {2}.&\cr}$$
\endproclaim
\demo{Proof} We will derive the stated results from Theorem 3.2. For that
purpose introduce, for
$\epsilon >0,$ the systems $(N_{i,\epsilon},
\partial_-N_{i,\epsilon},\partial_+N_{i,\epsilon},
\Cal G_{i,\epsilon}), i=1,2,$ defined by
$$
N_{1,\epsilon}:=W_1\times [-\epsilon,0],
\ \partial_-N_{1,\epsilon}:=W_1\times \{-\epsilon\},\
\partial_+N_{1,\epsilon}:=W_1\times \{0 \}
;$$
$$
N_{2,\epsilon}:=V_1\times [0,\epsilon],
\ \partial_-N_{2,\epsilon}:=V_1\times \{0 \},
\ \partial_+N_{2,\epsilon}:=V_2\times \{\epsilon \};
$$
$$
\Cal G_{1,\epsilon}:=(\Cal F_2\restriction _{W_1})_{[-\epsilon ,0]},\
\Cal G_{2, \epsilon}:= (\Cal F_1\restriction _{V_1})_{[0,\epsilon]},$$
and let $g_{i,\epsilon}$ be the metrics on $N_{1,\epsilon}$
defined by
$g_{1,\epsilon}:= g_2\restriction _{W_1}\oplus g_0$ and
$g_{2,\epsilon}:=g_1\restriction _{V_1}\oplus g_0$
(where $g_0$ is the standard Euclidean  metric) (see Figure 2).
\newline Denote by $(M_{i,\epsilon},
\partial_-M_{i,\epsilon},\partial_+M_{i,\epsilon}, \Cal F_{i,\epsilon})$
the disjoint union of the systems
\newline $(M_i,
\partial_-M_i,\partial_+M_i, \Cal F_i)$
and $(N_{i,\epsilon},
\partial_-N_{i,\epsilon},\partial_+N_{i,\epsilon}, \Cal G_{i,\epsilon})$
and let $\omega_\epsilon: \partial_+M_{1,\epsilon} :=V_0\sqcup V_1\sqcup W_1
\to \partial_-M_{2,\epsilon} := W_0\sqcup V_1\sqcup W_1$
be the diffeomorphism which is equal to
$\omega_0$ on $V_0$ and to the identity on $V_1\sqcup W_1.$ Moreover
$\omega_{0} $ and $\Phi_{0}$ induce an isomorphism of parallel flat
bundles
$$\Phi_\epsilon:\Cal F_{1,\epsilon}\restriction_{\partial_+M_{1,\epsilon}}\to
\Cal F_{2,\epsilon}\restriction_{\partial_-M_{2,\epsilon}}.$$
Following the glueing instructions described in subsection 3.2, one
forms the system
$(M_\epsilon,\partial_-M_\epsilon,\partial_+M_\epsilon,
\Cal F_\epsilon)$
and defines the Riemannian metric $g_\epsilon$ given by
$g_i$ on $M_i$ and $g_{i,\epsilon}$ on $N_{i,\epsilon}.$
\newline There exists a
commutative diagram connecting the cohomology sequence
$\Cal H$, defined in (3.19),
with the sequence (3.10) for the bordism $(M_\epsilon,\partial_-M_\epsilon,
\partial_+M_\epsilon),$ which we denote by $\Cal H_\epsilon,$
{\eightpoint
$$
\matrix &\rightarrow &H^q(M_2,\partial_-M_2,g_2;\Cal F_2)&
\rightarrow &
H^q(M,\partial_-M,g;\Cal F)
&\rightarrow &H^q(M_1,\partial_-M_1,g_1;\Cal F_1) &\cdots\\
&&\qquad\downarrow (i_{2,\epsilon})_q&&\quad\qquad\qquad\downarrow
(i_{\epsilon})_q &&\qquad\downarrow (i_{1,\epsilon})_q\\
& \rightarrow &H^q(M_{2,\epsilon},\partial_-M_{2,\epsilon},
g_{2,\epsilon};\Cal F_{2,\epsilon}) & \rightarrow&
H^q(M_{\epsilon},\partial_-M_{\epsilon},g_{\epsilon};\Cal F_{\epsilon})
&\rightarrow &
H^q(M_{1,\epsilon},\partial_-M_{1,\epsilon},g_{1,\epsilon};
\Cal F_{1,\epsilon}) & \cdots.\endmatrix$$ }

In this diagram, the horizontal lines represent the cohomology
sequence (3.19)
and (3.10) respectively and the vertical arrows are
isomorphisms of Hilbert
modules of finite type. The isomorphisms $(i_{1,\epsilon})_q,\
(i_{\epsilon})_q \ \text { and } (i_{2,\epsilon})_q $ can be described
as the composition
$$
(i_{1,\epsilon})_q:= (i'_{1,\epsilon})_q \cdot
((i''_{1,\epsilon})_q)^{-1},\
(i_{2,\epsilon})_q := (i'_{2,\epsilon})_q\cdot
((i''_{2,\epsilon})_q)^{-1},$$

$$(i_{\epsilon})_q= (i'_{\epsilon})_q\cdot
((i''_{\epsilon})_q)^{-1}
$$
where $(i'_{k,\epsilon})_q$ and $(i''_{k,\epsilon})_q$
are induced by the inclusions (k=1,2)
$$
\leqalignno{
i'_{k,\epsilon}&:(M_k,\partial_-M_k)\to (M_{k,\epsilon},
\partial_-M_{k,\epsilon}\cup N_{k,\epsilon})\ \text{and}\cr
&i''_{k,\epsilon}:(M_{k,\epsilon},\partial_-M_{k,\epsilon})\to
(M_{k,\epsilon},
\partial_-M_{k,\epsilon}\cup N_{k,\epsilon})&\cr}
$$
and
where $(i'_{\epsilon})_q$ and $(i''_{\epsilon})_q$
are induced by the inclusions
$$
\leqalignno{
i'_{\epsilon}:&(M_,\partial_-M)\to (M_{\epsilon},
\partial_-M_{\epsilon}\cup N_{1,\epsilon}) \text{ and }\cr&
i''_{\epsilon}:(M_{\epsilon},\partial_-M_{\epsilon})\to (M_{\epsilon},
\partial_-M_{\epsilon}\cup N_{1,\epsilon}).&\cr}
$$
Since the harmonic forms, which represent the cohomology classes in
\newline $H^{\ast}(M_k,\partial_-M_k,g_k ;\Cal F_k)$ and
$H^{\ast}(M_{k,\epsilon},\partial_-M_{k,\epsilon},g_{k,\epsilon}
;\Cal F_{k,\epsilon})$ are in fact the same, the homomorphisms
$(i_{k,\epsilon})_q$ are isometries for any $\epsilon$ (k=1,2).
Notice that one can provide
a smooth family of diffeomorphisms $\varphi_\epsilon:M \to M_{\epsilon},$
with $\varphi_0=id$ so that
$\varphi_\epsilon^{\#}(g_\epsilon)=\tilde g_\epsilon$ is a smooth
family of Riemannian metrics on $M$ with $\tilde g_0=g.$ This implies that
$$\lim_{\epsilon\to 0}\log T_{an}(M_\epsilon,\partial_-M_\epsilon,g_\epsilon,
\Cal F_\epsilon)=
\log T_{an}(M,\partial_-M,g,\Cal F).$$
Using $\varphi_\epsilon,$ one can also show that
$$\lim _{\epsilon\to 0}\log Vol ((i_{\epsilon})_q)=0.$$
In view of Proposition 1.3 A(iii) $\Cal H$ is of
determinant class iff $\Cal H_{\epsilon}$ is of determinant class
(for one and then for any $\epsilon > 0 $)
and if so
$$\log T(\Cal H)= \log T(\Cal H_\epsilon) +
\sum_q (-1)^q\log Vol (i_\epsilon)_q.$$
Note that if the systems $(M_i,\partial_-M_i,\partial_+M_i,
\Cal F_i), i=1,2,$ are of determinant class, then the systems
$(M_{i,\epsilon},\partial_-M_{i,\epsilon},\partial_+M_{i,\epsilon},
\Cal F_{i,\epsilon}), i=1,2,$ are of determinant class as well.
{}From (3.2) one obtains
$$
\log T_{an}(M_{i,\epsilon},\partial_-M_{i,\epsilon},
g_{i,\epsilon},\Cal F_{i,\epsilon})=
\log T_{an}(M_i,\partial_-M_i,g_i,\Cal F_i) +\chi_i\cdot \frac{\log 2}{2}$$
with $\chi_1=\chi(W_1;\Cal F_2\restriction_{W_1})$ and
$\chi_2=\chi(V_1;\Cal F_1\restriction_{V_1}).$ Then by Theorem 3.2
$$
\leqalignno{
\log T_{an}(M_{\epsilon},\partial_-M_{\epsilon},
g_{\epsilon},&\Cal F_\epsilon)=
\sum_{i=1,2}\log T_{an}(M_i,\partial_-M_i, g_i,\Cal F_i) +\cr &+\log T(\Cal
H_\epsilon)
- \chi(V_0,\Cal F_1\restriction_{V_0})\cdot \frac{\log 2}{2}&\cr}$$
The result follows by passing to the limit $\epsilon \to 0.$
\hfill $\square$
\enddemo

\proclaim {3.4. Comparison of the analytic torsion for different
boundary conditions}\endproclaim

Suppose $(M,g)$ is a Riemannian manifold whose boundary $\partial M$
is a disjoint union of three components $\partial_1M$, $\partial_2M$,
and $\partial_3M$ and $\Cal F$ is a parallel flat bundle
of $\Cal A-$Hilbert modules on $M$.
In this subsection we will compare the analytic
torsions of  $(M,\partial_1M , g, \Cal F)$ and
$(M,\partial_1M \cup \partial_3M, g, \Cal F).$

Introduce $g_3:= g\restriction_{\partial_3M}$ and
$\Cal F_3:= \Cal F\restriction_{\partial_3M}$ and denote by
 $\Cal H'$ the cohomology sequence

$$\leqalignno{
&\cdots
\rightarrow H^q(M, \partial_1 M\cup\partial_3M,g;\Cal F)
\rightarrow  H^q (M,\partial_1M ,g; \Cal F)
\rightarrow \cr & H^{q} (\partial_3 M, g_3; \Cal F_3)
\rightarrow  H^{q+1} (M,\partial_1 M\cup \partial_3 M,g; \Cal F)
\rightarrow \cdots &(3.21)\cr}$$
which is induced by the short exact sequence
$$
0\rightarrow \Lambda^{\ast}(M,\partial_1M\cup \partial_3M; \Cal F)\rightarrow
\Lambda^{\ast}(M,\partial_1M; \Cal F)\rightarrow
\Lambda^{\ast}(\partial_3M; \Cal F_3)\rightarrow 0.$$

\proclaim{Theorem 3.3} The system
$(M,\partial_1M, \partial_2M \cup \partial_3M, \Cal F)$
is of determinant
class iff the system $(M, \partial_1M \cup \partial_3M, \partial_2M, \Cal F)$
is of determinant class and if so
$$
\leqalignno{
\log T_{an} (M,&\partial_1M, g, \Cal F)
= \log T_{an} (M,\partial_1M \cup \partial_3M, g, \Cal F)
+ \log T(\Cal H') +
\cr \log &T_{an}(\partial_3M, g_3, \Cal F_3)
.&(3.22)\cr}$$
\endproclaim

\proclaim{Remark} For $\Cal A = \Bbb C,$ the result, as it stands, is
implicit in Vishik [Vi3]. Again his proof is very different from ours.
\endproclaim

\demo{Proof} The first statement follows from Proposition 2.5 (4).
To prove (3.22), we apply Theorem 3.2$'$ (partial gluing theorem)
to the systems
$$(M_1,\partial_-M_1, \partial_+M_1, \Cal F_1):=
(M,\partial_1M, \partial_2M \cup \partial_3M, \Cal F)$$
and, for $\epsilon >0,$
$$(M_2,\partial_-M_2, \partial_+M_2, \Cal F_2) :=
(\partial_3M\times [0,\epsilon],\partial_3M\times\{0,\epsilon\},
\emptyset ,\Cal F_{[0,\epsilon]})$$
where $\Cal F_{[0,\epsilon]}$ is the pull back of $\Cal F
\restriction_{\partial_{3}M}$ by the projection
$\partial_3M\times [0,\epsilon] \to \partial_3M.$
We regard $M_1$ equipped with the metric g and $M_2$ with
the metric
$g_3\oplus g_0$ ($g_0$ the Euclidean metric).
Further, using the same notation as in subsection 3.3,

$$V_0= W_0:= \partial_3M,\ V_1:= \partial_2M\ \text{and}
\ W_1:=\partial_3M.$$
In view of Proposition 2.3 and (3.1$'$) we have

$$
\leqalignno{
\log T_{an}(M_2,\partial_-M_2,
g_2, \Cal F_2)&=
-\log T_{an}(\partial_3M, g_3, \Cal F_3)
\cr + \frac 1 2 \chi(\partial_3M&;\Cal F_3)
\cdot (\log 2 + \log \epsilon). &(3.23)\cr}$$
Denote the system obtained by partial gluing by
$(N_{\epsilon},\partial_-N_{\epsilon},\partial_+N_{\epsilon}
,\Cal F_{\epsilon})$ and the resulting metric by $g_\epsilon.$
Then, by Theorem 3.2' and (3.23),

$$
\leqalignno{
\log T_{an}&(N_{\epsilon},\partial_-N_{\epsilon},
g_{\epsilon},\Cal F_{\epsilon})=\log T_{an}(M,\partial_1M, g,\Cal F)
+ \log T((\Cal H(\epsilon)) &(3.24)\cr -
\log T_{an}(\partial_3M, & g_3, \Cal F_3)
+\frac 1 2 \chi(\partial_3M,\Cal F_3)
\cdot \log \epsilon
,\cr}$$
where $\Cal H(\epsilon)$ is the cohomology sequence (3.19) for
the partial gluing of
\newline $(M_1,\partial_-M_1,\partial_+M_1, \Cal F_1)$
and $(M_2,\partial_-M_2,\partial_+M_2,\Cal F_2).$
Notice that there exists a commutative diagram connecting
the cohomology
sequence  $\Cal H'$ (cf (3.21))
and the cohomology sequence $\Cal H(\epsilon)$

{\eightpoint
$$
\matrix &\rightarrow &H^{q-1}(\partial_3M, g_3; \Cal F_3)&
\rightarrow &
H^q(M,\partial_1M \cup\partial_3M ,g;\Cal F)
&\rightarrow &H^q(M, \partial_1M ,g;\Cal F) &\cdots\\
&&\qquad\downarrow (r_{2,\epsilon})_q&&\quad\qquad\qquad\downarrow
(r_{\epsilon})_q &&\qquad\downarrow Id\\
& \rightarrow &H^q(M_2,\partial_- M_2,
g_3\oplus g_0;\Cal F_{[0,\epsilon]}) & \rightarrow&
H^q(N_\epsilon,\partial_- N_\epsilon ,g_\epsilon;
\Cal F_\epsilon)
&\rightarrow &
H^q(M,\partial_1 M, g;
\Cal F) & \cdots.\endmatrix$$ }

In this commutative diagram the vertical arrows are isomorphisms
of $\Cal A-$Hilbert
modules.  The maps $(r_\epsilon)_q$ are induced by a family of
smooth retractions
\newline $r_{\epsilon}:(N_\epsilon,\partial_-N_\epsilon)
\to (M, \partial_1M \cup \partial_3 M)$
with $r_{\epsilon = 0} = id$ which depends smoothly on $\epsilon$
and has the property that the restriction of $r_\epsilon$ to $M\setminus
U_\epsilon$ is the identity where $U_\epsilon$ is a collar neighborhood
of $\partial_3M$ of size $\epsilon.$
The maps $((r_{2,\epsilon})_q),$ $(q \ge 1),$ are isomorphisms

and induced by the maps, assigning
to $\omega \in \Lambda^{q-1}(\partial_3M;\Cal F_3)$ the wedge product $dt
\wedge
\omega.$

One can easily see that
$$\leqalignno{
\lim_{\epsilon \to 0} \log Vol (r_\epsilon)_q =0
&&(3.25)\cr}$$
and
$$\leqalignno{
\log Vol ((r_{2,\epsilon})_q) = \frac 1 2  \log \epsilon\cdot dim_N
H^{q-1}(\partial_3M, g_3; \Cal F_3).&&(3.26)\cr}$$
{}From Proposition 1.3 A(iii), using (3.26), we conclude that
$$
\leqalignno{
- \log T(\Cal H(\epsilon))
=\log T(\Cal H')
+\frac 1 2 \chi(\partial_3M;\Cal F_3)  \log\epsilon
+ \sum (-1)^q \log Vol(r_{\epsilon})_q.
&&(3.27)\cr}
$$
Notice that one can provide
a smooth family of diffeomorphisms $\varphi_\epsilon:M_\epsilon\to M,
\varphi_0=id$ so that
$\varphi_\epsilon^{\#}(g_\epsilon)=\tilde g_\epsilon$ is a smooth
family of Riemannian metrics with $\tilde g_0=g.$ This implies that
$$\lim_{\epsilon\to 0}\log T_{an}(N_\epsilon,\partial_-N_\epsilon,
g_\epsilon,
\Cal F_\epsilon)=
\log T_{an}(M,\partial_1M \cup \partial_3M, g, \Cal F)$$
The result then  follows from (3.24), (3.25) and (3.27)
by passing to the limit $\epsilon\to 0.$ \hfill $\square$
\enddemo

\vfill \eject

\proclaim {Appendix:  On manifolds of determinant class}
\endproclaim

To state our result let us introduce the following notation.
Let $M$ be a compact, connected manifold (possibly with boundary) with
fundamental group $\pi_1(M).$ Let $\Gamma$ be a group and $\theta$ an
arbitrary group homomorphism, $\theta : \pi_1(M) \to \Gamma.$
Denote by $\Cal W_0$ the $(\Cal N(\Gamma), \pi_1(M)^{op})-$Hilbert module
$\ell_2(\Gamma)$ with the
 $\pi_1(M)^{op}-$action induced by $\theta$ and the right
translation of $\Gamma.$ Denote by $\Cal F_{\theta}$ the parallel flat
bundle of $\Cal N(\Gamma)-$Hilbert modules of finite type, induced by
$\Cal W_0.$  Notice that this is the parallel flat bundle induced by the
principal covering $\tilde{M} \to M,$ defined by $\theta.$

\proclaim{Theorem A}
If $\Gamma$ is residually finite then
for any bordism $(M, \partial_-M,\partial_+M),$
the pair $\{(M,\partial_-M,\partial_+M), \Cal F_{\theta}\}$
is of determinant class.
As a consequence, the system $(M,\partial_-M,\partial_+M,\Cal F_\theta)$
is of determinant class.
\endproclaim

Although L\"uck does not use the notion of determinant class, Theorem A is
implicit in his paper [L\"u3].  Unfortunately, there are a number of
misleading  misprints and his definition of the
$L^2$-determinant [L\"u3, Definition, p. 471]
is incorrect.  For the convenience of the reader we present an outline of
L\"uck's arguments to prove Theorem A.

We recall that a group $\Gamma$ is residually finite if
there exists a sequence $(\Gamma_m)_{m \ge 1}$ of nested
normal subgroups of $\Gamma,$
$\Gamma_0 = \Gamma \supseteq \Gamma_1 \supseteq \Gamma_2 \supseteq \dots,$
such that

(i) $\underset {m \ge 1} \to \cap \Gamma_m = \{ e \};$
(ii) $\gamma_m := \#(\Gamma: \Gamma_m) < \infty \,\,(m \ge 1).$

Further recall that
$\big( \Cal C^* (M, \tau, \Cal O_h, \Cal F_{\theta}) , \delta_*)$
is of determinant class if, for $0 \le q \le d,$
$$ - \infty < \int^1_{0^+} \log \lambda d N_q (\lambda),$$
where $N_q (\lambda) = N_q^{comb} (\lambda)$ denotes the $L_2$-combinatorial
spectral distribution
function of the combinatorial q-Laplacian $\Delta_q = \Delta_q^{comb}.$

 Let $M_m$ be the principal $(\Gamma/ \Gamma_m)$- finite cover of $M$
induced by the composition of $\theta$ with the
projection $\Gamma\to \Gamma /\Gamma_m$ and denote by $\Delta_{m; q} \equiv
\Delta_{m; q}^{comb}$ the
combinatorial $q$-Laplacians on $M_m$.  Further let $P_{m; q} ( \lambda)
\equiv
P_{m; q}^{comb} (\lambda)$ be the spectral family
associated to $\Delta_{m; q}$,
defined in such a way that it is right continuous, and
introduce the corresponding  normalized
spectral distribution function

$$
N_{m; q} (\lambda) = \frac 1 {\gamma_m} tr P_{m; q} (\lambda).$$

We point out $N_{m;q}(\lambda)$ are step functions.
Denote by
$det' \Delta_{m; q}$ the modified determinant of $\Delta_{m, q}$, i.e. the
product
of all \underbar {non zero} eigenvalues of $\Delta_{m; q}$.  Let $a_{m; q}$
be
the
smallest non zero eigenvalue and $b_{m; q}$ the largest eigenvalue of
$\Delta_{m; q}$.  Then, for any $a$ and $b$, such that $0 < a < a_{m; q}$ and
$b > b_{m; q}$,

$$\leqalignno{
\frac 1 {\gamma_m} \log det' \Delta_{m; q} = \int_a^b \log \lambda d N_{m; q}
(\lambda).&&(A.1)\cr}$$

Integrating by parts, the Stieltjes integral $\int_a^b \log \lambda d N_{m;
q}
(\lambda)$ can be written as

$$\leqalignno{
\int_a^b \log \lambda d N_{m;q} (\lambda) = &(\log b) \big( N_{m; q} (b)
- N_{m; q} (0) \big)&(A.2)\cr
&- \int_a^b \frac {N_{m; q} (\lambda) - N_{m; q} (0)} \lambda d
\lambda.\cr}$$

Recall that the spectral
distribution function $N_q (\lambda) \equiv N_q^{comb} (\lambda)$
is given by
$$
N_q(\lambda) := tr_{\Cal N(\Gamma)} P_q (\lambda)$$
with $tr_{\Cal N(\Gamma)} P_q (\lambda)$ denoting the von Neumann trace
where $P_q(\lambda)$ is the spectral family
corresponding to $\Delta_q$ defined in such a way
that it is right continuous with respect to $\lambda.$
Notice that $N_q(\lambda)$ is continuous to the right in $\lambda$.  Denote
by $det'\Delta_q$ the modified determinant of $\Delta_q$, given by the
following Stieltjes integral,

$$
\log \, det' \Delta_q = \int^b_{0^+} \log \lambda d N_q (\lambda) :=
\lim_{\epsilon \rightarrow 0^+} \int_\epsilon^b \log \lambda d N_q
(\lambda)$$
with $1 \le b < \infty$ chosen in such a way that
$ |||  \Delta_q |||  < b$ (operator norm).

Notice that $\log \, det' (\Delta_q) \in [-\infty, \infty)$ and recall that
$\Delta_q$ is said to be of determinant class if
$
\int^b_{0^+} \log \lambda d N_q (\lambda) \in (- \infty, \infty).$

Integrating by parts, one obtains

$$\leqalignno{
\qquad\log \, det' (\Delta_q) = & \log b \big( N_q (b) - N_q (0) \big)
&(A.3)\cr
&+ \lim_{\epsilon \rightarrow 0^+} \{(- \log \epsilon) \big( N_q (\epsilon)
-N_q
(0) \big) - \int_\epsilon^b \frac {N_q (\lambda) - N_q (0)} \lambda d \lambda
\}.\cr}$$

Using that $\underset {\epsilon \rightarrow 0^+}\to \liminf (- \log \epsilon)
\big( N_q (\epsilon) - N_q (0) \big) \ge 0$
(in fact, this limit exists and is zero)
and $\frac {N_q (\lambda) - N_q (0)}
\lambda \ge 0$ for $\lambda > 0,$ one sees that

$$\leqalignno{
\log \, det' (\Delta_q) \ge ( \log b) \big(N_q (b) - N_q (0) \big) - \int_0^b
\frac {N_q (\lambda) - N_q (0)} \lambda d \lambda.&&(A.4)\cr}$$

\demo{Proof of Theorem A.2}

The main ingredient is L\"uck's estimate of $\log \, det'(\Delta_q)$ in terms
of
$\log \, det' \Delta_{m, q}$ combined with the fact that $\log \, det'
\Delta_{m, q}\ge 0$ as the determinant $det' \Delta_{m, q}$ is integer
valued.
By
[L\"u, Lemma 2.5], there exists $1 \le b < \infty$ such that, for $m \ge 1$,

$$
|||  \Delta_{m, q} |||  \le b;\quad |||  \Delta_q |||  \le b.$$
\enddemo

By [L\"u, Lemma 3.3.1],

$$\leqalignno{
\int_0^b \frac {N_q (\lambda) - N_q (0)} \lambda d \lambda \le
\liminf_{m \rightarrow \infty} \int_0^b \frac
{N_{m; q} (\lambda) - N_{m; q} (0)} \lambda d \lambda.&&(A.5)\cr}$$

Combining (A.1) and (A.2) with the inequalities $\log \, det' \Delta_{m; q}
\ge
0$, it follows that

$$\leqalignno{
\int_0^b \frac {N_{m; q} (\lambda) - N_{m; q} (0)} \lambda d \lambda
\le (\log b) \big(N_{m; q} (b) - N_{m; q} (0)\big).&&(A.6)\cr}$$

{}From (A.4) - (A.6) we then conclude that

$$\leqalignno{
\log \, det' \Delta_q \ge &(\log b) \big( N_q (b) - N_q (0) \big)&(A.7)\cr
&- \liminf_{m \rightarrow \infty}(\log b) \big( N_{m; q} (b) - N_{m; q}
(0) \big).\cr}$$

The estimate (A.7) together with the identities ([L\"u3, Theorem 2.3.1])

$$
N_q (\lambda) = \lim_{\epsilon \rightarrow 0^+} \liminf_{m \rightarrow
\infty}
N_{m; q} (\lambda + \epsilon)\, \text { ([L\"u3, Theorem 2.3.1])}$$
and

$$
N_q (0) = \lim_{m \rightarrow \infty} N_{m; q} (0)\, \text
{ ([L\"u3, Theorem 2.3.2])}$$
yield that
$
\log \, det' \Delta_q \ge 0$
which proves that $\Delta_q$ is of determinant class.
Since $q$ is arbitrary this
shows that the pair
$\{(M, \partial M_-, \partial M_+),\Cal F_{\theta}\}$
is of determinant class.

\hfill $\square$

\vfill \eject

\end

\enddocument

\enddocument